\documentclass[11pt]{article}
\usepackage{amssymb,amsmath,color,amsfonts,float,amscd,amsthm,bm,mathrsfs} 

\allowdisplaybreaks[4]

\catcode`!=11
\let\!int\int \def\int{\displaystyle\!int}
\let\!lim\lim \def\lim{\displaystyle\!lim}
\let\!sum\sum \def\sum{\displaystyle\!sum}
\let\!sup\sup \def\sup{\displaystyle\!sup}
\let\!inf\inf \def\inf{\displaystyle\!inf}
\let\!cap\cap \def\cap{\displaystyle\!cap}
\let\!max\max \def\max{\displaystyle\!max}
\let\!min\min \def\min{\displaystyle\!min}

\let\oldsection\section
\renewcommand\section{\setcounter{equation}{0}\oldsection}

\newtheorem{theorem}{Theorem}[section]
\newtheorem{lemma}[theorem]{Lemma}

\theoremstyle{definition}

\theoremstyle{remark}

\begin{document}

\title{Uniform regularity estimates and invisicid limit for the compressible non-resistive magnetohydrodynamics system}

\author{Xiufang Cui\footnote{School of Mathematical Sciences, Shanghai Jiao Tong University, Shanghai 200240, P. R. China. Email: cuixiufang@sjtu.edu.cn}
\and
Shengxin Li\footnote{School of Mathematical Sciences, Shanghai Jiao Tong University, Shanghai 200240, P. R. China. Email: lishengxin@sjtu.edu.cn}
\and
Feng Xie
\footnote{School of Mathematical Sciences, Shanghai Jiao Tong University, Shanghai 200240, P. R. China.
Email: tzxief@sjtu.edu.cn}
}
\date{}
\maketitle

\begin{abstract}
We are concerned with the uniform regularity estimates of solutions to the two dimensional compressible non-resistive magnetohydrodynamics (MHD) equations with the no-slip boundary condition on velocity in the half plane. Under the assumption that the initial magnetic field is transverse to the boundary, the uniform conormal energy estimates are established for the solutions to compressible MHD equations with respect to small viscosity coefficients. As a direct consequence, we proved the inviscid limit of solutions from viscous MHD systems to the ideal MHD systems in $L^\infty$ sense. It shows that the transverse magnetic field can prevent the boundary layers from occurring in some physical regime.
\end{abstract}

%
%
\maketitle





	\section{Introduction}
The system of compressible magnetohydrodynamics (MHD) equations is one of the most fundamental models in magneto-fluid mechanics, which  describes the dynamics of electrically conducting fluids arising from plasmas and other related physical phenomena.
In this paper, we  are concerned with the inviscid limit of the solutions to the
two dimensional compressible non-resistive magnetohydrodynamics equations:
\begin{align}\label{1.1}
\begin{cases}
\partial_t\rho^\varepsilon+\nabla\cdot(\rho^\varepsilon {\bf{v}}^\varepsilon)=0,\\
\partial_t(\rho^\varepsilon {\bf{v}}^\varepsilon)+\mathrm{div}(\rho^\varepsilon {\bf{v}}^\varepsilon\otimes {\bf{v}}^\varepsilon)-\varepsilon\mu\Delta {\bf{v}}^\varepsilon-\varepsilon(\mu+\lambda)\nabla(\nabla\cdot {\bf{v}}^\varepsilon)+\nabla p ^{\varepsilon}=(\nabla\times {\bf{B}}^\varepsilon)\times {\bf{B}}^\varepsilon,\\
\partial_t{\bf{B}}^\varepsilon-\nabla\times({\bf{v}}^\varepsilon\times {\bf{B}}^\varepsilon)=0,\quad \mathrm{div}\  {\bf{B}}^\varepsilon=0,
\end{cases}
\end{align}	
where $\rho^\varepsilon $ is the density of the fluid, ${\bf{v}}^\varepsilon=(v_1^\varepsilon, v_2^\varepsilon )$ is the fluid velocity,  ${\bf{B}}^\varepsilon=(b_1^\varepsilon, b_2^\varepsilon )$ is the magnetic field.  The  shear viscosity coefficient  $\varepsilon \mu$ and  the bulk viscosity coefficient $\varepsilon\lambda$  satisfy $\mu> 0$ and $\mu+\lambda>0$ with $\varepsilon\in (0, 1)$ being a small parameter. Here $\nabla=(\partial_x, \partial_y)$ and $\Delta=\partial_x^2+\partial_y^2$. The pressure $p^\varepsilon$ is the function of $\rho^\varepsilon$,  which takes the following form:
\begin{align}\label{1.2}
p^\varepsilon=(\rho^\varepsilon)^\gamma,\quad \gamma\geq1,
\end{align}
where $\gamma$ is some constant.
The initial data is given by
\begin{align*}
(\rho^\varepsilon, {\bf{v}}^\varepsilon, {\bf{B}}^\varepsilon)(t, x, y)|_{t=0}=(\rho^\varepsilon_0, {\bf{v}}^\varepsilon_0,  {\bf{B}}^\varepsilon_0)(x, y).
\end{align*}	

In the presence of boundary, we impose  the no-slip boundary condition on the velocity field:
\begin{align}\label{1.3}
{\bf{v}}^\varepsilon|_{y=0}=0.
\end{align}	
Based on the  above no-slip boundary condition, we do not need to impose any boundary condition for the magnetic field. 	
In fact, since the equation of $b_2^\varepsilon$ is
\begin{align*}
\partial_t b_2^\varepsilon+b_2^\varepsilon\partial_x v_1^\varepsilon -b^\varepsilon_1\partial_x v_2^\varepsilon -v_2^\varepsilon\partial_x b_1^\varepsilon+v_1^\varepsilon\partial_x b_2^\varepsilon=0,
\end{align*}
and ${\bf{v}}^\varepsilon|_{y=0}=0$, it is direct to check that
$$ b_2^\varepsilon(t, x, y)=b_2^\varepsilon(0, x, y), \quad \mathrm{on}\quad y=0.$$
Hence,  to simplify the calculation,  one sets $$b_2^\varepsilon(0, x, y)|_{y=0}=1,\quad  \tilde b_1^\varepsilon(t, x, y)=b_1^\varepsilon(t, x, y), \quad\tilde b_2^\varepsilon(t, x, y)=b_2^\varepsilon(t, x, y)-1, $$
then  $$\tilde b_2^\varepsilon(t,x, y)|_{y=0}=0.$$ 	

Before stating the main results in this paper, it is helpful to review some related works. As an important model in physics and because of its mathematical challenges, the system of the MHD equations attracts many attentions from physicists and mathematicians. Readers can refer to \cite{CW02, CW03, DF06, HT05, Wang03}  for the study of compressible MHD equations and \cite{AZ17, CRW13, CW11, LZ14, LXZ15} for incompressible MHD equations. And the inviscid limit problem is a challenging but important problem in hydrodynamics, which is one of the hot topics in mechanics and applied mathematics. For example, \cite{AD04, CP10, CPW95, HWY12, RWX14, RXZ16}. Moreover, when the inviscid limit process is considered in a domain with physical boundaries, the mathematical analysis will become much more difficult due to the possible presence of strong  boundary layers. As is known that the strong boundary layer will disappear, provided that the velocity is imposed the so-called Navier-slip boundary condition. Thus, under this kind of slip boundary condition, it is reasonable to establish the inviscid limit for both compressible Navier-Stokes equations \cite{WW12} and incompressible Navier-Stokes equations  \cite{IS11, MR12} without specifying the boundary layer corrector. And we also refer to \cite{DXX20, WZ17,  XXW09} for the inviscid limit of incompressible MHD system with Navier-slip boundary conditions.
However, when the velocity is given the no-slip boundary condition, in general the strong boundary layer must occur \cite{GGW,LWXY,O,Ole,prandtl,WXY}. Consequently, the inviscid limit in $L^\infty$ sense becomes dramatically difficult due to the uncontrollability of the vorticity of strong boundary layer as the coefficient of viscosity goes to zero.
Recently, in \cite{LXY191, LXY192}, Liu, Xie and Yang established the well-posedness of solutions to MHD boundary layer equations and proved the validity of Prandtl boudary layer expansion in the Sobolev spaces under the conditions that the tangential component of magnetic filed does not degenerate near the physical boundary initially.
However, it is noted that there are very few results about the inviscid limit of solutions to
both compressible Navier-Stokes equations and compressible MHD equations with the no-slip boundary condition. This is one of key motivations of this paper.
The main goal of this paper is to prove the inviscid limit from compressible MHD equations \eqref{1.1} to the following ideal MHD equations in finite regularity functional spaces as $\varepsilon\rightarrow 0$ under the no-slip boundary condition \eqref{1.3}.	
\begin{align} \label{1.4}
\begin{cases}
\partial_t\rho^0+\nabla\cdot(\rho^0 {\bf{v}}^0)=0,\\
\partial_t(\rho^0 {\bf{v}}^0)+\mathrm{div}(\rho^0 {\bf{v}}^0\otimes {\bf{v}}^0) +\nabla p^0=(\nabla\times {\bf{B}}^0)\times {\bf{B}}^0,\\
\partial_t{\bf{B}}^0-\nabla\times({\bf{v}}^0\times {\bf{B}}^0)=0,\quad \mathrm{div}\ {\bf{B}}^0=0.
\end{cases}
\end{align}	

It is believed that  the suitable magnetic field has a stabilizing effect on the boundary layer  which could prevent some kind of instability or singularity being occurring.  Hence,  it is possible to justify inviscid limit for the compressible MHD equations under the no-slip boundary condition provided some suitable magnetic field is given initially. Precisely, we will verify the inviscid limit from viscous MHD equations (\ref{1.1}) to ideal MHD case (\ref{1.4}) by establishing the uniform energy estimates of solutions to (\ref{1.1}) with the no-slip boundary condition (\ref{1.3}) when the normal component of the initial magnetic field does not degenerate. It is pointed out that the similar problem was considered in \cite{LXY20} for the incompressible MHD equations.
There are two novelties of the methods used in this paper. First, we have to derive the estimates of the normal derivatives for all components of both velocity and magnetic field due to the lack of divergence free condition of velocity. Second, the arguments in deriving the estimates of normal derivatives for pressure are totally different from those used in both \cite{LXY20} and \cite{MR12}.

To formulate the main results, it is necessary to introduce some notations. Throughout this paper,  we take the following conormal derivatives of functions depending on $(t,\bf{x})$:
\begin{align*}
\mathcal{Z}_0=\partial_t,\quad \mathcal{Z}_1=\partial_x,\quad \mathcal{Z}_2=\phi(y)\partial_y, \quad \mathcal{Z}^\alpha=\mathcal{Z}_0^{\alpha_0}\mathcal{Z}_1^{\alpha_1}\mathcal{Z}_2^{\alpha_2},
\end{align*}
where  the spatial variables ${\bf{x}}=(x,  y)$, the multi-index $\alpha=(\alpha_0, \alpha_1, \alpha_2)$ and  $|\alpha|=\alpha_0+\alpha_1+\alpha_2$.  The weight $\phi(y)$  is a smooth and bounded function of $y $ satisfying $\phi(y)|_{y=0}=0$ and $\phi'|_{y=0}>0$, typically, one can choose $\phi(y)=\frac{y}{1+y}$.

For any integer $m\in \mathbb{N}$, we denote the conormal Sobolev space
\begin{align*}
H_{co}^m([0, T] \times \mathbb{R}_+^2)=\{f(t, {\bf{x}}): \mathcal{Z}^\alpha f\in L^2([0, T]\times \mathbb{R}_+^2), \quad |\alpha|\leq m\}.
\end{align*}

For any $t\geq 0$, we set the norms
\begin{align*}
\|f(t)\|_{m}^2=\sum_{|\alpha|\leq m}\|\mathcal{Z}^\alpha f(t,\cdot)\|^2,
\end{align*}
then
\begin{align*}
\|f\|_{H_{co}^m}^2=\int_0^t\|f(s)\|_{m}^2 ds.
\end{align*}

Similarly,
\begin{align*}
W_{co}^{m, \infty}([0, T] \times \mathbb{R}_+^2)=\{f(t, {\bf{x}}): \mathcal{Z}^\alpha f\in L^\infty([0, T]\times \mathbb{R}_+^2), \quad |\alpha|\leq m\}
\end{align*}
with the norm
\begin{align*}
\|f\|_{m, \infty}=\sum_{|\alpha|\leq m}\|\mathcal{Z}^\alpha f \|_{L_{t, {\bf{x}}}^\infty}.
\end{align*}

It is also convenient to introduce the function setting
\begin{align*}
\Lambda^m(t)=\{(\rho, {\bf{v}}, {\bf B}): \partial_y^i(\rho-1, {\bf{v}}, {\bf B}-\overset{\rightarrow}{e_y})\in H_{co}^{m-i},  i=0, 1, 2 \}
\end{align*}
with $\overset{\rightarrow}{e_y}=(0, 1)$.

To derive the uniform conormal estimates of the classical solution  $(\rho^\varepsilon, {\bf{v}}^\varepsilon, {\bf{B}}^\varepsilon)$ to compressible MHD \eqref{1.1}-\eqref{1.3}, the  following energy norm is needed:
\begin{align*}
N_m(t)=&\sum_{|\alpha|\leq m}\sup_{0\leq s\leq t}\int_{\mathbb{R}_+^2}\left(\rho^\varepsilon(s)|\mathcal{Z}^\alpha {\bf{v}}^\varepsilon(s)|^2+|\mathcal{Z}^\alpha ({\bf{B}}^\varepsilon-\overset{\rightarrow}{e_y})(s)|^2+\gamma^{-1}(p^\varepsilon)^{-1}|\mathcal{Z}^\alpha (p^\varepsilon(s)-1)|^2\right) d{\bf{x}} \\
&+ \|\partial_y({\bf{v}}^\varepsilon, b_1^\varepsilon, p^\varepsilon) \|_{H_{co}^{m-1}}^2
+ \|\partial_y^2 ({\bf{v}}^\varepsilon, b_1^\varepsilon, p^\varepsilon)\|_{H_{co}^{m-2}}^2+\varepsilon\mu \|\nabla {\bf{v}}^\varepsilon\|_{H_{co}^m}^2   +\varepsilon(\mu+\lambda) \|\nabla\cdot {\bf{v}}^\varepsilon\|_{H_{co}^m}^2\\
&+\varepsilon^2\mu^2\left(\|\partial_y^2 v_1^\varepsilon\|_{H_{co}^{m-1}}^2+ \|\partial_y^3 v_1^\varepsilon\|_{H_{co}^{m-2}}^2\right) +\varepsilon^2(2\mu+\lambda)^2\left(\|\partial_y^2 v_2^\varepsilon\|_{H_{co}^{m-1}}^2 + \|\partial_y^3 v_2^\varepsilon\|_{H_{co}^{m-2}}^2\right).
\end{align*}
Now, it is position to state the main results of this paper.
\begin{theorem}\label{Th1}
(Uniform regularity and inviscid limit) Let  the integer $m\geq 9$. Suppose the initial data $ (\rho^\varepsilon_0, \bf{v}^\varepsilon_0, \bf{B}^\varepsilon_0)$ satisfies
	\begin{align}\label{1.5}
	\sum_{i=0}^2\|\partial_y^i(\rho_0^\varepsilon-1, {\bf{v}}_0^\varepsilon, {\bf{B}}_0^\varepsilon-\overset{\rightarrow}{e_y})\|_{m-i}^2 \leq\sigma,
	\end{align}
	where $\sigma>0$ is some sufficiently small  constant. Then for the classical solution $U^\varepsilon=(\rho^\varepsilon, {\bf{v}}^\varepsilon, {\bf{B}}^\varepsilon)\in \Lambda^m(T)$ to  compressible MHD  \eqref{1.1}-\eqref{1.3}, there exists  a time $T>0$  indepent  of $\varepsilon$ such that for any   $t\in [0, T]$,  the following regularity estimate holds:
	\begin{align*}
	N_m(t)+ &\sum_{|\alpha|+i\leq m\atop i=1,2}\int_{\mathbb{R}_+^2}
	\left(\gamma^{-1}\varepsilon(2\mu+\lambda) (p^\varepsilon)^{-1}(t)|\mathcal{Z}^\alpha \partial_y^i (p^\varepsilon(t)-1)|^2 +\varepsilon\mu |\mathcal{Z}^\alpha\partial_y^i b_1^\varepsilon(t)|^2 \right)d{\bf{x}}
	\le C\sigma
	\end{align*}
	where $C>0$ is  some   constant, which is independent of $\varepsilon$.

	Moreover, there exists a  unique solution   $U^0=(\rho^0, {\bf{v}}^0, {\bf{B}}^0)\in \Lambda^m(T)$  to the ideal compressible MHD equations \eqref{1.4}, such that
	\begin{align*}
	\lim_{\varepsilon\rightarrow 0}\sup_{t\in[0, T]}\|(U^{\varepsilon}-U^0, \partial_y(U^{\varepsilon}-U^0)(t, \cdot)\|_{L^\infty(\mathbb{R}_+^2)}=0.
	\end{align*}
	
\end{theorem}

Below, it is helpful to introduce the strategy of proof of Theorem \ref{Th1}. We  first derive the   uniform conormal estimates  for  the classical solution $(\rho^\varepsilon, {\bf{v}}^\varepsilon, {\bf{B}}^\varepsilon)$ to compressible MHD equations \eqref{1.1}-\eqref{1.3} in Section 3. To close the energy estimates in Section 3,  by the anistropic Sobolev embedding inequality in the conormal Sobolev space,  it suffices to give the conormal estimates for both the first order and the second order normal derivatives of $(\rho^\varepsilon, {\bf{v}}^\varepsilon, {\bf{B}}^\varepsilon)$, which are given in details in Section 4. It is remarked that the equations of magnetic field and mass play an essential role in  presenting  the normal derivatives of velocity field.  It is indeed the key point for us to obtain the uniform conormal estimates for the  normal derivatives of $(\rho^\varepsilon, {\bf{v}}^\varepsilon, {\bf{B}}^\varepsilon)$.
Based on these energy estimates achieved in Section 3 and Section 4, we prove that the solutions to compressible MHD equations \eqref{1.1}-\eqref{1.3} are uniformly bounded in the conormal Sobolev space in a fixed time interval, which is independent of the small viscosity coefficient. As a direct consequence,  we justify the inviscid limit of solutions to the compressible MHD equations \eqref{1.1}-\eqref{1.3} as $\varepsilon\rightarrow 0$ by some compactness arguments.

The paper is organized as follows.  In Section 2, we give  some preliminaries.  The uniform conormal estimates of solutions to the compressible MHD equations \eqref{1.1}-\eqref{1.3} are established in Section 3. We derive the uniform conormal estimates for normal derivatives of the classical solution to \eqref{1.1}-\eqref{1.3}  in Section 4. In Section 5, we  prove Theorem \ref{Th1} based on the results in both Section 3 and Section 4.

The notation $A\lesssim B$ stands for $A\leq C B$ for some generic constant $C>0$ independent of $\varepsilon$. And we denote the polynomial functions by $\mathcal{P}(\cdot)$, which may vary from line to line,  and the commutator is expressed by $[\cdot, \cdot]$.

\section{Preliminaries}
In this section,  we present some properties of the conormal Sobolev spaces, which will be used frequently in the rest of this paper.
The first one is the Sobolev-Gagliardo-Nirenberg-Moser type inequality for the conormal Sobolev space and its proof can be found in \cite{G90}.
\begin{lemma}\label{lem 2.2}
	For the functions $f, g \in L^\infty([0,T]\times\mathbb{R}_+^2)\cap H_{co}^m([0,T]\times\mathbb{R}_+^2)$ with $m\in \mathbb{N}$, it holds that for any  $\alpha, \beta\in \mathbb{N}^3$ with $|\alpha|+|\beta|=m$,
	\begin{align}\label{2.1}
	\int_0^t\|(\mathcal{Z}^\alpha f\mathcal{Z}^\beta g) (s)\|^2 ds\lesssim \|f\|_{L_{t,\bf{x}}^\infty}^2\int_0^t\|g(s)\|_{m}^2 ds+\|g\|_{L_{t,\bf{x}}^\infty}^2\int_0^t\|f(s)\|_{m}^2 ds.
	\end{align}
\end{lemma}
Next, we give the anisotropic Sobolev embedding property in the conormal  Sobolev space, whose proof can be found in \cite{P16}.
\begin{lemma}
	Suppose that  $f(t,{\bf{x}})\in H_{co}^3([0, t]\times \mathbb{R}^2_+)$ and $\partial_y f\in H_{co}^2([0, t]\times \mathbb{R}^2_+)$, then
	\begin{align*}
	\|f\|_{L_{t, {\bf{x}}}^\infty}^2\lesssim\|f(0)\|_{2}^2+\|\partial_yf(0)\|_{1}^2+\int_0^t\left(\|f(s)\|_{3}^2+\|\partial_yf(s)\|_{2}^2\right)ds.
	\end{align*}
\end{lemma}
To deal with the commutator involving conormal derivatives, the following properties of commutators will be frequently used later, and their proof can be found in \cite{P16}.	

Notice that
\begin{align*}
[\mathcal{Z}^\alpha, \partial_t]=	[\mathcal{Z}^\alpha, \partial_x]=0.
\end{align*}
However, $\mathcal{Z}_2$ does not commute with $\partial_y$. A direct calculation by induction shows that for any integer $m\geq 1$, there exist two families of bounded smooth functions $\{\phi_{k,m}(y)\}_{0\leq k\leq m-1}$ and $\{\phi^{k,m}(y)\}_{0\leq k\leq m-1}$ depending only on $\phi(y)$, such that
\begin{align}\label{2.2}
[\mathcal{Z}_2^m, \partial_y]=\sum_{k=0}^{m-1}\phi_{k,m}(y)\mathcal{Z}_y^k\partial_y=\sum_{k=0}^{m-1}\phi^{k,m}(y)\partial_y\mathcal{Z}_y^k.
\end{align}
Similarly, there  exist other two families of bounded smooth functions $\{(\phi_{1,k,m}(y), \phi_{2,k,m}(y))\}_{0\leq k\leq m-1}$
and $\{(\phi^{1,k,m}(y), \phi^{2,k,m}(y))\}_{0\leq k\leq m-1}$ depending only on $\phi(y)$, such that
\begin{align*}
[\mathcal{Z}_2^m, \partial_y^2]=&\sum_{k=0}^{m-1}\left(\phi_{1,k,m}(y)\mathcal{Z}_2^k\partial_y+\phi_{1,k,m}(y)\mathcal{Z}_2^k\partial_y^2\right)\\
=&\sum_{k=0}^{m-1}\left(\phi^{1,k,m}(y)\partial_y\mathcal{Z}_2^k+\phi^{2,k,m}(y)\partial_y^2\mathcal{Z}_2^k\right).
\end{align*}
In addition, for a suitable function $f$, there exist two  families of bounded smooth functions $\{\varphi_{ k, m}(y)\}_{0\leq k\leq m-1}$ and   $\{\varphi^{ k, m}(y)\}_{0\leq k\leq m-1}$  depending only on $\phi(y)$, such that
\begin{align}\label{2.3}
[\mathcal{Z}_2^m, \phi^{-1}]f=\sum_{k=0}^{m-1}\varphi_{k, m}(y)\mathcal{Z}_y^k(\phi^{-1}f),
\end{align}
and
\begin{align}\label{2.4}
[\mathcal{Z}_2^m, \phi ]f=\sum_{k=0}^{m-1}\varphi^{ k, m}(y)\mathcal{Z}_y^k(\phi f).
\end{align}

\section{Conormal Energy Estimates}
This section is devoted to the conormal  energy estimates of the classical solution $(\rho^\varepsilon, {\bf{v}}^{\varepsilon}, {\bf{B}}^\varepsilon)$ to  compressible  MHD   \eqref{1.1}-\eqref{1.3}.
The equations of \eqref{1.1}  can be rewritten as follows
\begin{align}\label{3.1}
\begin{cases}
\partial_t \rho+{\bf{v}}\cdot\nabla\rho +\rho\nabla\cdot {\bf{v}}=0,\\
\rho\partial_t {\bf{v}}+\rho {\bf{v}}\cdot\nabla {\bf{v}}-\varepsilon \mu \Delta  {\bf{v}}-\varepsilon(\mu+\lambda)\nabla(\nabla\cdot {\bf{v}})+ \nabla p =(\nabla\times {\bf{B}})\times {\bf{B}},\\
\partial_t{\bf{B}}-\nabla\times({\bf{v}}\times {\bf{B}}) =0, \quad \mathrm{div}\ {\bf{B}}=0.
\end{cases}
\end{align}
The main result in this section is listed in the following lemma:
\begin{lemma}\label{lem1}
	Under the assumption in Theorem \ref{Th1}, for any integer $m>0$ and $|\alpha|\leq m$, the classical solution $(\rho, {\bf{v}}, {\bf{B}})$  to the compressible MHD equations \eqref{3.1} with no-slip boundary condition \eqref{1.3}  satisfies
	\begin{align*}
	&\int_{\mathbb{R}_+^2}\left(\rho(t)|\mathcal{Z}^\alpha {\bf{v}}(t)|^2+|\mathcal{Z}^\alpha ({\bf{B}}(t)-\overset{\rightarrow}{e_y})|^2  +  \gamma^{-1} p^{-1}(t)|\mathcal{Z}^\alpha (p(t)-1)|^2\right) d{\bf{x}}\\
	&+\varepsilon\mu\int_0^t\|\nabla {\bf{v}}\|_{m}^2 ds+\varepsilon(\mu+\lambda)\int_0^t \|\nabla \cdot {\bf{v}}\|_{m}^2 ds\\
	\lesssim & \int_{\mathbb{R}_+^2}\left(\rho_0|\mathcal{Z}^\alpha {\bf{v}}_0|^2+|\mathcal{Z}^\alpha ({\bf{B}}_0-\overset{\rightarrow}{e_y})|^2  +  \gamma^{-1} p_0^{-1}|\mathcal{Z}^\alpha (p_0-1)|^2\right) d{\bf{x}}\\
	&+\left(1+ \|({\bf{v}}, {\bf{B}}, p)\|_{2, \infty}^2+\|\partial_y ({\bf{v}}, {\bf{B}},   p)\|_{1, \infty}^2 \right) \int_0^t\left(\|({\bf{v}}, {\bf{B}}-\overset{\rightarrow}{e_y}, p-1)(s)\|_{m}^2+\|\partial_y( {\bf{v}},   {\bf{B}},  p)(s)\|_{m-1}^2\right) ds.
	\end{align*}
\end{lemma}
\begin{proof}
	For any multi-index $\alpha$ satisfying $|\alpha|\leq m$,  by applying the conormal derivatives $\mathcal{Z}^\alpha$ to the last two equations in \eqref{3.1}, then multiplying $(\mathcal{Z}^\alpha {\bf{v}}, \mathcal{Z}^\alpha ({\bf{B}}-\overset{\rightarrow}{e_y}))$ on both sides of the resulting equalities and integrating it over $[0, t]\times\mathbb{R}_+^2$, we get
	\begin{equation}\label{3.2}
	\begin{split}
	&\frac12	 \int_{\mathbb{R}_+^2}\left(\rho(t)|\mathcal{Z}^\alpha {\bf{v}}(t)|^2  +	 |\mathcal{Z}^\alpha ({\bf{B}}(t)-\overset{\rightarrow}{e_y})|^2 \right)d{\bf{x}}-\frac12	 \int_{\mathbb{R}_+^2}\left(\rho_0|\mathcal{Z}^\alpha {\bf{v}}_0|^2  +	 |\mathcal{Z}^\alpha ({\bf{B}}_0-\overset{\rightarrow}{e_y})|^2 \right)d{\bf{x}}  \\
	=&  -\int_0^t\int_{\mathbb{R}_+^2} \mathcal{Z}^\alpha\nabla p\cdot \mathcal{Z}^\alpha {\bf{v}} \ d{\bf{x}}ds+\varepsilon \mu\int_0^t\int_{\mathbb{R}_+^2}  \mathcal{Z}^\alpha \Delta {\bf{v}} \cdot \mathcal{Z}^\alpha {\bf{v}} \ d{\bf{x}}ds\\
	& +\varepsilon(\mu+\lambda)\int_0^t\int_{\mathbb{R}_+^2} \mathcal{Z}^\alpha\nabla(\nabla\cdot {\bf{v}})\cdot \mathcal{Z}^\alpha {\bf{v}} \ d{\bf{x}}ds +\int_0^t\int_{\mathbb{R}_+^2} \left(\mathcal{C}_1^\alpha+\mathcal{C}_2^\alpha\right)\cdot \mathcal{Z}^\alpha {\bf{v}}\  d{\bf{x}}ds\\
	& +\int_0^t\int_{\mathbb{R}_+^2} \mathcal{Z}^\alpha [(\nabla\times {\bf{B}})\times {\bf{B}}]\cdot \mathcal{Z}^\alpha {\bf{v}} \ d{\bf{x}}ds+\int_0^t\int_{\mathbb{R}_+^2} \mathcal{Z}^\alpha[\nabla\times({\bf{v}}\times {\bf{B}})]\cdot\mathcal{Z}^\alpha ({\bf{B}}-\overset{\rightarrow}{e_y})\  d{\bf{x}}ds,
	\end{split}
	\end{equation}
	where
	\begin{align*}
	&	\mathcal{C}_1^\alpha=-[Z^\alpha, \rho\partial_t]{\bf v} =-\sum_{\beta+\gamma=\alpha\atop |\beta|\geq 1} C_\alpha^\beta \mathcal{Z}^\beta \rho\mathcal{Z}^\gamma\partial_t {\bf{v}},
	\end{align*}
	and 	
	\begin{align*}
	\mathcal{C}_2^\alpha =-[Z^\alpha, \rho{\bf v}\cdot\nabla]{\bf v}=-\sum_{\beta+\gamma=\alpha\atop |\beta|\geq 1}C_\alpha^\beta \mathcal{Z}^\beta(\rho {\bf{v}})\cdot\mathcal{Z}^\gamma \nabla {\bf{v}}-\rho {\bf{v}}\cdot[\mathcal{Z}^\alpha, \nabla]{\bf{v}}.
	\end{align*}
	For the first  term on the right hand side of \eqref{3.2}, from integration by parts, we have
	\begin{equation}\label{3.3}
	\begin{split}
	& -\int_0^t\int_{\mathbb{R}_+^2} \mathcal{Z}^\alpha\nabla p\cdot \mathcal{Z}^\alpha {\bf{v}}\  d{\bf{x}}ds\\ =&-\int_0^t\int_{\mathbb{R}_+^2}\nabla \mathcal{Z}^\alpha (p-1)\cdot \mathcal{Z}^\alpha {\bf{v}}\ d{\bf{x}}ds -\int_0^t\int_{\mathbb{R}_+^2}[ \mathcal{Z}^\alpha,  \nabla] (p-1)\cdot \mathcal{Z}^\alpha {\bf{v}}\  d{\bf{x}}ds\\
	=&\int_0^t\int_{\mathbb{R}_+^2}  \mathcal{Z}^\alpha (p-1)\cdot \mathcal{Z}^\alpha(\nabla\cdot {\bf{v}}) \ d{\bf{x}}ds+\int_0^t\int_{\mathbb{R}_+^2} \mathcal{Z}^\alpha (p-1)\cdot [\nabla\cdot, \mathcal{Z}^\alpha] {\bf{v}} \ d{\bf{x}}ds\\
	& - \int_0^t\int_{\mathbb{R}_+^2}[ \mathcal{Z}^\alpha,  \nabla] (p-1)\cdot \mathcal{Z}^\alpha {\bf{v}}\ d{\bf{x}}ds.
	\end{split}
	\end{equation}
	Since
	\begin{align}\label{3.4}
	\nabla\cdot  {\bf{v}}= - \gamma^{-1} p^{-1}\partial_t (p-1) - \gamma^{-1} p^{-1}{\bf{v}}\cdot\nabla (p-1),
	\end{align}	
	by inserting  (\ref{3.4}) into the first term on the right hand side of \eqref{3.3},  one obtains
	\begin{equation}\label{3.5}
	\begin{split}
	&\int_0^t\int_{\mathbb{R}_+^2} \mathcal{Z}^\alpha (p-1)\cdot \mathcal{Z}^\alpha (\nabla\cdot {\bf{v}}) \ d{\bf{x}}ds\\
	=&-\gamma^{-1}\int_0^t\int_{\mathbb{R}_+^2} \mathcal{Z}^\alpha (p-1)\cdot ( p^{-1}\partial_t	\mathcal{Z}^\alpha   (p-1))\ d{\bf{x}}ds\\
	&-\gamma^{-1}\int_0^t\int_{\mathbb{R}_+^2} \mathcal{Z}^\alpha (p-1)\cdot ( p^{-1}{\bf{v}}\cdot	\mathcal{Z}^\alpha  \nabla (p-1))\ d{\bf{x}}ds\\
	& -\gamma^{-1}\sum_{\beta+\gamma=\alpha\atop |\beta|\geq1} C_\alpha^\beta\int_0^t\int_{\mathbb{R}_+^2} \mathcal{Z}^\alpha (p-1)\cdot \mathcal{Z}^\beta   p^{-1}\mathcal{Z}^\gamma  \partial_t (p-1)\ d{\bf{x}}ds\\ &-\gamma^{-1}\sum_{\beta+\gamma=\alpha\atop |\beta|\geq1} C_\alpha^\beta\int_0^t\int_{\mathbb{R}_+^2} \mathcal{Z}^\alpha (p-1) \cdot \mathcal{Z}^\beta  (  p^{-1}{\bf{v}})\cdot\mathcal{Z}^\gamma \nabla (p-1)\ d{\bf{x}}ds.
	\end{split}
	\end{equation}
	Next, we handle term by term on the right hand side of \eqref{3.5}. First,
	\begin{align*}
	&-\gamma^{-1}\int_0^t\int_{\mathbb{R}_+^2} \mathcal{Z}^\alpha (p-1)\cdot (p^{-1} \partial_t	\mathcal{Z}^\alpha  (p-1)) \  d{\bf{x}}ds\\
	=&-\frac{1}{2\gamma}\frac{d}{dt}\int_0^t\int_{\mathbb{R}_+^2}    p^{-1}|\mathcal{Z}^\alpha (p-1)|^2\ d{\bf{x}}ds +\frac{1}{2\gamma}\int_0^t\int_{\mathbb{R}_+^2} \partial_t  p^{-1}|\mathcal{Z}^\alpha (p-1)|^2\ d{\bf{x}}ds\\
	\lesssim &-\frac{1}{2\gamma}\int_{\mathbb{R}_+^2} p^{-1}(t)|\mathcal{Z}^\alpha (p(t)-1)|^2 d{\bf{x}}+\frac{1}{2\gamma}\int_{\mathbb{R}_+^2}  p_0^{-1}|\mathcal{Z}^\alpha (p_0-1)|^2 d{\bf{x}}+  \|  p^{-1}\|_{1, \infty}\int_0^t \|p-1\|_{m}^2 ds.
	\end{align*}
	Second,
	\begin{equation}\label{3.6}
	\begin{split}
	&-\gamma^{-1}\int_0^t\int_{\mathbb{R}_+^2} \mathcal{Z}^\alpha (p-1)\cdot  (  p^{-1}{\bf{v}}\cdot	\mathcal{Z}^\alpha  \nabla (p-1)  ) \ d{\bf{x}}ds\\
	=&	-\gamma^{-1}\int_0^t\int_{\mathbb{R}_+^2}  \mathcal{Z}^\alpha (p-1)\cdot (  p^{-1}{\bf{v}}\cdot[\mathcal{Z}^\alpha , \nabla] (p-1)  )\  d{\bf{x}}ds\\
	&-\gamma^{-1}\int_0^t\int_{\mathbb{R}_+^2} \mathcal{Z}^\alpha (p-1)\cdot(  p^{-1}{\bf{v}}\cdot\nabla\mathcal{Z}^\alpha  (p-1))\ d{\bf{x}}ds.
	\end{split}
	\end{equation}
	For the first term on the right hand side of \eqref{3.6}, by \eqref{2.2} and the a prior assumption that $p^{-1}$ has positive lower bound, we have
	\begin{align*}
	&-\gamma^{-1}\int_0^t\int_{\mathbb{R}_+^2}  \mathcal{Z}^\alpha (p-1)\cdot( p^{-1}{\bf{v}}\cdot[\mathcal{Z}^\alpha , \nabla] (p-1)) \ d{\bf{x}}ds\\
	= &-\gamma^{-1}\sum_{k=0}^{m-1}\int_0^t\int_{\mathbb{R}_+^2} \phi_{k,m}(y) \mathcal{Z}^\alpha (p-1)\cdot  (p^{-1}{\bf{v}}\cdot\mathcal{Z}_2^k \partial_y (p-1) )\ d{\bf{x}}ds\\
	\lesssim&\;\|{\bf{v}}\|_{L_{t,{\bf{x}}}^\infty}\|  p^{-1}\|_{L_{t,{\bf{x}}}^\infty}\left(\int_0^t\|p-1\|_{m}^2 ds\right)^\frac12\left(\int_0^t\|\partial_y p\|_{m-1}^2ds\right)^\frac12.
	\end{align*}
	For the second term on the right hand side of \eqref{3.6}, since ${\bf{v}}|_{\partial\Omega}=0$, with integration by parts, we have
	\begin{align*}
	&\gamma^{-1}\int_0^t\int_{\mathbb{R}_+^2} \mathcal{Z}^\alpha (p-1)\cdot (  p^{-1}{\bf{v}}\cdot\nabla\mathcal{Z}^\alpha  (p-1)) \ d{\bf{x}}ds\\
	=&	-\frac{1}{2\gamma}\int_0^t\int_{\mathbb{R}_+^2}\nabla\cdot(  p^{-1}{\bf{v}})|\mathcal{Z}^\alpha (p-1)|^2 d{\bf{x}}ds\\
	\lesssim &\left(\|\nabla p^{-1}\|_{L_{t,{\bf{x}}}^\infty} \|{\bf{v}}\|_{L_{t,{\bf{x}}}^\infty}+\| p^{-1}\|_{L_{t,{\bf{x}}}^\infty} \|\nabla\cdot {\bf{v}}\|_{L_{t,{\bf{x}}}^\infty}\right)\int_0^t\|p-1\|_{m}^2 ds.
	\end{align*}
	The third term on the right hand side of \eqref{3.5} can be estimated as follows. By \eqref{2.1}, we have
	\begin{align*}
	&\gamma^{-1}\sum_{\beta+\gamma=\alpha\atop |\beta|\geq1} C_\alpha^\beta\int_0^t\int_{\mathbb{R}_+^2} \mathcal{Z}^\alpha (p-1)\cdot \mathcal{Z}^\beta  p^{-1}\mathcal{Z}^\gamma  \partial_t (p-1)\ d{\bf{x}}ds\\
	\lesssim&\;\|  p ^{-1}\|_{1,\infty} \int_0^t \|p-1\|_{m}^2 ds+\|  p\|_{1, \infty}\left(\int_0^t  \| p^{-1}-1 \|_{m}^2 ds\right)^\frac12\left( \int_0^t \| p-1\|_{m}^2 ds\right)^\frac12.
	\end{align*}
	Similarly, for the last term on the right hand side of \eqref{3.5},  it follows that
	\begin{align*}
	&\gamma^{-1}\sum_{\beta+\gamma=\alpha\atop |\beta|\geq1} C_\alpha^\beta\int_0^t\int_{\mathbb{R}_+^2} \mathcal{Z}^\alpha (p-1)  \cdot \mathcal{Z}^\beta (p^{-1}{\bf{v}})\cdot\mathcal{Z}^\gamma \nabla (p-1)\  d{\bf{x}}ds\\
	\lesssim &\;\| {\bf{v}}\|_{1, \infty}\|p^{-1}\|_{1, \infty}\left[\left(\int_0^t \|p-1\|_{m}^2 ds\right)^\frac12\left( \int_0^t \|\partial_y p\|_{m-1}^2 ds\right)^\frac12 + \int_0^t \|p-1\|_{m}^2 ds\right]\\
	&+ \|\nabla p\|_{L_{t,{\bf{x}}}^\infty}\|({\bf{v}}, p^{-1})\|_{1,\infty}\left(\int_0^t \|({\bf{v}},p^{-1}-1)\|_{m}^2ds\right)^\frac12\left(\int_0^t \|p-1\|_{m}^2ds\right)^\frac12.
	\end{align*}  	
Now, we come back to estimate the second term on the right hand side of \eqref{3.3}. By \eqref{2.2},  we obtain
	\begin{align*}
	-\int_0^t\int_{\mathbb{R}_+^2} \mathcal{Z}^\alpha (p-1)\cdot [\nabla\cdot, \mathcal{Z}^\alpha] {\bf{v}}\ d{\bf{x}}ds
	&=-\sum_{k=0}^{m-1}\int_0^t\int_{\mathbb{R}_+^2}\phi_{k, m}(y) \mathcal{Z}^\alpha (p-1)\cdot  \mathcal{Z}_2^k\partial_y {\bf{v}}\ d{\bf{x}}ds\\
	&\lesssim \left(\int_0^t\|p-1\|_{m}^2 ds\right)^\frac12\left(\int_0^t\|\partial_y {\bf{v}}\|_{m-1}^2 ds\right)^\frac12.
	\end{align*}
By the similar arguments, the third term on the right hand side of \eqref{3.3} is estimated as follows.
	\begin{align*}
	\int_0^t\int_{\mathbb{R}_+^2}[ \mathcal{Z}^\alpha,  \nabla] (p-1)\cdot \mathcal{Z}^\alpha {\bf{v}}\ d{\bf{x}}ds\lesssim &\left(\int_0^t\|\partial_y p\|_{m-1}^2 ds\right)^\frac12\left(\int_0^t\|{\bf{v}}\|_{m}^2ds\right)^\frac12.
	\end{align*}
	
    Next, we continue to estimate the second term on the right hand side of \eqref{3.2}. Since ${\bf{v}}|_{\partial\Omega}=0$, by integration by parts, we obtain
	\begin{equation}\label{3.7}
	\begin{split}
	\varepsilon \mu\int_0^t\int_{\mathbb{R}_+^2} \mathcal{Z}^\alpha \Delta {\bf{v}}\cdot \mathcal{Z}^\alpha {\bf{v}}\ d{\bf{x}}ds
	= &-\varepsilon \mu\int_0^t \|Z^\alpha\nabla {\bf{v}}\|^2 ds+\varepsilon \mu\int_0^t\int_{\mathbb{R}_+^2}[\mathcal{Z}^\alpha, \nabla\cdot]\nabla {\bf{v}}\cdot\mathcal{Z}^\alpha {\bf{v}}\ d{\bf{x}}ds\\
	&+\varepsilon \mu\int_0^t\int_{\mathbb{R}_+^2}\mathcal{Z}^\alpha\nabla {\bf{v}}\cdot[\mathcal{Z}^\alpha, \nabla] {\bf{v}} \ d{\bf{x}}ds.
	\end{split}
	\end{equation}
	For the second term on the right hand side of \eqref{3.7}, by \eqref{2.2} and  integration by parts, it follows that
	\begin{align*}
	 \varepsilon \mu\int_0^t\int_{\mathbb{R}_+^2}[\mathcal{Z}^\alpha, \nabla\cdot]\nabla {\bf{v}}\cdot\mathcal{Z}^\alpha {\bf{v}}\ d{\bf{x}}ds
	&=\;\varepsilon \mu\sum_{k=0}^{m-1}\int_0^t\int_{\mathbb{R}_+^2}\phi^{k,m}(y)\partial_y\mathcal{Z}_2^k\nabla {\bf{v}}\cdot\mathcal{Z}^\alpha {\bf{v}}\ d{\bf{x}}ds\\
	&=-\varepsilon \mu\sum_{k=0}^{m-1}\int_0^t\int_{\mathbb{R}_+^2}\partial_y\phi^{k,m}(y)\mathcal{Z}_2^k\nabla {\bf{v}}\cdot\mathcal{Z}^\alpha {\bf{v}}\ d{\bf{x}}ds\\
	&\quad-\varepsilon \mu\sum_{k=0}^{m-1}\int_0^t\int_{\mathbb{R}_+^2}\phi^{k,m}(y)\mathcal{Z}_2^k\nabla {\bf{v}}\cdot\partial_y\mathcal{Z}^\alpha {\bf{v}}\ d{\bf{x}}ds\\
	&\leq \varepsilon^2 \mu^2\int_0^t\|\partial_y{\bf{v}}\|_{m}^2 ds+C\left(\int_0^t\| {\bf{v}}\|_{m}^2 ds+\int_0^t\|\partial_y {\bf{v}}\|_{m-1}^2 ds\right).
	\end{align*}
	For the third term on the right hand side of \eqref{3.7}, by \eqref{2.2}, we conclude that
	\begin{align*}
	\varepsilon \mu\int_0^t\int_{\mathbb{R}_+^2}\mathcal{Z}^\alpha\nabla {\bf{v}}\cdot[\mathcal{Z}^\alpha, \nabla] {\bf{v}} \ d{\bf{x}}ds&=\varepsilon \mu\sum_{k=0}^{m-1}\int_0^t\int_{\mathbb{R}_+^2}\phi_{k, m}(y)\mathcal{Z}^\alpha\nabla {\bf{v}}\cdot\mathcal{Z}_2^k\partial_y{\bf{v}} \ d{\bf{x}}ds\\
	&\leq \varepsilon^2 \mu^2 \int_0^t\|\nabla {\bf{v}}\|_{m}^2 ds+C \int_0^t\|\partial_y {\bf{v}}\|_{m-1}^2 ds.
	\end{align*}
	
	Similarly, for the third term on the right hand side of \eqref{3.2},  by \eqref{2.2}, and  integration by parts,  one arrives at
	\begin{align*}	
	&\varepsilon(\mu+\lambda)\int_0^t\int_{\mathbb{R}_+^2} \mathcal{Z}^\alpha\nabla(\nabla\cdot {\bf{v}})\cdot \mathcal{Z}^\alpha {\bf{v}} \ d{\bf{x}}ds\\
	=&\;\varepsilon(\mu+\lambda)\int_0^t\int_{\mathbb{R}_+^2} \nabla\mathcal{Z}^\alpha (\nabla\cdot {\bf{v}})  \cdot\mathcal{Z}^\alpha {\bf{v}}\ d{\bf{x}}ds+\varepsilon(\mu+\lambda) \int_0^t\int_{\mathbb{R}_+^2} [\mathcal{Z}^\alpha, \nabla](\nabla\cdot {\bf{v}})\cdot \mathcal{Z}^\alpha {\bf{v}} \ d{\bf{x}}ds\\
	=&-\varepsilon(\mu+\lambda)\int_0^t\int_{\mathbb{R}_+^2}|\mathcal{Z}^\alpha (\nabla\cdot {\bf{v}})|^2\ d{\bf{x}}ds-\varepsilon(\mu+\lambda)\int_0^t\int_{\mathbb{R}_+^2} \mathcal{Z}^\alpha (\nabla\cdot {\bf{v}})\cdot[\nabla\cdot , \mathcal{Z}^\alpha] {\bf{v}}\ d{\bf{x}}ds\\
	&+\varepsilon(\mu+\lambda)\sum_{k=0}^{m-1}\int_0^t\int_{\mathbb{R}_+^2} \phi^{k,m}(y)\partial_y\mathcal{Z}_2^k(\nabla\cdot {\bf{v}}) \mathcal{Z}^\alpha {\bf{v}} \ d{\bf{x}}ds\\
	\leq &-\frac{\varepsilon(\mu+\lambda)}{2}\int_0^t\| \nabla\cdot {\bf{v}}\|_{m}^2 ds+\varepsilon^2 (\mu+\lambda)^2\int_0^t\|\partial_y{\bf{v}}\|_{m}^2 ds+ C\left(\int_0^t\|{\bf{v}}\|_{m}^2ds
	+ \int_0^t\|\partial_y {\bf{v}}\|_{m-1}^2ds\right).
	\end{align*}	

	As for the fourth term on the right hand side of \eqref{3.2}, by \eqref{2.1} and \eqref{2.2}, we have	
	\begin{align*}
	\int_0^t\int_{\mathbb{R}_+^2} \mathcal{C}_1^\alpha \cdot \mathcal{Z}^\alpha {\bf{v}} \ d{\bf{x}}ds& =-\sum_{\beta+\gamma=\alpha\atop |\beta|\geq 1} C_\alpha^\beta\int_0^t\int_{\mathbb{R}_+^2} \mathcal{Z}^\beta \rho\mathcal{Z}^\gamma\partial_t {\bf{v}}\cdot\mathcal{Z}^\alpha {\bf{v}}  \ d{\bf{x}}ds\\
	&\lesssim \| \rho\|_{1, \infty} \int_0^t\|  {\bf{v}}\|_{m}^2\ ds+\|  {\bf{v}}\|_{1, \infty}\left(\int_0^t\| \rho-1\|_{m}^2\ ds\right)^\frac12\left(\int_0^t \|{\bf{v}}\|_{m} ^2\ ds\right)^\frac12,
	\end{align*}
	and
	\begin{align*}
	\int_0^t\int_{\mathbb{R}_+^2}  \mathcal{C}_2^\alpha \cdot \mathcal{Z}^\alpha {\bf{v}} \ d{\bf{x}}ds
	=&-\sum_{\beta+\gamma=\alpha\atop |\beta|\geq 1}\int_0^t\int_{\mathbb{R}_+^2} \left(C_\alpha^\beta\mathcal{Z}^\beta(\rho {\bf{v}})\cdot\mathcal{Z}^\gamma \nabla {\bf{v}}\cdot \mathcal{Z}^\alpha {\bf{v}}
	- \rho {\bf{v}}\cdot[\mathcal{Z}^\alpha, \nabla]{\bf{v}}\cdot\mathcal{Z}^\alpha {\bf{v}}\right) d{\bf{x}}ds \\\
	\lesssim &\;\|(\rho, {\bf{v}})\|_{1,\infty}^2 \left[\left(\int_0^t \|\partial_y {\bf{v}}\|_{m-1}^2 ds\right)^\frac12+ \left(\int_0^t \| {\bf{v}}\|_{m}^2 ds\right)^\frac12 \right]\left(\int_0^t\|{\bf{v}}\|_{m}^2 ds\right)^\frac12\\
	&+\|\nabla {\bf{v}}\|_{L_{t,{\bf{x}}}^\infty} \|(\rho, {\bf{v}})\|_{1,\infty}\left(\int_0^t\|(\rho-1, {\bf{v}})\|_{m}^2 ds\right)^\frac12\left(\int_0^t\|{\bf{v}}\|_{m}^2 ds\right)^\frac12.
	\end{align*}

	Then we turn to deal with the penultimate term on the right hand side of \eqref{3.2}.
	\begin{equation}\label{3.8}
	\begin{split}
	&\int_0^t\int_{\mathbb{R}_+^2} \mathcal{Z}^\alpha[ (\nabla\times {\bf{B}})\times {\bf{B}}]\cdot \mathcal{Z}^\alpha {\bf{v}} \ d{\bf{x}}ds\\
	=&\int_0^t\int_{\mathbb{R}_+^2}  [\mathcal{Z}^\alpha(\nabla\times {\bf{B}})]\times {\bf{B}}\cdot\mathcal{Z}^\alpha {\bf{v}}\ d{\bf{x}}ds
	+\sum_{\beta+\gamma=\alpha\atop |\gamma|\geq 1}\int_0^t\int_{\mathbb{R}_+^2} [\mathcal{Z}^\beta (\nabla\times {\bf{B}})]\times\mathcal{Z}^\gamma {\bf{B}}\cdot \mathcal{Z}^\alpha {\bf{v}}\ d{\bf{x}}ds.
	\end{split}
	\end{equation}
	Then, for the second term on the right hand side of \eqref{3.8}, by \eqref{2.1}, it follows that
	\begin{align*}
	&\sum_{\beta+\gamma=\alpha\atop |\gamma|\geq 1}\int_0^t\int_{\mathbb{R}_+^2} \mathcal{Z}^\beta (\nabla\times {\bf{B}})\times\mathcal{Z}^\gamma {\bf{B}}\cdot \mathcal{Z}^\alpha {\bf{v}}\ d{\bf{x}}ds \\
	\lesssim &\left[\|{\bf{B}}\|_{1,\infty}\left(\int_0^t\|\nabla  {\bf{B}}\|_{m-1}^2 ds\right)^\frac12+\|\nabla  {\bf{B}}\|_{L_{t,{\bf{x}}}^\infty}\left(\int_0^t\|{\bf{B}}-\overset{\rightarrow}{e_y}\|_{m}^2 ds\right)^\frac12\right]\left(\int_0^t\|{\bf{v}}\|_{m}^2ds\right)^\frac12.
	\end{align*}
	
Finally, it holds for the last  term  on the right hand side of \eqref{3.2} that

	\begin{equation}\label{3.9}
	\begin{split}
	&\int_0^t\int_{\mathbb{R}_+^2}  \mathcal{Z}^\alpha[\nabla\times({\bf{v}}\times {\bf{B}})]\cdot\mathcal{Z}^\alpha ({\bf{B}}-\overset{\rightarrow}{e_y}) \ d{\bf{x}}ds \\
	=& \sum_{1\leq|\beta|,|\gamma|\leq |\alpha|-1}\int_0^t\int_{\mathbb{R}_+^2} \nabla\times( \mathcal{Z}^\beta {\bf{v}}\times \mathcal{Z}^\gamma {\bf{B}}) \cdot \mathcal{Z}^\alpha ({\bf{B}}-\overset{\rightarrow}{e_y}) \ d{\bf{x}}ds\\
	&+\int_0^t\int_{\mathbb{R}_+^2}  [\mathcal{Z}^\alpha, \nabla\times]({\bf{v}}\times {\bf{B}})\cdot\mathcal{Z}^\alpha ({\bf{B}}-\overset{\rightarrow}{e_y})\  d{\bf{x}}ds\\
	&+\int_0^t\int_{\mathbb{R}_+^2}  \nabla\times(\mathcal{Z}^\alpha {\bf{v}}\times {\bf{B}})\cdot \mathcal{Z}^\alpha ({\bf{B}}-\overset{\rightarrow}{e_y})\ d{\bf{x}}ds\\\
	&+\int_0^t\int_{\mathbb{R}_+^2}  \nabla\times( {\bf{v}}\times \mathcal{Z}^\alpha {\bf{B}})\cdot \mathcal{Z}^\alpha ({\bf{B}}-\overset{\rightarrow}{e_y})\ d{\bf{x}}ds.
	\end{split}
	\end{equation}	
	Where for the first term on the right hand side of \eqref{3.9},  by \eqref{2.1},   one has
	\begin{align*}
	&\sum_{1\leq|\beta|,|\gamma|\leq |\alpha|-1}\int_0^t\int_{\mathbb{R}_+^2} \nabla\times( \mathcal{Z}^\beta {\bf{v}}\times \mathcal{Z}^\gamma {\bf{B}}) \cdot \mathcal{Z}^\alpha ({\bf{B}}-\overset{\rightarrow}{e_y}) \ d{\bf{x}}ds\\
	\lesssim& \left(\|( {\bf{v}},  {\bf{B}})\|_{2, \infty}+\|\partial_y({\bf{v}},  {\bf{B}})\|_{1, \infty}\right)\left(\int_0^t\| ({\bf{v}}, {\bf{B}}-\overset{\rightarrow}{e_y})\|_{m-1}^2 ds\right)^\frac12\left(\int_0^t\| {\bf{B}}-\overset{\rightarrow}{e_y}\|_{m}^2 ds\right)^\frac12\\
	&+\|({\bf{v}}, {\bf{B}})\|_{1, \infty}\left[\left(\int_0^t\| ( {\bf{v}},   {\bf{B}}-\overset{\rightarrow}{e_y})\|_{m}^2ds\right)^\frac12+\left(\int_0^t\| \partial_y( {\bf{v}},   {\bf{B}})\|_{m-1}^2ds\right)^\frac12\right]\cdot\left(\int_0^t\| {\bf{B}}-\overset{\rightarrow}{e_y}\|_{m}^2ds\right)^\frac12.
	\end{align*}  		
	And for the second term on the right hand side of \eqref{3.9}, by \eqref{2.1} and \eqref{2.2}, we have
	\begin{align*}
	\int_0^t\int_{\mathbb{R}_+^2} [\mathcal{Z}^\alpha, \nabla\times](&{\bf{v}}\times {\bf{B}})\cdot\mathcal{Z}^\alpha ({\bf{B}}-\overset{\rightarrow}{e_y})d{\bf{x}}ds
	\lesssim \left(\int_0^t\|\partial_y ({\bf{v}}\times {\bf{B}})\|_{m-1}^2ds\right)^\frac12\left(\int_0^t\|{\bf{B}}-\overset{\rightarrow}{e_y}\|_{m}^2 ds\right)^\frac12\\
	&\lesssim \|\partial_y ({\bf{v}},  {\bf{B}}) \|_{L_{t,{\bf{x}}}^\infty}\left(\int_0^t\|({\bf{v}}, {\bf{B}}-\overset{\rightarrow}{e_y})\|_{m-1}^2ds\right)^\frac12\left(\int_0^t\|{\bf{B}}-\overset{\rightarrow}{e_y}\|_{m}^2ds\right)^\frac12\\
	&+ \|({\bf{v}}, {\bf{B}})\|_{L_{t,{\bf{x}}}^\infty}\left(\int_0^t\|\partial_y( {\bf{v}},   {\bf{B}})\|_{m-1}^2ds\right)^\frac12\left(\int_0^t\|{\bf{B}}-\overset{\rightarrow}{e_y}\|_{m}^2ds\right)^\frac12.
	\end{align*}
	By integration by parts, the third term  on the right hand side of \eqref{3.9} is handled as follows.
	\begin{equation}\label{3.10}
	\begin{split}
	\int_0^t\int_{\mathbb{R}_+^2}  \nabla\times(\mathcal{Z}^\alpha {\bf{v}}\times {\bf{B}})\cdot \mathcal{Z}^\alpha ({\bf{B}}-\overset{\rightarrow}{e_y})\ d{\bf{x}}ds
	&=\int_0^t\int_{\mathbb{R}_+^2}  (\mathcal{Z}^\alpha {\bf{v}}\times {\bf{B}})\cdot\mathcal{Z}^\alpha(\nabla\times {\bf{B}}) \ d{\bf{x}}ds\\
	&+\int_0^t\int_{\mathbb{R}_+^2}  \mathcal{Z}^\alpha {\bf{v}}\times {\bf{B}}\cdot [\nabla\times, \mathcal{Z}^\alpha] {\bf{B}}\  d{\bf{x}}ds.
	\end{split}
	\end{equation}
First, notice that the sum of the  first term on the right hand side of \eqref{3.10} and the  first term on the right hand side of \eqref{3.8} equals to zero.
For the second term on the right hand side of \eqref{3.10},  by \eqref{2.1} and \eqref{2.2}, we have
	\begin{align*}
	\int_0^t\int_{\mathbb{R}_+^2} \mathcal{Z}^\alpha {\bf{v}}\times {\bf{B}}\cdot [\nabla\times, \mathcal{Z}^\alpha] {\bf{B}} \ d{\bf{x}}ds
	\lesssim \|{\bf{B}}\|_{L_{t,{\bf{x}}}^\infty}\left(\int_0^t\| {\bf{v}}\|_{m}^2ds\right)^\frac12\left(\int_0^t\| \partial_y{\bf{B}}\|_{m-1}^2ds\right)^\frac12.
	\end{align*}

	It is left to estimate the last term on the right hand side of \eqref{3.9}. Since ${\bf{v}}|_{y=0}=0$ and $\nabla\cdot {\bf{B}}=0$, by integration by parts, we obtain
	\begin{equation}\label{3.11}
	\begin{split}
	&\int_0^t\int_{\mathbb{R}_+^2}   \nabla\times({\bf{v}}\times \mathcal{Z}^\alpha {\bf{B}})\cdot \mathcal{Z}^\alpha ({\bf{B}}-\overset{\rightarrow}{e_y})\ d{\bf{x}}ds\\
	=&\sum_{i, j=1,2\atop i\neq j} \int_0^t\int_{\mathbb{R}_+^2}  v_i\partial_j\mathcal{Z}^\alpha b_j\mathcal{Z}^\alpha \tilde b_i \ d{\bf{x}}ds-\sum_{i, j=1,2\atop i\neq j} \int_0^t\int_{\mathbb{R}_+^2} v_i\partial_i\mathcal{Z}^\alpha b_j \mathcal{Z}^\alpha \tilde b_j  d{\bf{x}}ds\\
	&+\sum_{i, j=1,2\atop i\neq j} \int_0^t\int_{\mathbb{R}_+^2} \left(\partial_j v_i\mathcal{Z}^\alpha b_i\mathcal{Z}^\alpha \tilde b_j -\partial_i v_i\mathcal{Z}^\alpha b_j\mathcal{Z}^\alpha \tilde b_j \right)d{\bf{x}}ds.
	\end{split}
	\end{equation}	
	For the first term on the right hand side of \eqref{3.11}, by the fact of $\nabla\cdot {\bf{B}}=0$, integration by parts and  \eqref{2.2} imply that
	\begin{align*}
	&\sum_{i, j=1,2\atop i\neq j} \int_0^t\int_{\mathbb{R}_+^2}  v_i\partial_j\mathcal{Z}^\alpha b_j\mathcal{Z}^\alpha \tilde b_i \ d{\bf{x}}ds\\
	=& \sum_{i, j=1,2\atop i\neq j}\int_0^t\int_{\mathbb{R}_+^2} \left( v_i[\partial_j, \mathcal{Z}^\alpha] \tilde b_j\mathcal{Z}^\alpha \tilde b_i d{\bf{x}}ds- v_i [\mathcal{Z}^\alpha, \partial_i]\tilde b_i \mathcal{Z}^\alpha \tilde b_i -v_i \partial_i\mathcal{Z}^\alpha \tilde b_i \mathcal{Z}^\alpha \tilde b_i\right)d{\bf{x}}ds\\
	\lesssim &\;\|{\bf{v}}\|_{L_{t,{\bf{x}}}^\infty}\left(\int_0^t\|\partial_y {\bf{B}}\|_{m-1}^2ds\right)^\frac12\left(\int_0^t\|{\bf{B}}-\overset{\rightarrow}{e_y}\|_{m}^2ds\right)^\frac12+ \|\nabla {\bf{v}}\|_{L_{t,{\bf{x}}}^\infty}\int_0^t\|{\bf{B}}-\overset{\rightarrow}{e_y}\|_{m}^2 ds.
	\end{align*}
	For the second term on the right hand side of \eqref{3.11}, by integration by parts, one has
	\begin{align*}
	-\sum_{i, j=1,2\atop i\neq j} \int_0^t\int_{\mathbb{R}_+^2} v_i\partial_i\mathcal{Z}^\alpha b_j \mathcal{Z}^\alpha \tilde b_j d{\bf{x}}ds\lesssim \|\nabla {\bf{v}}\|_{L_{t,{\bf{x}}}^\infty}\int_0^t\|{\bf{B}}-\overset{\rightarrow}{e_y}\|_{m}^2 ds.
	\end{align*}
	Similarly, for the last two terms on  the right hand side of \eqref{3.11}, we have
	\begin{align*}
	\sum_{i, j=1,2\atop i\neq j} \int_0^t\int_{\mathbb{R}_+^2} \left(\partial_j v_i\mathcal{Z}^\alpha b_i\mathcal{Z}^\alpha \tilde b_j -\partial_i v_i\mathcal{Z}^\alpha b_j\mathcal{Z}^\alpha \tilde b_j  \right)d{\bf{x}}ds
	\lesssim\|\nabla {\bf{v}}\|_{L_{t,{\bf{x}}}^\infty}\int_0^t\|{\bf{B}}-\overset{\rightarrow}{e_y}\|_{m}^2 ds.
	\end{align*}
	
	Collecting all of estimates above together, we arrive at the estimates in Lemma \ref{lem1}.
\end{proof}

\section{Normal derivative estimates}
In order to close the energy estimates established in Section 3,  based on the anisotropic Sobolev embedding property in the conormal  Sobolev space, it suffices to estimate  $\|\partial_y({\bf{v}}, {\bf{B}}, p)\|_{m-1}$ and $\|\partial_y^2({\bf{v}}, {\bf{B}}, p)\|_{m-2}$. Moreover, for the magnetic field, we only need to establish the conormal estimates of $\partial_yb_1$ and $\partial_y^2b_1$ due to the divergence free condition of $\nabla\cdot {\bf{B}}=0$.
Consequently, this section is devoted to the conormal estimates of the normal derivatives for the classical solutions $(\rho, {\bf{v}}, {\bf{B}})$ to compressible MHD  equations \eqref{3.1} with no-slip boundary condition \eqref{1.3}. Precisely, the conormal energy estimates for the  first  order and  the second order  normal derivatives of solutions $(\rho, {\bf{v}}, {\bf{B}})$ will be derived in details in the next subsections.

It is remarked that away from the boundary the estimates of $\|\partial_y({\bf{v}}, {\bf{B}}, p)\|_{m-1}$ and \\
$\|\partial_y^2({\bf{v}}, {\bf{B}}, p)\|_{m-2}$ are included in $\|({\bf{v}}, {\bf{B}}-\overset{\rightarrow}{e_y}, p-1)\|_m$. Consequently, we only need to derive estimates of normal derivatives near the boundary. In general, one can use the cut-off function technique. However, for the simplicity of presentation, we will omit the cut-off function in the analysis later.

Notice that if there exists a strong boundary layer, one can not expect the following uniform conormal estimates. This is exact the reason we claim the transverse magnetic field can prevent the strong boundary layer from occurring for compressible non-resistive MHD system with the no-slip boundary condition.
\subsection{The first order  normal derivatives} In this subsection, we consider the conormal estimate for  the first order normal derivative of  the classical solutions $(\rho, {\bf{v}}, {\bf{B}})$ to compressible  MHD equations \eqref{3.1} with no-slip boundary condition \eqref{1.3}.
\begin{lemma}\label{lem2}
	Under the assumption in Theorem \ref{Th1}, for $|\alpha|\leq m-1$, there exists a sufficiently small $\varepsilon_0>0$ such that for any $\varepsilon\in(0, \varepsilon_0)$,    the classical solutions $(\rho, {\bf{v}}, {\bf{B}})$ to compressible  MHD \eqref{3.1} with no-slip boundary condition \eqref{1.3} satisfies
	\begin{equation*}
	\begin{split}
	&\int_0^t\|\partial_y({\bf{v}}, b_1, p)\|_{m-1}^2 ds+ \varepsilon^2\mu^2 \int_0^t\|\partial_y^2 v_1\|_{m-1}^2ds+\varepsilon^2(2\mu+\lambda)^2\int_0^t\|\partial_y^2 v_2\|_{m-1}^2 ds\\
	&+\varepsilon(2\mu+\lambda)\gamma^{-1}\sum_{|\alpha|\leq m-1}\int_{\mathbb{R}_+^2} p^{-1}(t)|\mathcal{Z}^\alpha \partial_y p(t)|^2 d{\bf{x}}+\varepsilon\mu\|\partial_y b_1(t)\|_{m-1}^2\\
	\lesssim&\;\varepsilon(2\mu+\lambda)\gamma^{-1}\sum_{|\alpha|\leq m-1}\int_ {\mathbb{R}_+^2} p_0^{-1}|\mathcal{Z}^\alpha \partial_y p_0|^2 d{\bf{x}}+\varepsilon\mu\|\partial_y b_1(0)\|_{m-1}^2\\
	&+\left(1+\sum_{i=0}^1\|\partial_y^i(\rho, {\bf{v}},  {\bf{B}}, p, p^{-1})\|_{[(m-1)/2], \infty}^2\right)^2\int_0^t\|(\rho-1, {\bf{v}}, {\bf{B}}-\overset{\rightarrow}{e_y}, p-1, p^{-1}-1)\|_{m}^2 ds\\
	&+\varepsilon(2\mu+\lambda)\left[\left(1+\sum_{i =0}^1\|\partial_y^i( {\bf{v}},b_1, p, p^{-1})\|_{[(m-1)/2]+1, \infty}^2\right)^2+\|\partial_y^2v_2(0)\|_{[(m-1)/2]+1}^2\right.\\
	&\left.+\int_0^t\|\partial_y^2v_2\|_{[(m-1)/2]+2}^2 ds\right] \int_0^t\left(\| ({\bf{v}},b_1, p-1, p^{-1}-1)\|_{m}^2 +\|\partial_y({\bf{v}},b_1, p, p^{-1})\|_{m-1}^2 \right)ds\\
	&+ \varepsilon^2(2\mu+\lambda)^2\int_0^t\| \nabla{\bf{v}}\|_{m}^2 ds+\varepsilon^2(\mu+\lambda)^2\int_0^t\|\nabla\cdot{\bf{v}}\|_{m}^2 ds.
	\end{split}
	\end{equation*}
\end{lemma}

\subsubsection{Conormal Estimate of $\partial_y v_1$}  Rewrite the equation of $b_1$ in (\ref{3.1}) as follows
\begin{align}\label{4.1}
\partial_y v_1=	\partial_tb_1-\tilde{b}_2\partial_y v_1+v_1\partial_x b_1+v_2\partial_y b_1+b_1\partial_y v_2.
\end{align}
For any multi-index $\alpha$ satisfying $|\alpha|\leq m-1$, by applying $\mathcal{Z}^\alpha  $ to the above equality and taking  the $L^2$ inner product of the  resulting equality, we have
\begin{equation}\label{4.2}
\begin{split}
\int_0^t\int_{\mathbb{R}_+^2}| \mathcal{Z}^\alpha \partial_y v_1|^2\ d{\bf{x}}ds
&\lesssim \int_0^t\int_{\mathbb{R}_+^2}|	\partial_t\mathcal{Z}^\alpha  b_1|^2 d{\bf{x}}ds+\int_0^t\int_{\mathbb{R}_+^2}|\mathcal{Z}^\alpha  (\tilde{b}_2\partial_y v_1)|^2 d{\bf{x}}ds \\
&+\int_0^t\int_{\mathbb{R}_+^2}|\mathcal{Z}^\alpha  (v_1\partial_x b_1)|^2 d{\bf{x}}ds+\int_0^t\int_{\mathbb{R}_+^2}|\mathcal{Z}^\alpha  (v_2\partial_y b_1 )|^2 d{\bf{x}}ds\\ &+\int_0^t\int_{\mathbb{R}_+^2}|\mathcal{Z}^\alpha  (b_1\partial_y v_2)|^2 d{\bf{x}}ds.
\end{split}
\end{equation}

For the first term on the right hand side of \eqref{4.2},  one has
\begin{align*}
\int_0^t\int_{\mathbb{R}_+^2}|	\partial_t\mathcal{Z}^\alpha  b_1|^2 \ d{\bf{x}}ds
\lesssim\int_0^t\|  b_1\|_{m}^2 ds.
\end{align*}

For the second term on the right hand side of \eqref{4.2}, since $\nabla\cdot {\bf{B}}=0$ and\\ $\|\phi^{-1}\tilde b_2\|_{L_{t,{\bf{x}}}^\infty}\lesssim\|\partial_y\tilde b_2\|_{L_{t, {\bf{x}}}^\infty}$, by \eqref{2.1}, it holds that
\begin{equation*}
\begin{split}
&\int_0^t\int_{\mathbb{R}_+^2}|\mathcal{Z}^\alpha  (\tilde{b}_2\partial_y v_1)|^2 \ d{\bf{x}}ds\\
\lesssim &\;\|\phi^{-1}\tilde b_2\|_{L_{t, {\bf{x}}}^\infty}^2\int_0^t\|\phi\partial_y v_1\|_{m-1}^2 ds+\|\phi\partial_y v_1\|_{L_{t, {\bf{x}}}^\infty}^2\int_0^t\|\phi^{-1}\tilde b_2\|_{m-1}^2 ds\\
\lesssim&\;\|\partial_y\tilde b_2\|_{L_{t, {\bf{x}}}^\infty}^2\int_0^t\|  v_1\|_{m}^2 ds+\| v_1\|_{1, \infty}^2\int_0^t\| \partial_y\tilde b_2\|_{m-1}^2 ds\\
\lesssim&\;\|b_1\|_{1, \infty}^2\int_0^t\|v_1\|_{m}^2 ds+\|v_1\|_{1, \infty}^2\int_0^t\|b_1\|_{m}^2 ds.
\end{split}
\end{equation*}

For the third term on the right hand side of \eqref{4.2},  by \eqref{2.1}, we get
\begin{align*}
\int_0^t\int_{\mathbb{R}_+^2}|\mathcal{Z}^\alpha  (v_1\partial_x b_1)|^2 d{\bf{x}}ds
\lesssim&\;\|v_1\|_{L_{t, {\bf{x}}}^\infty}^2\int_0^t\|\partial_x b_1\|_{m-1}^2 ds+\|\partial_x b_1\|_{L_{t, {\bf{x}}}^\infty}^2\int_0^t\|v_1\|_{m-1}^2 ds\\
\lesssim&\;\|v_1\|_{L_{t, {\bf{x}}}^\infty}^2\int_0^t\| b_1\|_{m}^2 ds+\|  b_1\|_{ 1,\infty}^2\int_0^t\|v_1\|_{m-1}^2 ds.
\end{align*}

For the fourth  term on the right hand side of \eqref{4.2}, by   \eqref{2.3} and \eqref{2.4},  we have

\begin{align*}
\begin{split}
&\int_0^t\int_{\mathbb{R}_+^2}|\mathcal{Z}^\alpha  (v_2\partial_y b_1 )|^2 d{\bf{x}}ds \\
\lesssim&\;\| \phi\partial_y b_1\|_{L_{t, {\bf{x}}}^\infty}^2\int_0^t \|\phi^{-1}  v_2\|_{m-1} ds + \|\phi^{-1}  v_2\|_{L_{t, {\bf{x}}}^\infty}^2\int_0^t \|\phi \partial_y b_1\|_{m-1}^2 ds\\
\lesssim&\;\|b_1\|_{1, \infty}^2\int_0^t\|\partial_yv_2\|_{m-1}^2+\|\partial_yv_2\|_ {   L_{t, {\bf{x}}}^\infty}^2\int_0^t\|  b_1\|_{m}^2 ds.
\end{split}
\end{align*}

As for the last  term on the right hand side of \eqref{4.2},  by \eqref{2.1},  it follows that
\begin{align*}
\int_0^t\int_{\mathbb{R}_+^2}|\mathcal{Z}^\alpha  (b_1\partial_y v_2)|^2 d{\bf{x}}ds
\lesssim \|b_1\|_{L_{t, {\bf{x}}}^\infty}^2\int_0^t\|\partial_y v_2\|_{m-1}^2 ds+\|\partial_y v_2\|_{L_{t, {\bf{x}}}^\infty}^2\int_0^t\|b_1\|_{m-1}^2 ds.
\end{align*}

Combining all of the estimates above together, we obtain
\begin{align*}
\int_0^t\|  \partial_y v_1\|_{m-1}^2 ds
\lesssim& \left(1+  \sum_{i=0}^1\|\partial_y^i({\bf{v}}, b_1)\|_{ 1,\infty}^2\right)\int_0^t\|( {\bf{v}}, b_1 )\|_{m}^2 ds +\|b_1\|_{1, \infty}^2\int_0^t\|\partial_y v_2\|_{m-1}^2 ds.
\end{align*}
\subsubsection{Conormal Estimate of $  \partial_y v_2$}  This subsection is devoted to the conormal estimate for $  \partial_y v_2$.   With the help of the equation of the density in \eqref{3.1}, one has
\begin{align}\label{4.3}
\partial_y v_2=-\partial_x v_1-  \gamma^{-1}p^{-1}\partial_t p- \gamma^{-1}p^{-1}{\bf{v}}\cdot \nabla p.
\end{align}
For any multi-index $\alpha$ satisfying $|\alpha|\leq m-1$, by applying  $\mathcal{Z}^\alpha $ to the above equality and taking  the $L^2$ inner product on both sides of the resulting equality, it follows that
\begin{equation}\label{4.4}
\begin{split}
\int_0^t\int_{\mathbb{R}_+^2}|\mathcal{Z}^\alpha\partial_y v_2|^2 d{\bf{x}}ds
\lesssim& \int_0^t\int_{\mathbb{R}_+^2}|\mathcal{Z}^\alpha \partial_x v_1|^2d{\bf{x}}ds+ \int_0^t\int_{\mathbb{R}_+^2} |\mathcal{Z}^\alpha ( p^{-1}\partial_t p) |^2d{\bf{x}}ds\\
&+ \int_0^t\int_{\mathbb{R}_+^2} | \mathcal{Z}^\alpha (p^{-1}{\bf{v}}\cdot \nabla p) |^2 d{\bf{x}}ds.
\end{split}
\end{equation}
It is direct to estimate the first term on the right hand side of \eqref{4.4} as
\begin{align*}
\int_0^t\int_{\mathbb{R}_+^2}|\mathcal{Z}^\alpha \partial_x v_1|^2 d{\bf{x}}ds\lesssim\int_0^t \|  v_1\|_{m}^2 ds.
\end{align*}
For the second term on the right hand side of \eqref{4.4}, by \eqref{2.1}, we have
\begin{align*}
\int_0^t\int_{\mathbb{R}_+^2} |\mathcal{Z}^\alpha (  p^{-1}\partial_t p) |^2 d{\bf{x}}ds
\lesssim  \|p^{-1}\|_{L_{t, {\bf{x}}}^\infty}^2\int_0^t\|p-1\|_{m}^2 ds+\|p\|_{1, \infty}^2\int_0^t\|p^{-1}-1\|_{m-1}^2ds.
\end{align*}
The third term on the right hand side of \eqref{4.4} can be estimated in the following way.
\begin{align}\label{4.5}
\int_0^t\int_{\mathbb{R}_+^2} | \mathcal{Z}^\alpha (p^{-1}{\bf{v}}\cdot \nabla p) |^2d{\bf{x}}ds
\lesssim\int_0^t\int_{\mathbb{R}_+^2} | \mathcal{Z}^\alpha (p^{-1}v_1\partial_x p) |^2d{\bf{x}}ds+\int_0^t\int_{\mathbb{R}_+^2} | \mathcal{Z}^\alpha (p^{-1}v_2\partial_y p) |^2d{\bf{x}}ds.
\end{align}
For the first term on the right hand side of \eqref{4.5},  by \eqref{2.1}, one has
\begin{align*}
&\int_0^t\int_{\mathbb{R}_+^2} | \mathcal{Z}^\alpha (p^{-1}v_1\partial_x p) |^2d{\bf{x}}ds\\
\lesssim&\;\|v_1\|_{L_{t, {\bf{x}}}^\infty}^2\|\partial_xp\|_{L_{t, {\bf{x}}}^\infty}^2\int_0^t\|p^{-1}-1\|_{m-1}^2 ds +\|p^{-1}\|_{L_{t, {\bf{x}}}^\infty}^2\|\partial_xp\|_{L_{t, {\bf{x}}}^\infty}^2\int_0^t\|v_1\|_{m-1}^2 ds\\
&+\|p^{-1}\|_{L_{t, {\bf{x}}}^\infty}^2\|v_1\|_{L_{t, {\bf{x}}}^\infty}^2\int_0^t\|\partial_xp\|_{m-1}^2 ds\\
\lesssim&\;\|(v_1, p, p^{-1} )\|_{1, \infty}^4\int_0^t\|(v_1, p-1, p^{-1}-1)\|_{m}^2 ds.
\end{align*}
For the second term on the right hand side of \eqref{4.5}, by  \eqref{2.3} and \eqref{2.4}, we have
\begin{align*}
&\int_0^t\int_{\mathbb{R}_+^2} | \mathcal{Z}^\alpha (p^{-1}v_2\partial_y p) |^2d{\bf{x}}ds\\
\lesssim&\sum_{\beta+\gamma+\iota=\alpha\atop |\beta|\geq|\gamma|,|\iota|}\|\mathcal{Z}^\gamma v_2\|_{L_{t, {\bf{x}}}^\infty}^2\|\mathcal{Z}^\iota\partial_y p\|_{L_{t, {\bf{x}}}^\infty}^2\int_0^t\int_ {\mathbb{R}_+^2}|\mathcal{Z}^\beta p^{-1}|^2 d{\bf{x}}ds\\
&+\sum_{\beta+\gamma+\iota=\alpha\atop |\gamma| \geq |\beta|,|\iota|}\|\mathcal{Z}^\beta p^{-1}\|_{L_{t, {\bf{x}}}^\infty}^2\|\mathcal{Z}^\iota\partial_y p\|_{L_{t, {\bf{x}}}^\infty}^2\int_0^t\int_{\mathbb{R}_+^2}|\mathcal{Z}^\gamma v_2|^2 d{\bf{x}}ds\\
&+\sum_{\beta+\gamma+\iota=\alpha\atop |\iota| \geq |\beta|,|\gamma|}\|\mathcal{Z}^\beta p^{-1}\|_{L_{t, {\bf{x}}}^\infty}^2\|\phi^{-1}\mathcal{Z}^\gamma v_2\|_{L_{t, {\bf{x}}}^\infty}^2\int_0^t\int_ {\mathbb{R}_+^2}|\phi\mathcal{Z}^\iota\partial_y p|^2 d{\bf{x}}ds\\
\lesssim& \|(v_2,  p^{-1},  \partial_y v_2, \partial_y p)\|_{[(m-1)/2], \infty}^4 \int_0^t\|(v_2, p-1, p^{-1}-1)\|_{m}^2 ds.
\end{align*}

Inserting all of the estimates above into \eqref{4.4}, we conclude that
\begin{align*}
\int_0^t\| \partial_y v_2\|_{m-1}^2  ds
\lesssim&\left(1+\sum_{i=0}^1\|\partial_y^i({\bf{v}}, p, p^{-1}) \|_{[(m-1)/2], \infty}^2\right)^2\int_0^t \|({\bf{v}}, p-1, p^{-1}-1)\|_{m}^2  ds.
\end{align*}

\subsubsection{Conormal Estimate of $ \partial_y p$}  This subsection is devoted to the conormal estimate of $\partial_y p$. Since
\begin{align*}
\partial_y p  -\varepsilon (2\mu+\lambda)\partial_y^2 v_2
=-\rho \partial_t v_2+b_1\partial_x b_2-\rho {\bf{v}}\cdot \nabla v_2-b_1\partial_y b_1+\varepsilon\mu\partial_x^2 v_2+\varepsilon (\mu+\lambda)\partial_y\partial_x v_1,
\end{align*}
which comes from the equation of $v_2$ in (\ref{3.1}).
For any  multi-index $\alpha$ satisfying $|\alpha|\leq m-1$, by applying $\mathcal{Z}^\alpha $ to the above equality and taking the  $L^2$ inner product on both sides of the resulting equation, we have
\begin{equation}\label{4.6}
\begin{split}
&\int_0^t\int_{\mathbb{R}_+^2}\left( |\mathcal{Z}^\alpha \partial_y p |^2+ \varepsilon ^2(2\mu+\lambda)^2 |\mathcal{Z}^\alpha \partial_y^2 v_2|^2\right)  \ d{\bf{x}}ds -2\varepsilon (2\mu+\lambda)\int_0^t\int_{\mathbb{R}_+^2}\mathcal{Z}^\alpha  \partial_y^2  v_2\cdot \mathcal{Z}^\alpha \partial_y p \ d{\bf{x}}ds \\
\lesssim &\int_0^t\int_{\mathbb{R}_+^2} |\mathcal{Z}^\alpha(\rho \partial_t v_2)|^2 d{\bf{x}}ds+\int_0^t\int_{\mathbb{R}_+^2} |\mathcal{Z}^\alpha (b_1\partial_x b_2)|^2 d{\bf{x}}ds\\
&+\int_0^t\int_{\mathbb{R}_+^2} |\mathcal{Z}^\alpha (\rho {\bf{v}}\cdot \nabla v_2)|^2 d{\bf{x}}ds +\int_0^t\int_{\mathbb{R}_+^2} |\mathcal{Z}^\alpha  (b_1\partial_y b_1)|^2 d{\bf{x}}ds\\
&+\varepsilon^2\mu^2\int_0^t\int_{\mathbb{R}_+^2} |\mathcal{Z}^\alpha \partial_x^2 v_2|^2 d{\bf{x}}ds  +\varepsilon^2 (\mu+\lambda)^2\int_0^t\int_{\mathbb{R}_+^2}|\mathcal{Z}^\alpha\partial_y\partial_x v_1|^2 d{\bf{x}}ds.
\end{split}
\end{equation}
The second  term on the left hand side of \eqref{4.6} is handled as follows. By \eqref{4.3},  one has
\begin{equation}\label{4.7}
\begin{split}
&2\varepsilon (2\mu+\lambda)\int_0^t\int_{\mathbb{R}_+^2}\mathcal{Z}^\alpha \partial_y^2 v_2\cdot\mathcal{Z}^\alpha \partial_y p\ d{\bf{x}}ds\\
=&-2\gamma^{-1}\varepsilon (2\mu+\lambda)\int_0^t\int_{\mathbb{R}_+^2}\mathcal{Z}^\alpha \partial_y( p^{-1}\partial_t p)\cdot\mathcal{Z}^\alpha \partial_y p\ d{\bf{x}}ds\\
& -2\gamma^{-1}\varepsilon (2\mu+\lambda)\int_0^t\int_{\mathbb{R}_+^2}\mathcal{Z}^\alpha \partial_y (p^{-1}{\bf{v}} \cdot\nabla p )\cdot\mathcal{Z}^\alpha \partial_y p\ d{\bf{x}}ds\\
& -2\varepsilon (2\mu+\lambda)\int_0^t\int_{\mathbb{R}_+^2}\mathcal{Z}^\alpha \partial_y\partial_x v_1\cdot\mathcal{Z}^\alpha \partial_y p \ d{\bf{x}}ds.
\end{split}
\end{equation}
Then, we begin to estimate each term on the right hand side of (\ref{4.7}). The first term on the right hand side of \eqref{4.7} is decomposed as following two terms.
\begin{equation}\label{4.8}
\begin{split}
&-2\gamma^{-1}\varepsilon (2\mu+\lambda)\int_0^t\int_{\mathbb{R}_+^2}\mathcal{Z}^\alpha \partial_y ( p^{-1}\partial_t p)\cdot\mathcal{Z}^\alpha \partial_y p\ d{\bf{x}}ds\\
= &-2\gamma^{-1}\varepsilon (2\mu+\lambda)  \int_0^t\int_{\mathbb{R}_+^2}\mathcal{Z}^\alpha( \partial_y  p^{-1}  \partial_t p) \cdot\mathcal{Z}^\alpha \partial_y p\ d{\bf{x}}ds\\
&-2\gamma^{-1}\varepsilon (2\mu+\lambda)\int_0^t\int_{\mathbb{R}_+^2} \mathcal{Z}^\alpha   (  p^{-1} \partial_t\partial_y p ) \cdot\mathcal{Z}^\alpha \partial_y p\ d{\bf{x}}ds.
\end{split}
\end{equation}
For the first  term on the right hand side of \eqref{4.8},  by \eqref{2.1}, we have
\begin{align*}
&-2\gamma^{-1}\varepsilon (2\mu+\lambda)  \int_0^t\int_{\mathbb{R}_+^2}\mathcal{Z}^\alpha( \partial_y p^{-1}  \partial_t p) \cdot\mathcal{Z}^\alpha \partial_y p\ d{\bf{x}}ds\\
\lesssim&\;\varepsilon (2\mu+\lambda) \|\partial_y p^{-1}\|_{L_{t, {\bf{x}}}^\infty}\left(\int_0^t\|p-1\|_{m}^2 ds\right)^\frac12\left(\int_0^t\|\partial_y p\|_{m-1}^2 ds\right)^\frac12\\
&+\varepsilon (2\mu+\lambda) \|p\|_{ 1,\infty}\left(\int_0^t\|\partial_y p^{-1}\|_{m-1}^2 ds\right)^\frac12\left(\int_0^t\|\partial_y p\|_{m-1}^2 ds\right)^\frac12.
\end{align*}
Similarly, for the second term on the right hand side of \eqref{4.8}, one has
\begin{equation}\label{4.9}
\begin{split}
&-2\gamma^{-1}\varepsilon (2\mu+\lambda)\int_0^t\int_{\mathbb{R}_+^2} \mathcal{Z}^\alpha(   p^{-1}\partial_t\partial_y p ) \cdot\mathcal{Z}^\alpha \partial_y p\ d{\bf{x}}ds\\
=&-2\gamma^{-1}\varepsilon (2\mu+\lambda)\sum_{\beta+\gamma=\alpha\atop |\beta|\geq1}\int_0^t\int_{\mathbb{R}_+^2} \mathcal{Z}^\beta   p^{-1}  \mathcal{Z}^\gamma \partial_t\partial_y p   \cdot\mathcal{Z}^\alpha \partial_y p\ d{\bf{x}}ds\\
&-2\gamma^{-1}\varepsilon (2\mu+\lambda) \int_0^t\int_{\mathbb{R}_+^2}     p^{-1}\mathcal{Z}^\alpha \partial_t\partial_y p   \cdot\mathcal{Z}^\alpha \partial_y p\ d{\bf{x}}ds.
\end{split}
\end{equation}
And for the first term on the right hand side of \eqref{4.9},  it is estimated as follows.
\begin{equation*}
\begin{split}
&-2\gamma^{-1}\varepsilon (2\mu+\lambda)\sum_{\beta+\gamma=\alpha\atop |\beta|\geq1}\int_0^t\int_{ \mathbb{R}_+^2} \mathcal{Z}^\beta   p^{-1}  \mathcal{Z}^\gamma \partial_t\partial_y p   \cdot\mathcal{Z}^\alpha \partial_y p\ d{\bf{x}}ds\\
\lesssim&\; \varepsilon (2\mu+\lambda)\sum_{\beta+\gamma=\alpha\atop |\beta|\geq|\gamma|}\|\mathcal{Z}^\gamma \partial_t\partial_y p\|_{L_{t, {\bf{x}}}^\infty}\left(\int_0^t\int_ {\mathbb{R}_+^2}|\mathcal{Z}^\beta p^{-1}|^2 d{\bf{x}}ds\right)^\frac12\left(\int_0^t\|\partial_y p\|_{m-1}^2 ds\right)^\frac12\\
&+ \varepsilon (2\mu+\lambda)\sum_{\beta+\gamma=\alpha\atop 1\leq |\beta|<|\gamma|<|\alpha|}\|\mathcal{Z}^\beta p^{-1}\|_{L_{t, {\bf{x}}}^\infty}\left(\int_0^t\int_{\mathbb{R}_+^2}|\mathcal{Z}^\gamma \partial_t\partial_y p|^2 d{\bf{x}}ds\right)^\frac12\left(\int_0^t\|\partial_y p\|_{m-1}^2 ds\right)^\frac12\\
\lesssim&\; \varepsilon (2\mu+\lambda)\|\partial_y p\|_{[(m-1)/2]+1, \infty}\left(\int_0^t\|p^{-1}-1\|_{m-1}^2 ds\right)^\frac12\left(\int_0^t\|\partial_yp\|_{m-1}^2 ds\right)^\frac12\\
&+ \varepsilon (2\mu+\lambda)\|  p^{-1}\|_{[(m-1)/2], \infty}\int_0^t\|\partial_y p \|_{m-1}^2 ds.
\end{split}
\end{equation*}
For the second term on the right hand side of \eqref{4.9}, we have
\begin{align*}
&-2\gamma^{-1}\varepsilon (2\mu+\lambda) \int_0^t\int_{\mathbb{R}_+^2}    p^{-1}\mathcal{Z}^\alpha \partial_t\partial_y p   \cdot\mathcal{Z}^\alpha \partial_y p\ d{\bf{x}}ds\\
= &-\gamma^{-1}\varepsilon (2\mu+\lambda)\frac{d}{dt}\int_0^t\int_{\mathbb{R}_+^2}  p^{-1}|\mathcal{Z}^\alpha  \partial_y p |^2  \ d{\bf{x}}ds +\gamma^{-1}\varepsilon (2\mu+\lambda)\int_0^t\int_{\mathbb{R}_+^2} \partial_t p^{-1}|\mathcal{Z}^\alpha  \partial_y p |^2  \ d{\bf{x}}ds\\
\lesssim &-\varepsilon (2\mu+\lambda)\gamma^{-1} \int_{\mathbb{R}_+^2}     p^{-1}(t)|\mathcal{Z}^\alpha  \partial_y p (t)|^2  \ d{\bf{x}}+ \varepsilon (2\mu+\lambda)\gamma^{-1} \int_{\mathbb{R}_+^2}     p_0^{-1}|\mathcal{Z}^\alpha  \partial_y p_0 |^2  \ d{\bf{x}}  \\
&+\varepsilon (2\mu+\lambda)\|p^{-1}\|_{1, \infty}\int_0^t\|\partial_y p\|_{m-1}^2 ds.
\end{align*}
By the similar line, we deal with the second  term on the right hand side of \eqref{4.7} as follows.
\begin{equation}\label{4.10}
\begin{split}
& -2\gamma^{-1}\varepsilon (2\mu+\lambda)\int_0^t\int_{\mathbb{R}_+^2}\mathcal{Z}^\alpha \partial_y  (p^{-1}{\bf{v}} \cdot\nabla p)\cdot\mathcal{Z}^\alpha \partial_y p\ d{\bf{x}}ds\\
=&   -2\gamma^{-1}\varepsilon (2\mu+\lambda) \int_0^t\int_{\mathbb{R}_+^2}\mathcal{Z}^\alpha [ \partial_y(p^{-1}{\bf{v}})\cdot\nabla  p]\cdot\mathcal{Z}^\alpha \partial_yp\ d{\bf{x}}ds\\
& -2\gamma^{-1}\varepsilon (2\mu+\lambda) \int_0^t\int_{\mathbb{R}_+^2}\mathcal{Z}^\alpha  (p^{-1}{\bf{v}} \cdot\nabla \partial_yp )\cdot\mathcal{Z}^\alpha \partial_y p\ d{\bf{x}}ds.
\end{split}
\end{equation}
For the first term on the right hand side of \eqref{4.10},  by \eqref{2.1}, we have
\begin{align*}
&  -2\gamma^{-1}\varepsilon (2\mu+\lambda) \int_0^t\int_{\mathbb{R}_+^2}\mathcal{Z}^\alpha [ \partial_y(p^{-1}{\bf{v}})\cdot\nabla  p]\cdot\mathcal{Z}^\alpha \partial_y p\ d{\bf{x}}ds\\
=& -2\gamma^{-1}\varepsilon (2\mu+\lambda) \int_0^t\int_{\mathbb{R}_+^2}\mathcal{Z}^\alpha (\partial_yp^{-1}{\bf{v}}\cdot\nabla  p)\mathcal{Z}^\alpha \partial_y p\ d{\bf{x}}ds\\
&-2\gamma^{-1}\varepsilon (2\mu+\lambda) \int_0^t\int_{\mathbb{R}_+^2}\mathcal{Z}^\alpha (p^{-1}\partial_y{\bf{v}}\cdot\nabla  p)\mathcal{Z}^\alpha \partial_y p\ d{\bf{x}}ds\\
\lesssim&\;\varepsilon (2\mu+\lambda) \sum_{i=0}^1\|\partial_y^i({\bf{v}}, p, p^{-1})\|_{1, \infty}^2\cdot \bigg[\int_0^t\Big(\sum_{j=0}^1\|\partial_y^j({\bf{v}}, p-1, p^{-1}-1)\|_{m-1}^2\\
&+\|p-1\|_{m}^2\Big)ds\bigg]^\frac12 \left(\int_0^t\|\partial_y p\|_{m-1}^2 ds\right)^\frac12.
\end{align*}
And the second term on the right hand side of \eqref{4.10} can be decomposed into two parts.
\begin{equation}\label{4.11}
\begin{split}
& -2\gamma^{-1}\varepsilon (2\mu+\lambda) \int_0^t\int_{\mathbb{R}_+^2}\mathcal{Z}^\alpha (p^{-1}{\bf{v}} \cdot\nabla \partial_yp)\cdot\mathcal{Z}^\alpha \partial_y p\ d{\bf{x}}ds\\
= & -2\gamma^{-1}\varepsilon (2\mu+\lambda)\sum_{\beta+\gamma+\iota=\alpha\atop \iota\neq \alpha} \int_0^t\int_{\mathbb{R}_+^2}\mathcal{Z}^\beta  p^{-1}\mathcal{Z}^\gamma {\bf{v}}\cdot\mathcal{Z}^\iota\nabla \partial_yp\cdot\mathcal{Z}^\alpha \partial_y p\ d{\bf{x}}ds\\
& -2\gamma^{-1}\varepsilon (2\mu+\lambda) \int_0^t\int_{\mathbb{R}_+^2} p^{-1}{\bf{v}}\cdot\mathcal{Z}^\alpha \nabla\partial_yp \cdot\mathcal{Z}^\alpha \partial_y p\ d{\bf{x}}ds.
\end{split}
\end{equation}
By \eqref{2.3} and \eqref{2.4}, the first term on the right hand side of \eqref{4.11} is estimated as follows.
\begin{align*}
& -2\gamma^{-1}\varepsilon (2\mu+\lambda)\sum_{\beta+\gamma+\iota=\alpha\atop \iota\neq \alpha} \int_0^t\int_{\mathbb{R}_+^2}\mathcal{Z}^\beta  p^{-1}\mathcal{Z}^\gamma {\bf{v}}\cdot\mathcal{Z}^\iota\nabla \partial_yp\cdot\mathcal{Z}^\alpha \partial_y p\ d{\bf{x}}ds\\
\lesssim&\;\varepsilon (2\mu+\lambda)\left[\sum_{\beta+\gamma+\iota=\alpha\atop \gamma|,|\iota|\leq|\beta|}\| \mathcal{Z}^\gamma v_1\|_{L_{t, {\bf{x}}}^\infty}\| \mathcal{Z}^\iota\partial_x\partial_yp\|_{L_{t, {\bf{x}}}^\infty}\left(\int_0^t\int_{\mathbb{R}_+^2}|\mathcal{Z}^\beta  p^{-1}|^2d{\bf{x}} ds\right)^\frac12\right.\\
&\left.+\sum_{\beta+\gamma+\iota=\alpha\atop \gamma|,|\iota|\leq|\beta|}\|\phi^{-1}\mathcal{Z}^\gamma v_2\|_{L_{t, {\bf{x}}}^\infty}\|\phi\mathcal{Z}^\iota \partial_y^2p\|_{L_{t, {\bf{x}}}^\infty}\left(\int_0^t\int_{\mathbb{R}_+^2}|\mathcal{Z}^\beta  p^{-1}|^2d{\bf{x}} ds\right)^\frac12\right.\\
&\left.+\sum_{\beta+\gamma+\iota=\alpha\atop |\beta|,|\iota|\leq |\gamma|  }\|\mathcal{Z}^\beta  p^{-1}\|_{L_{t, {\bf{x}}}^\infty}\|\mathcal{Z}^\iota\partial_x\partial_yp\|_{L_{t, {\bf{x}}}^\infty}\left(\int_0^t\int_{\mathbb{R}_+^2}| \mathcal{Z}^\gamma v_1|^2d{\bf{x}} ds\right)^\frac12\right.\\
&\left.+\sum_{\beta+\gamma+\iota=\alpha\atop |\beta|,|\iota|\leq |\gamma|  }\|\mathcal{Z}^\beta  p^{-1}\|_{L_{t, {\bf{x}}}^\infty}\|\phi\mathcal{Z}^\iota\partial_y^2p\|_{L_{t, {\bf{x}}}^\infty}\left(\int_0^t\int_{\mathbb{R}_+^2}|\phi^{-1}\mathcal{Z}^\gamma v_2|^2d{\bf{x}} ds\right)^\frac12\right.\\
&\left.+\sum_{\beta+\gamma+\iota=\alpha\atop |\beta|,|\gamma|\leq |\iota|<|\alpha| }\|\mathcal{Z}^\beta  p^{-1}\|_{L_{t, {\bf{x}}}^\infty}\|\phi^{-1}\mathcal{Z}^\gamma {\bf{v}}\|_{L_{t, {\bf{x}}}^\infty}\left(\int_0^t\int_{\mathbb{R}_+^2}|\phi\mathcal{Z}^\iota\nabla \partial_yp|^2d{\bf{x}} ds\right)^\frac12 \right]\\
&\cdot\left(\int_0^t\|\partial_y p\|_{m-1}^2 ds\right)^\frac12\\
\lesssim&\;\varepsilon (2\mu+\lambda)\|({\bf{v}}, p^{-1}, \partial_y {\bf{v}})\|_{[(m-1)/2], \infty}\|(\partial_y {\bf{v}},  \partial_y p)\|_{[(m-1)/2]+1, \infty}\\
&\cdot\left[\int_0^t\left(\sum_{i=0}^1\|\partial_y^i( {\bf{v}},   p-1, p^{-1}-1)\|_{m-1}^2+\|p-1\|_{m}^2\right)ds\right] ^\frac12\left(\int_0^t\|\partial_y p\|_{m-1}^2 ds\right)^\frac12.
\end{align*}
And the last term on the right hand side of \eqref{4.11} satisfies the following equality.
\begin{equation}\label{4.12}
\begin{split}
&-2\gamma^{-1}\varepsilon (2\mu+\lambda) \int_0^t\int_{\mathbb{R}_+^2} p^{-1}{\bf{v}}\cdot\mathcal{Z}^\alpha \nabla\partial_yp \cdot\mathcal{Z}^\alpha \partial_y p\ d{\bf{x}}ds\\
= &-2\gamma^{-1}\varepsilon (2\mu+\lambda) \int_0^t\int_{\mathbb{R}_+^2} p^{-1}{\bf{v}}\cdot[\mathcal{Z}^\alpha, \nabla] \partial_y p \cdot\mathcal{Z}^\alpha \partial_y p\ d{\bf{x}}ds\\
&-2\gamma^{-1}\varepsilon (2\mu+\lambda) \int_0^t\int_{\mathbb{R}_+^2} p^{-1}{\bf{v}}\cdot\nabla\mathcal{Z}^\alpha\partial_y p \cdot\mathcal{Z}^\alpha \partial_y p\ d{\bf{x}}ds.
\end{split}
\end{equation}
For the first term on the right hand side of \eqref{4.12}, by \eqref{2.2}, we have
\begin{equation*}
\begin{split}
&-2\gamma^{-1}\varepsilon (2\mu+\lambda) \int_0^t\int_{\mathbb{R}_+^2} p^{-1}{\bf{v}}\cdot[\mathcal{Z}^\alpha, \nabla]   \partial_yp \cdot\mathcal{Z}^\alpha \partial_y p\ d{\bf{x}}ds\\
=&-2\gamma^{-1}\varepsilon (2\mu+\lambda)\sum_{k=0}^{m-2} \int_0^t\int_{\mathbb{R}_+^2}\phi^{k, m-1} (y)p^{-1}\phi^{-1}{\bf{v}}\cdot\phi\partial_y\mathcal{Z}_2^k\partial_y p \cdot\mathcal{Z}^\alpha \partial_y p\ d{\bf{x}}ds\\
\lesssim&\;\varepsilon (2\mu+\lambda)\|p^{-1}\|_{L_{t, {\bf{x}}}^\infty}\|\partial_y {\bf{v}}\|_{L_{t, {\bf{x}}}^\infty} \int_0^t\|\partial_y p\|_{m-1}^2 ds.
\end{split}
\end{equation*}
For the second term on the right hand side of \eqref{4.12}, by integration by parts,  one has

\begin{align*}
\begin{split}
&-2\gamma^{-1}\varepsilon (2\mu+\lambda) \int_0^t\int_{\mathbb{R}_+^2} p^{-1}{\bf{v}}\cdot\nabla\mathcal{Z}^\alpha\partial_y p \cdot\mathcal{Z}^\alpha \partial_y p\ d{\bf{x}}ds\\
\lesssim&\;\varepsilon (2\mu+\lambda)\sum_{i =0}^1\|\partial_y^i({\bf{v}}, p^{-1})\|_{1,\infty}^2\int_0^t\|\partial_y p\|_{m-1}^2 ds.
\end{split}
\end{align*}
Finally, for the last term on the right hand side of \eqref{4.7}, we have
\begin{align*}
& -2\varepsilon (2\mu+\lambda)\int_0^t\int_{\mathbb{R}_+^2}\mathcal{Z}^\alpha \partial_y\partial_x v_1\cdot\mathcal{Z}^\alpha \partial_y p \ d{\bf{x}}ds\\
\leq &\;2\varepsilon^2 (2\mu+\lambda)^2\int_0^t\|\partial_y v_1\|_{m}^2 ds + \frac12 \int_0^t\|\partial_y p\|_{m-1}^2 ds.
\end{align*}

Next, we turn to estimate the terms on the right hand side of \eqref{4.6}. For the first two terms on the right hand side of \eqref{4.6},  by \eqref{2.1}, one has
\begin{align*}
&	\int_0^t\int_{\mathbb{R}_+^2} |\mathcal{Z}^\alpha  (\rho \partial_t v_2)|^2 d{\bf{x}}ds+ \int_0^t\int_{\mathbb{R}_+^2} |\mathcal{Z}^\alpha (b_1\partial_x b_2)|^2  d{\bf{x}}ds\\
\lesssim &\;\|\rho \|_{L_{t,{\bf{x}}}^\infty}^2\int_0^t\|  v_2\|_{m}^2 ds+\|  v_2 \|_{1,\infty}^2\int_0^t\| \rho-1\|_{m-1}^2 ds\\
&+ \|  b_1\|_{L_{t,{\bf{x}}}^\infty}^2\int_0^t\|  \tilde{b}_2\|_{m}^2ds+\| \tilde{b}_2\|_{1,\infty}^2\int_0^t\| b_1\|_{m-1}^2ds.
\end{align*}
The third and the fourth terms on the right hand side of \eqref{4.6} can be handled as follows.
\begin{align*}
&\int_0^t\int_{\mathbb{R}_+^2} |\mathcal{Z}^\alpha (\rho {\bf{v}}\cdot \nabla v_2)|^2 d{\bf{x}}ds+	\int_0^t\int_{\mathbb{R}_+^2} |\mathcal{Z}^\alpha (b_1\partial_y b_1)|^ 2 d{\bf{x}}ds \\
\lesssim &\;\|(\rho, v_1, v_2)\|_{1, \infty}^4 \int_0^t \|(\rho, v_1, v_2)\|_{m}^2 ds+ \sum_{\beta+\gamma+\iota=\alpha\atop |\beta|\geq |\gamma|, |\iota|}\|\mathcal{Z}^\gamma v_2\|_{L_{t, {\bf{x}}}^\infty}^2\|\mathcal{Z}^\iota\partial_y v_2\|_{L_{t, {\bf{x}}}^\infty}^2\int_0^t\int_{\mathbb{R}_+^2}|\mathcal{Z}^\beta \rho|^2 d{\bf{x}}ds\\
&+ \sum_{\beta+\gamma+\iota=\alpha\atop |\gamma|\geq |\beta|, |\iota|}\|\mathcal{Z}^\beta \rho\|_{L_{t, {\bf{x}}}^\infty}^2\|\mathcal{Z}^\iota\partial_y v_2\|_{L_{t, {\bf{x}}}^\infty}^2\int_0^t\int_{\mathbb{R}_+^2}|\mathcal{Z}^\gamma v_2|^2 d{\bf{x}}ds\\
&+ \sum_{\beta+\gamma+\iota=\alpha\atop |\iota|\geq |\beta|, |\gamma|}\|\mathcal{Z}^\beta \rho\|_{L_{t, {\bf{x}}}^\infty}^2\|\phi^{-1}\mathcal{Z}^\gamma v_2\|_{L_{t, {\bf{x}}}^\infty}^2\int_0^t\int_{\mathbb{R}_+^2}|\phi\mathcal{Z}^\iota\partial_y v_2|^2 d{\bf{x}}ds\\
&+\|b_1\|_{L_{t,{\bf{x}}}^\infty}^2\int_0^t\|\partial_y b_1\|_{m-1}^2 ds+ \|\partial_yb_1\|_{L_{t,{\bf{x}}}^\infty}^2\int_0^t\|b_1\|_{m-1}^2 ds\notag\\
\lesssim&\;\left(1+\sum_{i=0}^1\|\partial_y^i(\rho, {\bf{v}}, b_1)\|_{[(m-1)/2], \infty}^2\right)^2\int_0^t\|(\rho-1,  {\bf{v}}, b_1)\|_{m}^2 dx+\|b_1\|_{L_{t,{\bf{x}}}^\infty}^2\int_0^t\|\partial_y b_1\|_{m-1}^2 ds.
\end{align*}
For the last two terms on the right hand side of \eqref{4.6},  we have
\begin{align*}
&\varepsilon^2\mu^2 \int_0^t\int_{\mathbb{R}_+^2} |\mathcal{Z}^\alpha \partial_x^2 v_2|^2d{\bf{x}}ds+\varepsilon^2 (\mu+\lambda)^2\ \int_0^t\int_{\mathbb{R}_+^2}|\mathcal{Z}^\alpha\partial_y\partial_x v_1|^2 d{\bf{x}}ds\\
\lesssim &\;\varepsilon^2\mu^2\int_0^t\| \partial_x v_2\|_{m}^2 ds+\varepsilon^2 (\mu+\lambda)^2\int_0^t\|\partial_yv_1\|_{m}^2 ds.
\end{align*}
Combining all of the above estimates  in this subsection yields that
\begin{align*}
&\int_0^t \left(\|  \partial_y p \|_{m-1}^2+ \varepsilon ^2(2\mu+\lambda)^2 \|  \partial_y^{2} v_2\|_{m-1}^2\right)  ds  +\varepsilon(2\mu+\lambda)\gamma^{-1}\sum_{|\alpha|\leq m-1}\int_{\mathbb{R}_+^2} p^{-1}(t)|\mathcal{Z}^\alpha \partial_y p(t)|^2 d{\bf{x}}\\
\lesssim &\;\varepsilon(2\mu+\lambda)\gamma^{-1}\sum_{|\alpha|\leq m-1}\int_{\mathbb{R}_+^2} p_0^{-1}|\mathcal{Z}^\alpha \partial_y p_0|^2 d{\bf{x}}  +\varepsilon^2 (2\mu+\lambda)^2\int_0^t\|\nabla {\bf{v}}\|_{m}^2 ds +\|b_1\|_{L_{t,{\bf{x}}}^\infty}^2\int_0^t\|\partial_y b_1\|_{m-1}^2 ds\notag\\
&+\varepsilon(2\mu+\lambda) \sum_{i=0}^1 \|\partial_y^i({\bf{v}}, p, p^{-1})\|_{[(m-1)/2]+1, \infty}^2  \int_0^t\Big(\sum_{j=0}^1\|\partial_y^j({\bf{v}}, p-1, p^{-1}-1) \|_{m-1}^2\\
& +\|p-1\|_{m}^2  \Big)ds+\left(1+\sum_{i=0}^1\|\partial_y^i(\rho, {\bf{v}},  {\bf{B}})\|_{[(m-1)/2],\infty}^2 \right)^2 \int_0^t \|(\rho-1, {\bf{v}},  {\bf{B}}-\overset{\rightarrow}{e_y})\|_{m}^2  ds.
\end{align*}

It should be remarked that the smallness of $\|b_1\|_{L_{t,{\bf{x}}}^\infty}$ in the third term in the right hand side of the above inequality will be used later. This is the reason why the small condition (\ref{1.5}) is required in Theorem \ref{Th1}.

\subsubsection{The Conormal Estimate of  $  \partial_y b_1$ } In this subsection, we will derive the conormal estimates of $ \partial_y b_1$. Rewrite the equation of $v_1$ as follows.
\begin{align*}
\partial_y b_1+\varepsilon\mu\partial_y^2 v_1=\rho\partial_t v_1+\partial_x p+b_2\partial_x b_2-\tilde b_2\partial_y b_1+\rho {\bf{v}}\cdot\nabla v_1-\varepsilon \mu \partial_x^2 v_1-\varepsilon(\mu+\lambda)\partial_x(\nabla\cdot {\bf{v}}) .
\end{align*}
For any  multi-index $\alpha$ satisfying $|\alpha|\leq m-1$, applying  $\mathcal{Z}^\alpha $ to  the  above equality and taking the $L^2$ inner product of the resulting equality give that
\begin{equation}\label{4.14}
\begin{split}
&\int_0^t\int_{\mathbb{R}_+^2}\left( |\mathcal{Z}^\alpha\partial_y b_1|^2+\varepsilon^2\mu^2|\mathcal{Z}^\alpha\partial_y^{2}  v_1|^2 \right)d{\bf{x}}ds +2\varepsilon\mu\int_0^t\int_{\mathbb{R}_+^2}\mathcal{Z}^\alpha\partial_y^{2}  v_1 \cdot \mathcal{Z}^\alpha\partial_y b_1 d{\bf{x}}ds\\
\lesssim&\int_0^t\int_{\Omega} |\mathcal{Z}^\alpha (\rho\partial_t v_1)|^2 d{\bf{x}}ds+\int_0^t\int_{\mathbb{R}_+^2} |\mathcal{Z}^\alpha \partial_x p|^2 d{\bf{x}}ds+\int_0^t\int_{\mathbb{R}_+^2} |\mathcal{Z}^\alpha (b_2\partial_x b_2)|^2 d{\bf{x}}ds\\
&+\int_0^t\int_{\mathbb{R}_+^2} |\mathcal{Z}^\alpha (\tilde b_2\partial_y b_1)|^2 d{\bf{x}}ds   +\int_0^t\int_{\mathbb{R}_+^2} |\mathcal{Z}^\alpha (\rho {\bf{v}}\cdot\nabla v_1)|^2 d{\bf{x}}ds\\
&+\varepsilon^2 \mu ^2\int_0^t\int_{\mathbb{R}_+^2}|\mathcal{Z}^\alpha \partial_x^2 v_1|^2d{\bf{x}}ds+\varepsilon^2(\mu+\lambda)^2\int_0^t\int_{\mathbb{R}_+^2} |\mathcal{Z}^\alpha \partial_x (\nabla\cdot {\bf{v}})|^2 d{\bf{x}}ds.
\end{split}
\end{equation}
To handle the second term on the left hand side of \eqref{4.14}, we apply $\mathcal{Z}^\alpha\partial_y$ to the following equation of $b_1$:
\begin{align*}
-\partial_y v_1=-\partial_tb_1+\tilde{b}_2\partial_y v_1-v_2\partial_y b_1-b_1\partial_y v_2 -v_1\partial_x b_1
\end{align*}
and take the $L^2$ inner product with $2\varepsilon\mu \mathcal{Z}^\alpha\partial_y b_1$ on both sides of  the resulting equality to get
\begin{equation}\label{4.15}
\begin{split}
&-2\varepsilon\mu\int_0^t\int_{\mathbb{R}_+^2} \mathcal{Z}^\alpha\partial_y^{2} v_1 \cdot\mathcal{Z}^\alpha\partial_y b_1\ d{\bf{x}}ds\\
=&-2\varepsilon\mu\int_0^t\int_{\mathbb{R}_+^2} \partial_t \mathcal{Z}^\alpha\partial_y	b_1\cdot\mathcal{Z}^\alpha\partial_y b_1\ d{\bf{x}}ds+2\varepsilon\mu\int_0^t\int_{\mathbb{R}_+^2}  \mathcal{Z}^\alpha\partial_y(\tilde{b}_2\partial_y v_1)\cdot \mathcal{Z}^\alpha\partial_y b_1\ d{\bf{x}}ds\\
&-2\varepsilon\mu\int_0^t\int_{\mathbb{R}_+^2} \mathcal{Z}^\alpha\partial_y(v_2\partial_y b_1)\cdot\mathcal{Z}^\alpha\partial_y b_1\ d{\bf{x}}ds-2\varepsilon\mu\int_0^t\int_{\mathbb{R}_+^2} \mathcal{Z}^\alpha\partial_y(b_1\partial_y v_2)\cdot\mathcal{Z}^\alpha\partial_y b_1\ d{\bf{x}}ds\\
&-2\varepsilon\mu\int_0^t\int_{\mathbb{R}_+^2} \mathcal{Z}^\alpha\partial_y(v_1\partial_x b_1)\cdot \mathcal{Z}^\alpha\partial_y b_1\ d{\bf{x}}ds.
\end{split}
\end{equation}
In this way, it suffices to estimate the terms on the right hand side of \eqref{4.15}.  First,
\begin{align*}
-2\varepsilon\mu\int_0^t\int_{\mathbb{R}_+^2}\partial_t \mathcal{Z}^\alpha\partial_y	b_1 \cdot\mathcal{Z}^\alpha\partial_y b_1 \ d{\bf{x}}ds
=-\varepsilon\mu\int_{\mathbb{R}_+^2}|\mathcal{Z}^\alpha\partial_y b_1(t)|^2 d{\bf{x}}+\varepsilon\mu\int_{\mathbb{R}_+^2}|\mathcal{Z}^\alpha\partial_y b_1(0)|^2 d{\bf{x}}.
\end{align*}
Second,
\begin{equation}\label{4.16}
\begin{split}
2\varepsilon\mu\int_0^t\int_{\mathbb{R}_+^2} \mathcal{Z}^\alpha\partial_y(\tilde{b}_2\partial_y v_1)\cdot\mathcal{Z}^\alpha\partial_y b_1 \ d{\bf{x}}ds
=&2\varepsilon\mu \int_0^t\int_{\mathbb{R}_+^2} \mathcal{Z}^\alpha(\partial_y \tilde{b}_2\partial_y  v_1)\cdot\mathcal{Z}^\alpha\partial_y b_1\ d{\bf{x}}ds\\
&+2\varepsilon\mu\int_0^t\int_{\mathbb{R}_+^2} \mathcal{Z}^\alpha(\tilde{b}_2\partial_y^{2} v_1)\cdot\mathcal{Z}^\alpha\partial_y b_1\ d{\bf{x}}ds.
\end{split}
\end{equation}
Here, for the first  term on the right hand side of \eqref{4.16},  since $\nabla\cdot {\bf{B}}=0$, by \eqref{2.1}, we have
\begin{align*}
&2\varepsilon\mu \int_0^t\int_{\Omega}
\mathcal{Z}^\alpha(\partial_y \tilde{b}_2\partial_y  v_1)\cdot\mathcal{Z}^\alpha\partial_y b_1\ d{\bf{x}}ds\\
\lesssim&\;\varepsilon\mu\| b_1\|_{1, \infty}\left(\int_0^t\|\partial_y v_1\|_{m-1}^2 ds\right)^\frac12\left(\int_0^t\|\partial_y b_1\|_{m-1}^2 ds\right)^\frac12\\
&+\varepsilon\mu\|\partial_y v_1\|_{L_{t,{\bf{x}}}^\infty}\left(\int_0^t\|b_1\|_{m}^2 ds\right)^\frac12\left(\int_0^t\|\partial_y b_1\|_{m-1}^2 ds\right)^\frac12.
\end{align*}
And for the second term on the right hand of \eqref{4.16},  one has
\begin{equation}\label{4.48}
\begin{split}
&2\varepsilon\mu\int_0^t\int_{\mathbb{R}_+^2} \mathcal{Z}^\alpha(\tilde{b}_2\partial_y^{2} v_1)\cdot\mathcal{Z}^\alpha\partial_y b_1\ d{\bf{x}}ds\\
=&2\varepsilon\mu\sum_{\beta+\gamma=\alpha\atop |\beta|\geq 1} C_\alpha^\beta\int_0^t\int_{\mathbb{R}_+^2} \mathcal{Z}^\beta\tilde{b}_2\mathcal{Z}^\gamma\partial_y^{2} v_1\cdot\mathcal{Z}^\alpha\partial_y b_1\ d{\bf{x}}ds\\
&+2\varepsilon\mu\int_0^t\int_{\mathbb{R}_+^2} \tilde{b}_2\mathcal{Z}^\alpha\partial_y^{2} v_1\cdot\mathcal{Z}^\alpha\partial_y b_1\ d{\bf{x}}ds.
\end{split}
\end{equation}
Where the first term on the right hand side of \eqref{4.48} is estimated as follows due to \eqref{2.3}, \eqref{2.4} and the fact of $\nabla\cdot {\bf{B}}=0$.
\begin{align*}
&2\varepsilon\mu\sum_{\beta+\gamma=\alpha\atop |\beta|\geq 1} C_\alpha^\beta\int_0^t\int_{\mathbb{R}_+^2} \mathcal{Z}^\beta\tilde{b}_2\mathcal{Z}^\gamma\partial_y^{2} v_1\cdot\mathcal{Z}^\alpha\partial_y b_1\ d{\bf{x}}ds\\
\lesssim&\;\varepsilon\mu\sum_{\beta+\gamma=\alpha\atop |\beta|\geq |\gamma|}\|\phi\mathcal{Z}^\gamma \partial_y^2 v_1\|_{L_{t, {\bf{x}}}^\infty}\left(\int_0^t\int_\Omega|\phi^{-1}\mathcal{Z}^\beta\tilde b_2|^2 d{\bf{x}}ds\right)^\frac12\left(\int_0^t\|\partial_y b_1\|_{m-1}^2 ds\right)^\frac12\\
&+\varepsilon\mu\sum_{\beta+\gamma=\alpha\atop |\beta|< |\gamma|<|\alpha|}\|\phi^{-1}\mathcal{Z}^\beta\tilde b_2\|_{L_{t, {\bf{x}}}^\infty}\left(\int_0^t\int_\Omega|\phi\mathcal{Z}^\gamma \partial_y^2 v_1|^2 d{\bf{x}}ds\right)^\frac12\left(\int_0^t\|\partial_y b_1\|_{m-1}^2 ds\right)^\frac12\\
\lesssim&\;\varepsilon\mu\|\partial_y v_1\|_{[(m-1)/2]+1, \infty}\left(\int_0^t\|b_1\|_{m}^2 ds\right)^\frac12\left(\int_0^t\|\partial_y b_1\|_{m-1}^2 ds\right)^\frac12\\
&+\varepsilon\mu\|  b_1\|_{[(m-1)/2]+1, \infty}\left(\int_0^t\|\partial_y v_1\|_{m-1}^2 ds\right)^\frac12\left(\int_0^t\|\partial_y b_1\|_{m-1}^2 ds\right)^\frac12,
\end{align*}
here $\tilde{b}_2|_{y=0}=0$ is used.\\
For the second term on the right hand side of \eqref{4.48}, it satisfies that
\begin{align*}
&2\varepsilon\mu\int_0^t\int_{\mathbb{R}_+^2} \tilde{b}_2\mathcal{Z}^\alpha\partial_y^{2} v_1\cdot\mathcal{Z}^\alpha\partial_y b_1\ d{\bf{x}}ds\\
\leq&\;2\varepsilon\mu\|\phi^{-1}\tilde{b}_2\|_{L_{t, {\bf{x}}}^\infty}
\left(\int_0^t\|\partial_y^2 v_1\|_{m-1}^2 ds\right)^\frac12\left(\int_0^t\int_\Omega|\phi\mathcal{Z}^\alpha \partial_y b_1|^2d{\bf{x}}ds\right)^\frac12\\
\leq &\;\frac{\varepsilon^2\mu^2}{2}\int_0^t\|\partial_y^2 v_1\|_{m-1}^2 ds+2 \|b_1\|_{1, \infty}^2\int_0^t\|b_1\|_{m}^2 ds.
\end{align*}
Third,
\begin{equation}\label{4.17}
\begin{split}
&-2\varepsilon\mu\int_0^t\int_{\mathbb{R}_+^2} \mathcal{Z}^\alpha\partial_y(v_2\partial_y b_1)\cdot\mathcal{Z}^\alpha\partial_y b_1\ d{\bf{x}}ds\\
=&-2\varepsilon\mu\int_0^t\int_{\mathbb{R}_+^2} \mathcal{Z}^\alpha(\partial_yv_2\partial_y b_1)\cdot\mathcal{Z}^\alpha\partial_y b_1\ d{\bf{x}}ds\\
&-2\varepsilon\mu\sum_{\beta+\gamma=\alpha\atop |\beta|\geq1}\int_0^t\int_{\mathbb{R}_+^2} \mathcal{Z}^\beta v_2\mathcal{Z}^\gamma\partial_y^2 b_1\cdot\mathcal{Z}^\alpha\partial_y b_1\ d{\bf{x}}ds\\
&-2\varepsilon\mu \int_0^t\int_{\mathbb{R}_+^2}  v_2\mathcal{Z}^\alpha\partial_y^2 b_1\cdot\mathcal{Z}^\alpha\partial_y b_1\ d{\bf{x}}ds.
\end{split}
\end{equation}
For the first term on the right hand side of \eqref{4.17},  by \eqref{2.1}, we have
\begin{equation*}
\begin{split}
&-2\varepsilon\mu\int_0^t\int_{\mathbb{R}_+^2} \mathcal{Z}^\alpha(\partial_yv_2\partial_y b_1)\cdot\mathcal{Z}^\alpha\partial_y b_1\ d{\bf{x}}ds\\
\lesssim&\;\varepsilon\mu \|\partial_yb_1\|_{ L_{t, {\bf{x}}}^\infty}\left(\int_0^t\|\partial_y  v_2\|_{m-1}^2 ds\right)^\frac12\left(\int_0^t\|\partial_y b_1\|_{m-1}^2 ds\right)^\frac12+\varepsilon\mu \|\partial_yv_2\|_{  L_{t, {\bf{x}}}^\infty}\int_0^t\|\partial_y  b_1\|_{m-1}^2 ds.
\end{split}
\end{equation*}
For the second term on the right hand side of \eqref{4.17}, by \eqref{2.3} and \eqref{2.4}, it holds that
\begin{align*}
&-2\varepsilon\mu\sum_{\beta+\gamma=\alpha\atop |\beta|\geq1}\int_0^t\int_{\mathbb{R}_+^2} \mathcal{Z}^\beta v_2\mathcal{Z}^\gamma\partial_y^2 b_1\cdot\mathcal{Z}^\alpha\partial_y b_1\ d{\bf{x}}ds\\
\lesssim&\;\varepsilon\mu\sum_{\beta+\gamma=\alpha\atop|\beta|\geq|\gamma|}
\|\phi\mathcal{Z}^\gamma\partial_y^2 b_1 \|_{L_{t, {\bf{x}}}^\infty}\left(\int_0^t\int_ {\mathbb{R}_+^2}|\phi^{-1}\mathcal{Z}^\beta v_2|^2 d{\bf{x}} ds\right)^\frac12\left(\int_0^t\|\partial_y b_1\|_{m-1}^2 ds\right)^\frac12\\
&+\varepsilon\mu\sum_{\beta+\gamma=\alpha\atop|\beta|<|\gamma|<|\alpha|}
\|\phi^{-1}\mathcal{Z}^\beta v_2\|_{L_{t, {\bf{x}}}^\infty}\left(\int_0^t\int_{\mathbb{R}_+^2}|\phi\mathcal{Z}^\gamma\partial_y^2 b_1|^2 d{\bf{x}}ds\right)^\frac12\left(\int_0^t\|\partial_y b_1\|_{m-1}^2 ds\right)^\frac12\\
\lesssim&\;\varepsilon\mu \|\partial_yb_1\|_{[(m-1)/2]+1,\infty}\left(\int_0^t\|\partial_y  v_2\|_{m-1}^2 ds\right)^\frac12\left(\int_0^t\|\partial_y b_1\|_{m-1}^2 ds\right)^\frac12\\
&+\varepsilon\mu \|\partial_yv_2\|_{ [(m-1)/2],\infty}\int_0^t\|\partial_y  b_1\|_{m-1}^2 ds.
\end{align*}
For the last term on the right hand side of \eqref{4.17}, by \eqref{2.2} and integration by parts,  we have
\begin{align*}
&-2\varepsilon\mu \int_0^t\int_{\mathbb{R}_+^2}  v_2\mathcal{Z}^\alpha\partial_y^2 b_1\cdot\mathcal{Z}^\alpha\partial_y b_1\ d{\bf{x}}ds\\
=&-2\varepsilon\mu \sum_{k=0}^{m-2}\int_0^t\int_{\mathbb{R}_+^2} \phi^{k, m-1}(y) v_2\partial_y\mathcal{Z}_2^k\partial_y b_1\cdot\mathcal{Z}^\alpha\partial_y b_1\ d{\bf{x}}ds\\
&-2\varepsilon\mu \int_0^t\int_{\mathbb{R}_+^2}  v_2\partial_y\mathcal{Z}^\alpha\partial_y b_1\cdot\mathcal{Z}^\alpha\partial_y b_1\ d{\bf{x}}ds\\
\lesssim&\;\varepsilon\mu\sum_{k=0}^{m-2}\|\phi^{-1}v_2\|_{L_{t, {\bf{x}}}^\infty}\left(\int_0^t\int_{\mathbb{R}_+^2}|\phi\partial_y\mathcal{Z}_2^k\partial_y b_1|^2 d{\bf{x}}ds\right)^\frac12\left(\int_0^t\|\partial_y b_1\|_{m-1}^2 ds\right)^\frac12\\
&+\varepsilon\mu\|\partial_yv_2\|_{L_{t, {\bf{x}}}^\infty}\int_0^t\|\partial_y b_1\|_{m-1}^2 ds\\
\lesssim&\;\varepsilon\mu\|\partial_yv_2\|_{L_{t, {\bf{x}}}^\infty}\int_0^t\|\partial_y b_1\|_{m-1}^2 ds.
\end{align*}
We continue to handle the fourth term on the right hand side of \eqref{4.15}.
\begin{equation}\label{4.18}
\begin{split}
&-2\varepsilon\mu\int_0^t\int_{\mathbb{R}_+^2} \mathcal{Z}^\alpha\partial_y(b_1\partial_y v_2)\cdot\mathcal{Z}^\alpha\partial_y b_1\ d{\bf{x}}ds\\
&=-2\varepsilon\mu\int_0^t\int_{\mathbb{R}_+^2} \mathcal{Z}^\alpha(\partial_yb_1\partial_y v_2)\cdot\mathcal{Z}^\alpha\partial_y b_1\ d{\bf{x}}ds\\
&-2\varepsilon\mu\int_0^t\int_{\mathbb{R}_+^2} \mathcal{Z}^\alpha( b_1\partial_y^2 v_2)\cdot\mathcal{Z}^\alpha\partial_y b_1\ d{\bf{x}}ds.
\end{split}
\end{equation}
For the first term on the right hand side of \eqref{4.18},  by \eqref{2.1}, we have
\begin{align*}
&-2\varepsilon\mu\int_0^t\int_{\mathbb{R}_+^2} \mathcal{Z}^\alpha(\partial_yb_1\partial_y v_2)\cdot\mathcal{Z}^\alpha\partial_y b_1\ d{\bf{x}}ds\\
\lesssim&\;\varepsilon\mu\|\partial_y b_1\|_{L_{t, {\bf{x}}}^\infty}\left(\int_0^t\|\partial_y v_2\|_{m-1}^2 ds\right)^\frac12\left(\int_0^t\|\partial_y b_1\|_{m-1}^2 ds\right)^\frac12\\
&+\varepsilon\mu\|\partial_y v_2\|_{L_{t, {\bf{x}}}^\infty} \int_0^t\|\partial_y b_1\|_{m-1}^2 ds.
\end{align*}
For the second term on the right hand side of \eqref{4.18}, one has
\begin{equation}\label{4.19}
\begin{split}
&-2\varepsilon\mu\int_0^t\int_{\mathbb{R}_+^2} \mathcal{Z}^\alpha( b_1\partial_y^2 v_2)\cdot\mathcal{Z}^\alpha\partial_y b_1\ d{\bf{x}}ds\\
=&-2\varepsilon\mu\sum_{\beta+\gamma=\alpha\atop |\beta|\geq1}\int_0^t\int_{\Omega} \mathcal{Z}^\beta b_1\mathcal{Z}^\gamma\partial_y^2 v_2\cdot\mathcal{Z}^\alpha\partial_y b_1\ d{\bf{x}}ds-2\varepsilon\mu\int_0^t\int_{\Omega}b_1 \mathcal{Z}^\alpha\partial_y^2 v_2\cdot\mathcal{Z}^\alpha\partial_y b_1\ d{\bf{x}}ds.
\end{split}
\end{equation}
By the Sobolev embedding inequality and \eqref{4.3}, the first term on the right hand side of \eqref{4.19} can be estimated by following arguments.
\begin{equation}\label{ad1}
\begin{split}
&-2\varepsilon\mu\sum_{\beta+\gamma=\alpha\atop |\beta|\geq1}\int_0^t\int_{\Omega} \mathcal{Z}^\beta b_1\mathcal{Z}^\gamma\partial_y^2 v_2\cdot\mathcal{Z}^\alpha\partial_y b_1\ d{\bf{x}}ds\\
=&-2\varepsilon\mu\sum_{\beta+\gamma=\alpha\atop |\gamma|<|\beta|\leq |\alpha|}\int_0^t\int_{\Omega} \mathcal{Z}^\beta b_1\mathcal{Z}^\gamma\partial_y^2 v_2\cdot\mathcal{Z}^\alpha\partial_y b_1\ d{\bf{x}}ds\\
&-2\varepsilon\mu\sum_{\beta+\gamma=\alpha\atop 1\leq|\beta|\leq |\gamma|}\int_0^t\int_{\Omega} \mathcal{Z}^\beta b_1\mathcal{Z}^\gamma\partial_y^2 v_2\cdot\mathcal{Z}^\alpha\partial_y b_1\ d{\bf{x}}ds.
\end{split}
\end{equation}
For the first term on the right hand side of \eqref{ad1}, we have
\begin{align*}
&-2\varepsilon\mu\sum_{\beta+\gamma=\alpha\atop |\gamma|<|\beta|\leq |\alpha|}\int_0^t\int_{\Omega} \mathcal{Z}^\beta b_1\mathcal{Z}^\gamma\partial_y^2 v_2\cdot\mathcal{Z}^\alpha\partial_y b_1\ d{\bf{x}}ds\\
\lesssim&\; \varepsilon\mu\sum_{\beta+\gamma=\alpha\atop |\beta|\geq|\gamma|}\sup_{0\leq s\leq t}\|\mathcal{Z}^\gamma \partial_y^2 v_2(s)\|_{L_{x}^{\infty}(L_y^2)}\left(\int_0^t\|\mathcal{Z}^\beta b_1\|_{L_{x}^{2}(L_y^\infty)}^2ds\right)^\frac12\left(\int_0^t\|\partial_y b_1\|_{m-1}^2 ds\right)^\frac12\\
\lesssim&\;\varepsilon\mu\left[\|\partial_y^2 v_2(0)\|_{[(m-1)/2]+1}+\left(\int_0^t\|\partial_y^2 v_2\|_{[(m-1)/2]+2}^2ds\right)^\frac12\right]\left(\int_0^t\|(b_1,\partial_y b_1)\|_{m-1}^2 ds\right)^\frac12\\
& \cdot\left(\int_0^t\|\partial_y b_1\|_{m-1}^2 ds\right)^\frac12.
\end{align*}
For the second term on the right hand side of \eqref{ad1},  by \eqref{4.3}, we have
\begin{align*}
&-2\varepsilon\mu\sum_{\beta+\gamma=\alpha\atop 1\leq|\beta|\leq |\gamma|}\int_0^t\int_{\Omega} \mathcal{Z}^\beta b_1\mathcal{Z}^\gamma\partial_y^2 v_2\cdot\mathcal{Z}^\alpha\partial_y b_1\ d{\bf{x}}ds\\
\lesssim&\;\varepsilon\mu\sum_{\beta+\gamma=\alpha\atop |\beta|<|\gamma|<|\alpha|}\|\mathcal{Z}^\beta b_1\|_{L_{t, {\bf{x}}}^\infty}\left(\int_0^t\int_{\mathbb{R}_+^2}|\mathcal{Z}^\gamma \partial_y^2 v_2|^2 d{\bf{x}}ds\right)^\frac12\left(\int_0^t\|\partial_y b_1\|_{m-1}^2 ds\right)^\frac12\\
\lesssim&\;\varepsilon\mu\sum_{\beta+\gamma=\alpha\atop |\beta|<|\gamma|<|\alpha|}\|\mathcal{Z}^\beta b_1\|_{L_{t, {\bf{x}}}^\infty}\left[\left(\int_0^t\int_{\mathbb{R}_+^2}|\mathcal{Z}^\gamma \partial_y\partial_x v_1|^2 d{\bf{x}}ds\right)^\frac12 + \left(\int_0^t\int_{\mathbb{R}_+^2}|\mathcal{Z}^\gamma \partial_y (p^{-1} \partial_t p)|^2 d{\bf{x}}ds\right)^\frac12\right. \\
&\left.+ \left(\int_0^t\int_{\mathbb{R}_+^2}|\mathcal{Z}^\gamma \partial_y (p^{-1} {\bf{v}}\cdot\nabla p)|^2 d{\bf{x}}ds\right)^\frac12\right]\left(\int_0^t\|\partial_y b_1\|_{m-1}^2 ds\right)^\frac12\\
\lesssim &\;\varepsilon\mu \|b_1\|_{[(m-1)/2], \infty}\left(1+\sum_{i=0}^1\|\partial_y^i({\bf{v}}, p, p^{-1})\|_{1, \infty}\right)^2 \sum_{j=0}^1\left(\int_0^t \|\partial_y^j({\bf{v}}, p-1, p^{-1}-1)\|_{m-1}^2 ds \right)^\frac12\\
&\cdot\left(\int_0^t\|\partial_y b_1\|_{m-1}^2 ds\right)^\frac12.
\end{align*}
As for the second term on the right hand side of \eqref{4.19}, it follows from  Young's inequality that
\begin{align*}
-2\varepsilon\mu\int_0^t\int_{\mathbb{R}_+^2}b_1 \mathcal{Z}^\alpha\partial_y^2 v_2\cdot\mathcal{Z}^\alpha\partial_y b_1\ d{\bf{x}}ds
\leq 2\varepsilon^2\mu^2 \|b_1\|_{L_{t, {\bf{x}}}^\infty}^2\int_0^t\| \partial_y^2 v_2\|_{m-1}^2  ds+\frac12\int_0^t\|\partial_y b_1\|_{m-1}^2 ds.
\end{align*}
Finally, the last term on the right hand side of \eqref{4.15} can be estimated as follows.
\begin{equation}\label{4.20}
\begin{split}
&-2\varepsilon\mu\int_0^t\int_{\mathbb{R}_+^2} \mathcal{Z}^\alpha\partial_y(v_1\partial_x b_1)\cdot\mathcal{Z}^\alpha\partial_y b_1\ d{\bf{x}}ds\\
=&-2\varepsilon\mu\int_0^t\int_{\mathbb{R}_+^2} \mathcal{Z}^\alpha(\partial_yv_1\partial_x b_1)\cdot\mathcal{Z}^\alpha\partial_y b_1\ d{\bf{x}}ds\\
&-2\varepsilon\mu\int_0^t\int_{\mathbb{R}_+^2} \mathcal{Z}^\alpha( v_1\partial_x\partial_y b_1)\cdot\mathcal{Z}^\alpha\partial_y b_1\ d{\bf{x}}ds.
\end{split}
\end{equation}
By \eqref{2.1}, for the first term on the right hand side of \eqref{4.20}, it holds that
\begin{align*}
&-2\varepsilon\mu\int_0^t\int_{\mathbb{R}_+^2} \mathcal{Z}^\alpha(\partial_yv_1\partial_x b_1)\cdot\mathcal{Z}^\alpha\partial_y b_1\ d{\bf{x}}ds\\
\lesssim&\;\varepsilon\mu\|\partial_y v_1\|_{L_{t, {\bf{x}}}^\infty}\left(\int_0^t\|b_1\|_{m}^2 ds\right)^\frac12\left(\int_0^t\|\partial_y b_1\|_{m-1}^2 ds\right)^\frac12\\
&+\varepsilon\mu\|b_1\|_{ 1, \infty}\left(\int_0^t\|\partial_y v_1\|_{m-1}^2 ds\right)^\frac12\left(\int_0^t\|\partial_y b_1\|_{m-1}^2 ds\right)^\frac12.
\end{align*}
And for the second term on the right hand side of \eqref{4.20},  one has
\begin{equation}\label{4.21}
\begin{split}
&-2\varepsilon\mu\int_0^t\int_{\mathbb{R}_+^2} \mathcal{Z}^\alpha( v_1\partial_x\partial_y b_1)\cdot\mathcal{Z}^\alpha\partial_y b_1\ d{\bf{x}}ds\\
=&-2\varepsilon\mu\sum_{\beta+\gamma=\alpha\atop |\beta|\geq 1}\int_0^t\int_{\mathbb{R}_+^2} \mathcal{Z}^\beta v_1\mathcal{Z}^\gamma\partial_x\partial_y b_1\cdot\mathcal{Z}^\alpha\partial_y b_1\ d{\bf{x}}ds\\
&-2\varepsilon\mu\int_0^t\int_{\mathbb{R}_+^2} v_1\partial_x\mathcal{Z}^\alpha\partial_y b_1\cdot\mathcal{Z}^\alpha\partial_y b_1\ d{\bf{x}}ds.
\end{split}
\end{equation}
Then, the first term on the right hand side of \eqref{4.21} satisfies the following inequality.
\begin{align*}
&-2\varepsilon\mu\sum_{\beta+\gamma=\alpha\atop |\beta|\geq 1}\int_0^t\int_{\mathbb{R}_+^2} \mathcal{Z}^\beta v_1\mathcal{Z}^\gamma\partial_x\partial_y b_1\cdot\mathcal{Z}^\alpha\partial_y b_1\ d{\bf{x}}ds\\
 \lesssim&\;\varepsilon\mu\| v_1\|_{ 1, \infty} \int_0^t\|\partial_y b_1\|_{m-1}^2 ds  +\varepsilon\mu\| \partial_y b_1 \|_{1, \infty}\left(\int_0^t\| v_1\|_{m-1}^2 ds\right)^\frac12\left(\int_0^t\|\partial_y b_1\|_{m-1}^2 ds\right)^\frac12.
\end{align*}
For the second term on the right hand side of \eqref{4.21}, by integration by parts, it is direct to obtain
\begin{align*}
-2\varepsilon\mu\int_0^t\int_{\mathbb{R}_+^2} v_1\partial_x\mathcal{Z}^\alpha\partial_y b_1\cdot\mathcal{Z}^\alpha\partial_y b_1\ d{\bf{x}}ds
\lesssim\varepsilon\mu\|v_1\|_{1, \infty}\int_0^t\|\partial_y b_1\|_{m-1}^2 ds.
\end{align*}
Now it is left to estimate the terms on the right hand side of \eqref{4.14}. For the first three terms on the right hand side of \eqref{4.14},  by \eqref{2.1}, we have
\begin{align*}
&	\int_0^t\int_{\mathbb{R}_+^2} |\mathcal{Z}^\alpha (\rho\partial_t v_1)|^2 d{\bf{x}}ds+	\int_0^t\int_{\mathbb{R}_+^2} |\mathcal{Z}^\alpha \partial_x p|^2 d{\bf{x}}ds+\int_0^t\int_{\mathbb{R}_+^2} |\mathcal{Z}^\alpha (b_2\partial_x b_2)|^2 d{\bf{x}}ds\\
\lesssim&\;\|(\rho, v_1, b_2)\|_{1, \infty}^2 \int_0^t\|(\rho-1, v_1, \tilde{b}_2)\|_{m}^2 ds  +\int_0^t\|p-1\|_{m}^2 ds.
\end{align*}
For the fourth term on the right hand side of \eqref{4.14},  by \eqref{2.1}, we get
\begin{align*}
&\int_0^t\int_{\mathbb{R}_+^2} |\mathcal{Z}^\alpha (\tilde b_2\partial_y b_1)|^2 d{\bf{x}}ds \\
\lesssim&\; \|\phi\partial_y b_1\|_{L_{t, {\bf{x}}}^\infty}^2\int_0^t\|\phi^{-1}\tilde b_2\|_{m-1}^2 ds+\|\phi^{-1}\tilde b_2\|_{L_{t, {\bf{x}}}^\infty}^2\int_0^t\|\phi\partial_y b_1\|_{m-1}^2 ds\\
\lesssim&\;
\|b_1\|_{1, \infty}^2\int_0^t\|b_1\|_{m}^2 ds.
\end{align*}
For the fifth term on the right hand side of \eqref{4.14}, by \eqref{2.1}, \eqref{2.3} and \eqref{2.4},  we have
\begin{align*}
&\int_0^t\int_{\mathbb{R}_+^2} |\mathcal{Z}^\alpha (\rho {\bf{v}}\cdot\nabla v_1)|^2 d{\bf{x}}ds \\
\lesssim&\;\|(\rho, v_1)\|_{1,\infty}^4 \int_0^t \|(\rho-1, v_1)\|_{m}^2 ds\\\
&+\sum_{\beta+\iota+\gamma=\alpha\atop |\beta|\geq|\gamma|, |\iota|}\|\mathcal{Z}^\gamma v_2\|_{L_{t, {\bf{x}}}^\infty}^2\|\mathcal{Z}^\iota \partial_yv_1\|_{L_{t, {\bf{x}}}^\infty}^2 \int_0^t\int_\Omega|\mathcal{Z}^\beta\rho|^2d{\bf{x}}ds \\
&+\sum_{\beta+\iota+\gamma=\alpha\atop |\gamma|\geq|\beta|, |\iota|}\|\mathcal{Z}^\beta\rho\|_{L_{t, {\bf{x}}}^\infty}^2\|\mathcal{Z}^\iota \partial_yv_1\|_{L_{t, {\bf{x}}}^\infty}^2 \int_0^t\int_\Omega|\mathcal{Z}^\gamma v_2|^2d{\bf{x}}ds \\
&+\sum_{\beta+\iota+\gamma=\alpha\atop |\iota|\geq|\beta|,|\gamma| }\|\mathcal{Z}^\beta\rho\|_{L_{t, {\bf{x}}}^\infty}^2\|\phi^{-1}\mathcal{Z}^\gamma v_2\|_{L_{t, {\bf{x}}}^\infty}^2 \int_0^t\int_\Omega|\phi\mathcal{Z}^\iota \partial_yv_1|^2d{\bf{x}}ds \\
\lesssim&\; \|(\rho, {\bf{v}})\|_{[(m-1)/2], \infty}^2\|(\rho, {\bf{v}}, \partial_y {\bf{v}})\|_{[(m-1)/2], \infty}^2\int_0^t\|(\rho-1, {\bf{v}})\|_{m}^2 ds.
\end{align*}
The last two terms on the right hand side of \eqref{4.14} can be estimated directly.
\begin{align*}
&\varepsilon^2\mu^2\int_0^t\int_{\mathbb{R}_+^2}|\mathcal{Z}^\alpha \partial_x^2 v_1|^2\ d{\bf{x}}ds+\varepsilon^2(\mu+\lambda)^2\int_0^t\int_{\mathbb{R}_+^2} |\mathcal{Z}^\alpha \partial_x(\nabla\cdot {\bf{v}})|^2\  d{\bf{x}}ds\\
\lesssim&\; \varepsilon^2\mu^2\int_0^t\|\partial_xv_1\|_{m}^2 ds+\varepsilon^2(\mu+\lambda)^2\int_0^t\|\nabla\cdot {\bf{v}}\|_{m}^2 ds.
\end{align*}
Collecting all of the estimates in this subsection leads to the following inequality.
\begin{align*}
&\int_0^t\left(\| \partial_y b_1\|_{m-1}^2 +\varepsilon^2\mu^2\|\partial_y^2 v_1\|_{m-1}^2 \right)ds  +\varepsilon\mu\| \partial_y b_1(t)\|_{m-1}^2\\
\lesssim&\;\varepsilon\mu\|\partial_y b_1(0)\|_{m-1}^2   +\varepsilon^2\mu^2\int_0^t\|\nabla {\bf{v}}\|_{m}^2 ds+\varepsilon^2(\mu+\lambda)^2\int_0^t\|\nabla\cdot {\bf{v}}\|_{m}^2 ds\\
&+\left(1+\|(\rho,{\bf{v}}, {\bf{B}}, \partial_y {\bf{v}})\|_{[(m-1)/2], \infty}^2\right)^2\int_0^t\|(\rho-1, {\bf{v}},  {\bf{B}}-\overrightarrow{e_y}, p-1)\|_{m}^2ds\\
&+\varepsilon\mu \left[\left(1+\sum_{i=0}^1\|\partial_y^i({\bf{v}}, b_1, p, p^{-1})\|_{[(m-1)/2]+1, \infty}^2\right)^2+\|\partial_y^2 v_2(0)\|_{[(m-1)/2]+1}^2\right.\\
&\left.+\int_0^t\|\partial_y^2 v_2\|_{[(m-1)/2]+2}^2 ds\right] \int_0^t\left( \| ({\bf{v}}, b_1, p-1, p^{-1}-1)\|_{m}^2 +\|\partial_y ({\bf{v}}, b_1)\|_{m-1}^2\right)ds\\
&+ \varepsilon^2\mu^2\|b_1\|_{L_{t,{\bf{x}}}^\infty}^2
\int_0^t\| \partial_y^2 v_2\|_{m-1}^2  ds.
\end{align*}

Thus, the proof of Lemma \ref{lem2} is done by combining the conormal estimates of $\partial_y{\bf{v}}, \partial_yb_1$ and $ \partial_yp$ together. It is remarked that the smallness of $\|b_1\|_{L_{t,{\bf{x}}}^\infty}$ is required here.

\subsection{The second order  normal derivative}
In this subsection, we will establish the conormal estimates for the second order normal derivative of the classical solutions to compressible MHD equations \eqref{3.1} with no-slip boundary condition \eqref{1.3}.
\begin{lemma}\label{lem3}
	Under the assumption in Theorem \ref{Th1}, for any $|\alpha|\leq m-2$, there exists a sufficiently small $\varepsilon_0>0$ such that for any $\varepsilon\in(0, \varepsilon_0)$, the classical solutions $(\rho, {\bf{v}}, {\bf{B}})$ to compressible  MHD equations \eqref{3.1} with no-slip boundary condition \eqref{1.3} satisfies
\begin{align*}
	&\int_0^t\|\partial_y^2({\bf{v}}, b_1, p)\|_{m-2}^2 ds+ \varepsilon^2\mu^2 \int_0^t\|\partial_y^3 v_1\|_{m-2}^2ds+\varepsilon^2(2\mu+\lambda)^2\int_0^t\|\partial_y^3  v_2\|_{m-2}^2 ds\\
	&+\varepsilon(2\mu+\lambda)\gamma^{-1}\sum_{|\alpha|\leq m-2}\int_{\mathbb{R}_+^2} p^{-1}(t)|\mathcal{Z}^\alpha \partial_y^2 p(t)|^2 d{\bf{x}}+\varepsilon\mu\|\partial_y^2 b_1(t)\|_{m-2}^2  \\
	\lesssim&\;\varepsilon(2\mu+\lambda)\gamma^{-1}\sum_{|\alpha|\leq m-2}\int_{\mathbb{R}_+^2} p_0^{-1}|\mathcal{Z}^\alpha \partial_y^2 p_0|^2 d{\bf{x}} +\varepsilon\mu\|\partial_y^2 b_1(0)\|_{m-2}^2+\varepsilon^2(2\mu+\lambda)^2\int_0^t\left(\|\partial_y {\bf{v}}\|_{m}^2\right. \\
	&+\left.\|\partial_y^2 {\bf{v}}\|_{m-1}^2\right)  ds+\left(\|\partial_y^2 b_1(0)\|_{[m/2]}^2+\int_0^t\|\partial_y^2 b_1 \|_{[m/2]+1}^2\right)\int_0^t\|(b_1, \partial_y b_1)\|_{m-2}^2 ds\\
&+\left(1+\sum_{i=0}^1\|\partial_y^i(\rho, {\bf{v}},  {\bf{B}}, p, p^{-1})\|_{[m/2], \infty}^2\right)^2\cdot\sum_{j=0}^1\int_0^t\|\partial_y^j(\rho-1, {\bf{v}},  {\bf{B}}-\overset{\rightarrow}{e_y}, p-1, p^{-1}-1)\|_{m-j}^2  ds\\
	& +\varepsilon(2\mu+\lambda)\left[ \left(1+\sum_{i=0}^1\|\partial_y^i({\bf{v}}, b_1,p, p^{-1} )\|_{[m/2]+1, \infty}^2\right)^2 +\sum_{i=1}^2\|\partial_y^i({\bf{v}}, b_1, p, p^{-1})(0)\|_{[m/2]+2}^2\right. \\
	&\left. +\sum_{i=1}^2\int_0^t\|\partial_y^i({\bf{v}}, b_1, p-1, p^{-1}-1)\|_{[m/2]+3}^2ds\right]\cdot\sum_{j=0}^2 \int_0^t \|\partial_y^j({\bf{v}}, b_1, p-1, p^{-1}-1)\|_{m-j}^2 ds.	\end{align*}
\end{lemma}

\subsubsection{Conormal Estimate of $  \partial_y^2v_1$} By similar arguments as in \eqref{4.2},  we have
\begin{equation}\label{4.22}
\begin{split}
\int_0^t\int_{\mathbb{R}_+^2}| \mathcal{Z}^\alpha \partial_y^2 v_1|^2\ d{\bf{x}}ds
&\lesssim \int_0^t\int_{\mathbb{R}_+^2}|	\partial_t\mathcal{Z}^\alpha \partial_y b_1|^2 d{\bf{x}}ds+\int_0^t\int_{\mathbb{R}_+^2}|\mathcal{Z}^\alpha \partial_y(\tilde{b}_2\partial_y v_1)|^2 d{\bf{x}}ds\\
&+\int_0^t\int_{\mathbb{R}_+^2}|\mathcal{Z}^\alpha \partial_y (v_2\partial_y b_1 )|^2 d{\bf{x}}ds+\int_0^t\int_{\mathbb{R}_+^2}|\mathcal{Z}^\alpha \partial_y (v_1\partial_x b_1)|^2 d{\bf{x}}ds\\
&+\int_0^t\int_{\mathbb{R}_+^2}|\mathcal{Z}^\alpha \partial_y (b_1\partial_y v_2)|^2 d{\bf{x}}ds.
\end{split}
\end{equation}
It is direct to estimate the first term on the right hand side of \eqref{4.22} as
\begin{align*}
\int_0^t\int_{\mathbb{R}_+^2}|	\partial_t\mathcal{Z}^\alpha \partial_y b_1|^2 \ d{\bf{x}}ds
\lesssim\int_0^t\|\partial_y b_1\|_{m-1}^2 ds.
\end{align*}
By \eqref{2.1} and the fact of $\nabla\cdot {\bf{B}}=0$, the second term on the right hand side of \eqref{4.22} satisfies
\begin{equation*}
\begin{split}
&\int_0^t\int_{\mathbb{R}_+^2}|\mathcal{Z}^\alpha \partial_y (\tilde{b}_2\partial_y v_1)|^2 \ d{\bf{x}}ds\\
\lesssim&\int_0^t\int_{\mathbb{R}_+^2}|\mathcal{Z}^\alpha (\partial_y \tilde{b}_2\partial_y  v_1)|^2 \ d{\bf{x}}ds+\int_0^t\int_{\mathbb{R}_+^2}|\mathcal{Z}^\alpha (\tilde{b}_2\partial_y^2 v_1)|^2 \ d{\bf{x}}ds\\
\lesssim&\;\|\partial_x b_1\|_{L_{t, {\bf{x}}}^\infty}^2\int_0^t\|\partial_y v_1\|_{m-2}^2 ds+\|\partial_y v_1\|_{L_{t, {\bf{x}}}^\infty}^2\int_0^t\|\partial_x b_1\|_{m-2}^2 ds\\
& +\|\phi^{-1}\tilde{b}_2\|_{L_{t, {\bf{x}}}^ \infty}^2\int_0^t\|\phi\partial_y^2v_1\|_{m-2}^2ds+\|\phi\partial_y^2v_1\|_{L_{t, {\bf{x}}}^ \infty}^2\int_0^t\|\phi^{-1}\tilde{b}_2\|_{m-2}^2ds\\
\lesssim&\;\|b_1\|_{1, \infty}^2\int_0^t\|\partial_y v_1\|_{m-1}^2 ds +\|\partial_yv_1\|_{1, \infty}^2\int_0^t\| b_1\|_{m-1}^2ds.
\end{split}
\end{equation*}
Similarly, for the third term on the right hand side of \eqref{4.22}, we have
\begin{align*}
&\int_0^t\int_{\mathbb{R}_+^2}|\mathcal{Z}^\alpha \partial_y(v_2\partial_y b_1 )|^2 d{\bf{x}}ds\\
\lesssim&\int_0^t\int_{\Omega}|\mathcal{Z}^\alpha (\partial_yv_2\partial_y b_1 )|^2 d{\bf{x}}ds+\int_0^t\int_{\Omega}|\mathcal{Z}^\alpha ( v_2\partial_y^2 b_1 )|^2 d{\bf{x}}ds\\
\lesssim&\;\|\partial_y v_2\|_{L_{t, {\bf{x}}}^\infty}^2\int_0^t\|\partial_y b_1\|_{m-2}^2 ds+\|\partial_y b_1\|_{L_{t, {\bf{x}}}^\infty}^2\int_0^t\|\partial_y v_2\|_{m-2}^2 ds\\
&+\|\phi^{-1}v_2\|_{L_{t, {\bf{x}}}^\infty}^2\int_0^t\|\phi\partial_y^2 b_1\|_{m-2}^2 ds+\|\phi\partial_y^2 b_1\|_{L_{t, {\bf{x}}}^\infty}^2\int_0^t\|\phi^{-1}v_2\|_{m-2}^2 ds\\
\lesssim&\;\|\partial_y v_2\|_{L_{t, {\bf{x}}}^\infty}^2\int_0^t\|\partial_y b_1\|_{m-1}^2 ds+\|\partial_y b_1\|_{1, \infty}^2\int_0^t\|\partial_y v_2\|_{m-2}^2 ds.
\end{align*}
And the following inequality holds for the fourth term on the right hand side of \eqref{4.22}.
\begin{align*}
&\int_0^t\int_{\mathbb{R}_+^2}|\mathcal{Z}^\alpha \partial_y (v_1\partial_x b_1)|^2 d{\bf{x}}ds\\
\lesssim & \int_0^t\int_{\mathbb{R}_+^2}|\mathcal{Z}^\alpha  (\partial_yv_1\partial_x b_1)|^2  d{\bf{x}}ds+\int_0^t\int_{\mathbb{R}_+^2}|\mathcal{Z}^\alpha   (v_1\partial_x \partial_yb_1)|^2  d{\bf{x}}ds\\
\lesssim&\;\|\partial_y v_1\|_{L_{t, {\bf{x}}}^\infty}^2\int_0^t\|\partial_x b_1\|_{m-2}^2 ds+\|\partial_x b_1\|_{L_{t, {\bf{x}}}^\infty}^2\int_0^t\|\partial_y v_1\|_{m-2}^2 ds\\\
&+\| v_1\|_{L_{t, {\bf{x}}}^\infty}^2\int_0^t\|\partial_x\partial_y b_1\|_{m-2}^2 ds+\|\partial_x\partial_y b_1\|_{L_{t, {\bf{x}}}^\infty}^2\int_0^t\| v_1\|_{m-2}^2 ds\\
\lesssim&\;\|\partial_y v_1\|_{L_{t, {\bf{x}}}^\infty}^2\int_0^t\|  b_1\|_{m-1}^2 ds+\|  b_1\|_{1,\infty}^2\int_0^t\|\partial_y v_1\|_{m-2}^2 ds\\
&+\| v_1\|_{L_{t, {\bf{x}}}^\infty}^2\int_0^t\| \partial_y b_1\|_{m-1}^2 ds+\| \partial_y b_1\|_{ 1,\infty}^2\int_0^t\| v_1\|_{m-2}^2 ds.
\end{align*}
As for the last  term on the right hand side of \eqref{4.22}, we have
\begin{align}\label{4.23}
\int_0^t\int_{\mathbb{R}_+^2}|\mathcal{Z}^\alpha \partial_y (b_1\partial_y v_2)|^2 d{\bf{x}}ds=\int_0^t\int_{\mathbb{R}_+^2}|\mathcal{Z}^\alpha (\partial_y b_1\partial_y  v_2)|^2 d{\bf{x}}ds+\int_0^t\int_{\mathbb{R}_+^2}|\mathcal{Z}^\alpha (b_1\partial_y^2 v_2)|^2 d{\bf{x}}ds.
\end{align}
Where the first term on the right hand side of \eqref{4.23} is estimated as follows due to \eqref{2.1}.
\begin{align*}
\int_0^t\int_{\mathbb{R}_+^2}|\mathcal{Z}^\alpha (\partial_y b_1\partial_y  v_2)|^2 d{\bf{x}}ds
\lesssim \|\partial_y b_1\|_{L_{t, {\bf{x}}}^\infty}^2\int_0^t\|\partial_y v_2\|_{m-2}^2 ds+\|\partial_y v_2\|_{L_{t, {\bf{x}}}^\infty}^2\int_0^t\|\partial_y b_1\|_{m-2}^2 ds.
\end{align*}
And for the second  term on the right hand side of \eqref{4.23}, by \eqref{4.3}, one has
\begin{equation}\label{4.24}
\begin{split}
\int_0^t\int_{\mathbb{R}_+^2}|\mathcal{Z}^\alpha (b_1\partial_y^2 v_2)|^2 d{\bf{x}}ds
\lesssim &\int_0^t\int_{\mathbb{R}_+^2}\left(|\mathcal{Z}^\alpha (b_1\partial_y  \partial_xv_1)|^2  +   |\mathcal{Z}^\alpha [b_1\partial_y   (p^{-1}\partial_t p)]|^2\right) d{\bf{x}}ds\\
&+ \int_0^t\int_{\mathbb{R}_+^2}|\mathcal{Z}^\alpha [b_1\partial_y ( p^{-1}{\bf{v}}\cdot\nabla p)]|^2 d{\bf{x}}ds.
\end{split}
\end{equation}
Then, the following inequality holds true for the first  term on the right hand side of  \eqref{4.24}.
\begin{align*}
\int_0^t\int_{\mathbb{R}_+^2}|\mathcal{Z}^\alpha (b_1\partial_y  \partial_xv_1)|^2 d{\bf{x}}ds
\lesssim &\|b_1\|_{L_{t,{\bf{x}}}^\infty}^2\int_0^t\|\partial_y v_1\|_{m-1}^2ds+\|\partial_y v_1\|_{1, \infty}^2\int_0^t\|b_1\|_{m-2}^2 ds.
\end{align*}
Similarly,  for the last two terms on the right hand side of  \eqref{4.24}, we have
\begin{align*}
&\int_0^t\int_{\mathbb{R}_+^2}|\mathcal{Z}^\alpha [b_1\partial_y  (p^{-1}\partial_t p)]|^2 d{\bf{x}}ds\\
\lesssim& \int_0^t\int_{\mathbb{R}_+^2}|\mathcal{Z}^\alpha (b_1\partial_y p^{-1}\partial_t p)|^2 d{\bf{x}}ds+\int_0^t\int_{\mathbb{R}_+^2}|\mathcal{Z}^\alpha (b_1 p^{-1}\partial_t\partial_y  p)|^2 d{\bf{x}}ds\\
\lesssim&\;\sum_{i, j=0}^1\|\partial_y^i(b_1, p, p^{-1})\|_{1, \infty}^4 \int_0^t\|\partial_y^j(b_1, p-1, p^{-1}-1)\|_{m-1}^2 ds,
\end{align*}
and
\begin{align*}
& \int_0^t\int_{\mathbb{R}_+^2}|\mathcal{Z}^\alpha [b_1\partial_y ( p^{-1}{\bf{v}}\cdot\nabla p)]|^2 d{\bf{x}}ds\\
\lesssim&\int_0^t\int_{\mathbb{R}_+^2}|\mathcal{Z}^\alpha (b_1\partial_y p^{-1}{\bf{v}}\cdot\nabla p)|^2 d{\bf{x}}ds+\int_0^t\int_{\mathbb{R}_+^2a}|\mathcal{Z}^\alpha (b_1 p^{-1}\partial_y{\bf{v}}\cdot\nabla p)|^2 d{\bf{x}}ds\\
&+\int_0^t\int_{\mathbb{R}_+^2}|\mathcal{Z}^\alpha (b_1 p^{-1}\phi^{-1}{\bf{v}}\cdot\phi\partial_y\nabla p)|^2 d{\bf{x}}ds\\
\lesssim&\sum_{i =0}^1\|\partial_y^i({\bf{v}}, b_1, p,    p^{-1})\|_{2, \infty}^6\cdot\int_0^t\left(\sum_{j =0}^1\|\partial_y^j({\bf{v}},  b_1, p-1, p^{-1}-1)\|_{m-1}^2+\|p-1\|_{m}^2\right) ds.
\end{align*}
After combining all of estimates in this subsection, we obtain
\begin{align*}
\int_0^t\|  \partial_y^2 v_1\|_{m-2}^2 ds
\lesssim&\left(1+\sum_{i=0}^1\|\partial_y^i({\bf{v}}, b_1, p, p^{-1} )\|_{2, \infty}^2\right)^3\\ &\cdot\int_0^t\left(\sum_{j=0}^1\|\partial_y^j({\bf{v}}, b_1, p-1, p^{-1}-1 )\|_{m-1}^2+\|p-1\|_{m}^2\right) ds.
\end{align*}

\subsubsection{The Conormal Estimate of $  \partial_y^2 v_2$} Next, we will derive the conormal estimates of $\partial_y^2 v_2$. Similar as \eqref{4.4}, it follows that
\begin{align}\label{4.25}
\begin{split}
\int_0^t\int_{\mathbb{R}_+^2}|\mathcal{Z}^\alpha\partial_y^{2}  v_2|^2 d{\bf{x}}ds
\lesssim& \int_0^t\int_{\mathbb{R}_+^2}|\mathcal{Z}^\alpha\partial_y \partial_x v_1|^2d{\bf{x}}ds+ \int_0^t\int_{\mathbb{R}_+^2}\left|\mathcal{Z}^\alpha\partial_y ( p^{-1}\partial_t p)\right|^2d{\bf{x}}ds\\
&+ \int_0^t\int_{\mathbb{R}_+^2}\left| \mathcal{Z}^\alpha\partial_y (p^{-1}{\bf{v}}\cdot \nabla p)\right|^2 d{\bf{x}}ds.
\end{split}
\end{align}
The first and second terms on the right hand side of \eqref{4.25} are estimated as follows.
\begin{align*}
\int_0^t\int_{\mathbb{R}_+^2}|\mathcal{Z}^\alpha\partial_y \partial_x v_1|^2 d{\bf{x}}ds\lesssim\int_0^t \|\partial_y  v_1\|_{m-1}^2 ds,
\end{align*}
and
\begin{align*}
&\int_0^t\int_{\mathbb{R}_+^2}|\mathcal{Z}^\alpha\partial_y (  p^{-1}\partial_t p)|^2 d{\bf{x}}ds\\
=&\int_0^t\int_{\mathbb{R}_+^2}|\mathcal{Z}^\alpha(\partial_y p^{-1}\partial_t p)|^2 d{\bf{x}}ds+\int_0^t\int_{\mathbb{R}_+^2}|\mathcal{Z}^\alpha( p^{-1}\partial_t\partial_y p)|^2 d{\bf{x}}ds\\
\lesssim &\;\|\partial_y p^{-1}\|_{L_{t, {\bf{x}}}^\infty}^2\int_0^t\|p-1\|_{m-1}^2ds+\|p\|_{1, \infty}^2\int_0^t\|\partial_y p^{-1}\|_{m-2}^2 ds\\
&+\| p^{-1}\|_{L_{t, {\bf{x}}}^\infty}^2\int_0^t\|\partial_yp\|_{m-1}^2ds+\|\partial_yp\|_{1, \infty}^2\int_0^t\|  p^{-1}-1\|_{m-2}^2 ds.
\end{align*}
The third term on the right hand side of \eqref{4.25} can be split into two parts.
\begin{equation}\label{4.26}
\begin{split}
&\int_0^t\int_{\mathbb{R}_+^2}| \mathcal{Z}^\alpha\partial_y (p^{-1}{\bf{v}}\cdot \nabla p)|^2d{\bf{x}}ds\\
=&\int_0^t\int_{\mathbb{R}_+^2}| \mathcal{Z}^\alpha\partial_y (p^{-1}v_1\partial_x p)|^2d{\bf{x}}ds
+\int_0^t\int_{\mathbb{R}_+^2}| \mathcal{Z}^\alpha\partial_y (p^{-1}v_2\partial_y p)|^2d{\bf{x}}ds.
\end{split}
\end{equation}
And each term on the right hand side of \eqref{4.26} can estimated accordingly.
\begin{align*}
&\int_0^t\int_{\mathbb{R}_+^2}| \mathcal{Z}^\alpha\partial_y ( p^{-1}v_1\partial_x p)|^2d{\bf{x}}ds\\
\lesssim&\int_0^t\int_{\mathbb{R}_+^2}| \mathcal{Z}^\alpha (\partial_y p^{-1}v_1\partial_x p)|^2d{\bf{x}}ds+\int_0^t\int_{\mathbb{R}_+^2}| \mathcal{Z}^\alpha ( p^{-1}\partial_yv_1\partial_x p)|^2d{\bf{x}}ds\\
&+\int_0^t\int_{\mathbb{R}_+^2}| \mathcal{Z}^\alpha ( p^{-1}v_1\partial_x \partial_y p)|^2d{\bf{x}}ds\\
\lesssim& \sum_{i, j=0}^1\|\partial_y^i(v_1, p, p^{-1} )\|_{1, \infty}^4  \int_0^t\|\partial_y^j(v_1, p-1, p^{-1}-1)\|_{m-1}^2 ds,
\end{align*}
and
\begin{align*}
&\int_0^t\int_{\mathbb{R}_+^2}| \mathcal{Z}^\alpha\partial_y (p^{-1}v_2\partial_y p)|^2d{\bf{x}}ds\\
\lesssim&\int_0^t\int_{\mathbb{R}_+^2}| \mathcal{Z}^\alpha (\partial_y p^{-1}v_2\partial_y p)|^2d{\bf{x}}ds+\int_0^t\int_{\mathbb{R}_+^2}| \mathcal{Z}^\alpha ( p^{-1}\partial_yv_2\partial_y p)|^2d{\bf{x}}ds\\
&+\int_0^t\int_{\mathbb{R}_+^2}| \mathcal{Z}^\alpha ( p^{-1}\phi^{-1}v_2 \phi\partial_y^2 p)|^2d{\bf{x}}ds\\
\lesssim& \sum_{i,j=0}^1\|\partial_y^i(v_2, p, p^{-1})\|_{1, \infty}^4 \int_0^t\|\partial_y^j(v_2, p-1, p^{-1}-1)\|_{m-1}^2 ds.
\end{align*}
Combining all of the estimates in this subsection produces the following inequality.
\begin{align*}
\int_0^t\| \partial_y^2  v_2\|_{m-2}^2 ds
\lesssim&  \left(1+\sum_{i=0}^1\|\partial_y^i({\bf{v}}, p, p^{-1})\|_{1,\infty}^2\right)^2\cdot\sum_{j=0}^1\int_0^t\|\partial_y^j({\bf{v}}, p-1, p^{-1}-1)\|_{m-1}^2 ds.
\end{align*}

\subsubsection{Conormal Estimate of $ \partial_y^2 p$} In this subsection, we derive the conormal estimate for the second order normal derivatives of pressure. Similar to \eqref{4.6}, we have
\begin{align}\label{4.27}
&\int_0^t\int_{\mathbb{R}_+^2}\left(|\mathcal{Z}^\alpha \partial_y^2 p |^2  +\varepsilon ^2(2\mu+\lambda)^2  |\mathcal{Z}^\alpha \partial_y^3 v_2|^2\right) d{\bf{x}}ds\notag\\
&-2\varepsilon (2\mu+\lambda)\int_0^t\int_{\mathbb{R}_+^2}\mathcal{Z}^\alpha  \partial_y^3v_2\cdot \mathcal{Z}^\alpha \partial_y^2 p \ d{\bf{x}}ds\notag \\
\lesssim &\int_0^t\int_{\mathbb{R}_+^2} |\mathcal{Z}^\alpha\partial_y (\rho \partial_t v_2)|^2 d{\bf{x}}ds+\int_0^t\int_{\mathbb{R}_+^2} |\mathcal{Z}^\alpha\partial_y (b_1\partial_x b_2)|^2 d{\bf{x}}ds\\
& +\int_0^t\int_{\mathbb{R}_+^2} |\mathcal{Z}^\alpha\partial_y (\rho {\bf{v}}\cdot \nabla v_2)|^2 d{\bf{x}}ds+\int_0^t\int_{\mathbb{R}_+^2} |\mathcal{Z}^\alpha\partial_y (b_1\partial_y b_1)|^2 d{\bf{x}}ds\notag\\
&+\varepsilon^2\mu^2\int_0^t\int_{\mathbb{R}_+^2} |\mathcal{Z}^\alpha\partial_y \partial_x^2 v_2|^2 d{\bf{x}}ds +\varepsilon^2 (\mu+\lambda)^2\int_0^t\int_{\mathbb{R}_+^2}|\mathcal{Z}^\alpha\partial_y^2 \partial_x v_1|^2 d{\bf{x}}ds.\notag
\end{align}
By using \eqref{3.4}, the second term on the left hand side of \eqref{4.27} is decomposed into the three terms.
\begin{equation}\label{4.28}
\begin{split}
&2\varepsilon (2\mu+\lambda)\int_0^t\int_{\mathbb{R}_+^2}\mathcal{Z}^\alpha \partial_y^3v_2\cdot\mathcal{Z}^\alpha \partial_y^2p\ d{\bf{x}}ds\\
=&-2\gamma^{-1}\varepsilon (2\mu+\lambda)\int_0^t\int_{\mathbb{R}_+^2}\mathcal{Z}^\alpha \partial_y^2( p^{-1}\partial_t p)\cdot\mathcal{Z}^\alpha \partial_y^2 p\ d{\bf{x}}ds\\
& -2\gamma^{-1}\varepsilon (2\mu+\lambda)\int_0^t\int_{\mathbb{R}_+^2}\mathcal{Z}^\alpha \partial_y^2( p^{-1}{\bf{v}}\cdot\nabla p)\cdot\mathcal{Z}^\alpha \partial_y^2 p\ d{\bf{x}}ds\\
&-2\varepsilon (2\mu+\lambda)\int_0^t\int_{\mathbb{R}_+^2}\mathcal{Z}^\alpha \partial_y^2\partial_xv_1\cdot\mathcal{Z}^\alpha \partial_y^2p\ d{\bf{x}}ds.
\end{split}
\end{equation}
Let us handle term by term on the right hand side of \eqref{4.28}.  First,
\begin{equation}\label{4.29}
\begin{split}
&-2\gamma^{-1}\varepsilon (2\mu+\lambda)\int_0^t\int_{\mathbb{R}_+^2}\mathcal{Z}^\alpha \partial_y^{2} ( p^{-1}\partial_t p)\cdot\mathcal{Z}^\alpha \partial_y^2 p\ d{\bf{x}}ds\\
=&-4\gamma^{-1}\varepsilon (2\mu+\lambda)  \int_0^t\int_{\mathbb{R}_+^2}\mathcal{Z}^\alpha( \partial_yp^{-1} \partial_y\partial_t p) \cdot\mathcal{Z}^\alpha \partial_y^2 p\ d{\bf{x}}ds\\
&-2\gamma^{-1}\varepsilon (2\mu+\lambda)  \int_0^t\int_{\mathbb{R}_+^2}\mathcal{Z}^\alpha( \partial_y^2  p^{-1}  \partial_t p) \cdot\mathcal{Z}^\alpha \partial_y^2 p\ d{\bf{x}}ds\\
&-2\gamma^{-1}\varepsilon (2\mu+\lambda)\int_0^t\int_{\mathbb{R}_+^2} \mathcal{Z}^\alpha   (  p^{-1} \partial_t\partial_y^2 p ) \cdot\mathcal{Z}^\alpha \partial_y^2p\ d{\bf{x}}ds.
\end{split}
\end{equation}
For the first term on the right hand side of \eqref{4.29}, by  \eqref{2.1}, we have
\begin{align*}
&-4\gamma^{-1}\varepsilon (2\mu+\lambda) \int_0^t\int_{\mathbb{R}_+^2}\mathcal{Z}^\alpha( \partial_y p^{-1}  \partial_y \partial_t p) \cdot\mathcal{Z}^\alpha \partial_y^2 p\ d{\bf{x}}ds\\
\lesssim &\;\varepsilon (2\mu+\lambda) \|\partial_y  p^{-1}\|_{L_{t,{\bf{x}}}^\infty}\left(\int_0^t\|\partial_y  p\|^2_{m-1} ds\right)^\frac12\left(\int_0^t\|\partial_y^2  p\|^2_{m-2} ds\right)^\frac12\\
&+\varepsilon (2\mu+\lambda) \|\partial_y p\|_{1, \infty}\left(\int_0^t\|\partial_y  p^{-1}\|^2_{m-2}ds\right)^\frac12\left(\int_0^t\|\partial_y^2 p\|^2_{m-2}ds\right)^\frac12.
\end{align*}
By the Sobolev embedding inequality, the second term on the right hand side of \eqref{4.29} is estimated as follows.
\begin{align*}
&-2\gamma^{-1}\varepsilon (2\mu+\lambda)  \int_0^t\int_{\mathbb{R}_+^2}\mathcal{Z}^\alpha( \partial_y^2 p^{-1}  \partial_t p) \cdot\mathcal{Z}^\alpha \partial_y^2 p\ d{\bf{x}}ds\\
\lesssim&\;\varepsilon (2\mu+\lambda)\left[\sum_{\beta+\gamma=\alpha\atop |\beta|<|\gamma| }\sup_{0\leq s\leq t}\|\mathcal{Z}^\beta\partial_y^2 p^{-1}(s)\|_{L_{x}^{\infty}(L_y^2)}\left(\int_0^t\|\mathcal{Z}^\gamma \partial_t p\|_{L_{x}^{2}(L_y^\infty)}^2 ds\right)^\frac12\right.\\
&\left.+ \sum_{\beta+\gamma=\alpha\atop |\beta|\geq|\gamma| }\|\mathcal{Z}^\gamma \partial_t p\|_{L_{t, {\bf{x}}}^\infty}\left(\int_0^t\int_{\mathbb{R}_+^2}|\mathcal{Z}^\beta \partial_y^{2}  p^{-1} |^2 d{\bf{x}}ds\right)^\frac12\right]\left(\int_0^t\|\partial_y^2p \|_{m-2}^2 ds\right)^\frac12\\
\lesssim &\; \varepsilon (2\mu+\lambda)\left[\|\partial_y^2 p^{-1}(0)\|_{[m/2]}+\left(\int_0^t\|\partial_y^2 p^{-1} \|_{[m/2]+1}^2 ds\right)^\frac12\right]\\
&\cdot\left(\int_0^t\|(p-1, \partial_y p)\|_{m-1}^2 ds\right)^\frac12 \left(\int_0^t\|\partial_y^2 p\|_{m-2}^2 ds\right)^\frac12\\
&+\varepsilon (2\mu+\lambda)\|p\|_{[m/2], \infty}\left(\int_0^t\|\partial_y^2 p^{-1}\|_{m-2}^2 ds\right)^\frac12 \left(\int_0^t\|\partial_y^2 p\|_{m-2}^2 ds\right)^\frac12.
\end{align*}
As for the last term on the right hand side of \eqref{4.29}, one has
\begin{equation}\label{4.30}
\begin{split}
&-2\gamma^{-1}\varepsilon (2\mu+\lambda)\int_0^t\int_{\mathbb{R}_+^2} \mathcal{Z}^\alpha(   p^{-1}\partial_t\partial_y^2 p ) \cdot\mathcal{Z}^\alpha \partial_y^2 p\ d{\bf{x}}ds\\
=&-2\gamma^{-1}\varepsilon (2\mu+\lambda)\sum_{\beta+\gamma=\alpha\atop |\beta|\geq1}\int_0^t\int_{\mathbb{R}_+^2} \mathcal{Z}^\beta   p^{-1}  \mathcal{Z}^\gamma \partial_t\partial_y^2 p   \cdot\mathcal{Z}^\alpha \partial_y^2 p\ d{\bf{x}}ds\\
&-2\gamma^{-1}\varepsilon (2\mu+\lambda) \int_0^t\int_{\mathbb{R}_+^2}     p^{-1}\mathcal{Z}^\alpha \partial_t\partial_y^2 p   \cdot\mathcal{Z}^\alpha \partial_y^2 p\ d{\bf{x}}ds.
\end{split}
\end{equation}
Then, by the Sobolev embedding inequality, the first term on the right hand side of \eqref{4.30} satisfies the following inequality.
\begin{align*}
&-2\gamma^{-1}\varepsilon (2\mu+\lambda)\sum_{\beta+\gamma=\alpha\atop |\beta |\geq1}\int_0^t\int_{ \mathbb{R}_+^2} \mathcal{Z}^\beta p^{-1}\mathcal{Z}^\gamma \partial_t\partial_y^2 p   \cdot\mathcal{Z}^\alpha \partial_y^2 p\ d{\bf{x}}ds\\
\lesssim &\;\varepsilon (2\mu+\lambda)\sum_{\beta+\gamma=\alpha\atop |\beta|\geq|\gamma|} \sup_{0\leq s\leq t}\|\mathcal{Z}^\gamma\partial_t\partial_y^2 p(s)\|_{L_{x}^{\infty}(L_y^2)}\left(\int_0^t\|\mathcal{Z}^\beta p^{-1}\|_{L_{x}^{2}(L_y^\infty)}^2 ds\right)^\frac12 \left(\int_0^t\| \partial_y^2 p\|_{m-2}^2 ds\right)^\frac12\\
&+\varepsilon (2\mu+\lambda)\sum_{\beta+\gamma=\alpha\atop 1\leq |\beta|<|\gamma| <|\alpha|}\|\mathcal{Z}^\beta    p^{-1}\|_{L_{t, {\bf{x}}}^\infty}\left(\int_0^t\int_ {\mathbb{R}_+^2}|\mathcal{Z}^\gamma \partial_t\partial_y^2 p   |^2d{\bf{x}}ds\right)^\frac12\left(\int_0^t\|\partial_y^2 p\|_{m-2}^2 ds\right)^\frac12\\
\lesssim&\; \varepsilon (2\mu+\lambda)\left[\|\partial_y^2 p(0)\|_{[m/2]+1}+\left(\int_0^t\|\partial_y^2 p \|_{[m/2]+2}^2 ds\right)^\frac12\right] \left(\int_0^t\| ( p^{-1}-1,\partial_yp^{-1})\|_{m-2}^2 ds\right)^\frac12\\
&\cdot\left(\int_0^t\|\partial_y^2 p\|_{m-2}^2 ds\right)^\frac12 +\varepsilon (2\mu+\lambda)\|p^{-1}\|_{[m/2]-1,\infty}\int_0^t\|\partial_y^2 p\|_{m-2}^2 ds.
\end{align*}
And the second term on the right hand side of \eqref{4.30} is dealt with as follows.
\begin{align*}
&-2\gamma^{-1}\varepsilon (2\mu+\lambda) \int_0^t\int_{\mathbb{R}_+^2}    p^{-1}\mathcal{Z}^\alpha \partial_t\partial_y^2 p   \cdot\mathcal{Z}^\alpha \partial_y^2 p\ d{\bf{x}}ds\\
= &-\gamma^{-1}\varepsilon (2\mu+\lambda)\frac{d}{dt}\int_0^t\int_{\mathbb{R}_+^2}  p^{-1}|\mathcal{Z}^\alpha  \partial_y^2 p |^2  \ d{\bf{x}}ds +\gamma^{-1}\varepsilon (2\mu+\lambda)\int_0^t\int_{\mathbb{R}_+^2} \partial_t p^{-1}|\mathcal{Z}^\alpha  \partial_y^2 p |^2  \ d{\bf{x}}ds\\
\lesssim &-\gamma^{-1}\varepsilon (2\mu+\lambda) \int_{\mathbb{R}_+^2}     p^{-1}(t)|\mathcal{Z}^\alpha  \partial_y^2 p (t)|^2  \ d{\bf{x}} +\gamma^{-1}\varepsilon (2\mu+\lambda) \int_{\mathbb{R}_+^2}     p_0^{-1}|\mathcal{Z}^\alpha  \partial_y^2 p_0 |^2  \ d{\bf{x}} \\
&+\varepsilon (2\mu+\lambda)\|p^{-1}\|_{1, \infty}\int_0^t\|\partial_y^2 p\|_{m-2}^2 ds.
\end{align*}
After a direct calculation, the second  term on the right hand side of \eqref{4.28} is rewritten as following form.
\begin{equation}\label{4.31}
\begin{split}
& -2\gamma^{-1}\varepsilon (2\mu+\lambda)\int_0^t\int_{\mathbb{R}_+^2}\mathcal{Z}^\alpha \partial_y^{2}  (p^{-1}{\bf{v}}\cdot\nabla p)\cdot\mathcal{Z}^\alpha \partial_y^2 p\ d{\bf{x}}ds\\
=& -2\gamma^{-1} \varepsilon (2\mu+\lambda) \int_0^t\int_{\mathbb{R}_+^2}\mathcal{Z}^\alpha ( \partial_y^{2}p^{-1}{\bf{v}}\cdot\nabla  p+p^{-1}\partial_y^{2}{\bf{v}}\cdot\nabla  p\\
&+2\partial_y p^{-1}v_2\partial_y^2 p+2p^{-1}\partial_yv_2\partial_y^2 p)\cdot\mathcal{Z}^\alpha \partial_y^2 p\ d{\bf{x}}ds\\
&-4\gamma^{-1}\varepsilon (2\mu+\lambda) \int_0^t\int_{\mathbb{R}_+^2}\mathcal{Z}^\alpha  (\partial_y p^{-1}\partial_y{\bf{v}} \cdot\nabla  p +\partial_y p^{-1}v_1\partial_x\partial_yp\\
&+ p^{-1}\partial_yv_1\partial_x\partial_yp )\cdot\mathcal{Z}^\alpha \partial_y^2 p\ d{\bf{x}}ds\\
& -2\gamma^{-1}\varepsilon (2\mu+\lambda) \int_0^t\int_{\mathbb{R}_+^2}\mathcal{Z}^\alpha (p^{-1}{\bf{v}}\cdot\nabla \partial_y^{2}p  )\cdot\mathcal{Z}^\alpha \partial_y^2 p\ d{\bf{x}}ds\\
=:&I_1+I_2+I_3.
\end{split}
\end{equation}
The first term on the right hand side of \eqref{4.31} can be estimated as follows due to the Sobolev embedding inequality.
\begin{align*}
I_1\;\lesssim&\;\varepsilon (2\mu+\lambda)\sum_{\beta+\gamma+\iota=\alpha\atop |\beta|\geq|\gamma|, |\iota|}\|\mathcal{Z}^\gamma({\bf{v}}, p^{-1}, \partial_y v_2, \nabla p, \partial_y p^{-1})\|_{L_{t, {\bf{x}}}^\infty}\|\mathcal{Z}^\iota({\bf{v}}, p^{-1}, \partial_y v_2\nabla p, \partial_y p^{-1})\|_{L_{t, {\bf{x}}}^\infty}\\
&\cdot\left(\int_0^t\int_{\mathbb{R}_+^2}|\mathcal{Z}^\beta\partial_y^2( {\bf{v}},   p, p^{-1} )|^2 d{\bf{x}}ds\right)^\frac12\left(\int_0^t\|\partial_y^2 p\|_{m-2}^2 ds\right)^\frac12\\
&+\varepsilon (2\mu+\lambda)
\sum_{\beta+\gamma+\iota=\alpha\atop |\gamma|\geq |\beta|, |\iota|}\sup_{0\leq s\leq t}\|\mathcal{Z}^\beta\partial_y^2(  {\bf{v}}, p, p^{-1} )(s)\|_{   L_{x}^{\infty}(L_y^2)}\|\mathcal{Z}^\iota({\bf{v}}, p^{-1}, \partial_y v_2, \nabla p, \partial_y p^{-1})\|_{L_{t, {\bf{x}}}^\infty}\\
&\cdot \left(\int_0^t\|\mathcal{Z}^\gamma({\bf{v}}, p^{-1}, \partial_y v_2, \nabla p, \partial_y p^{-1})\|_{L_{x}^{2}(L_y^\infty)}^2  ds\right)^\frac12\left(\int_0^t\|\partial_y^2 p\|_{m-2}^2 ds\right)^\frac12\\
&+\varepsilon (2\mu+\lambda)\sum_{\beta+\gamma+\iota=\alpha\atop |\iota|\geq |\beta|,|\gamma| }\sup_{0\leq s\leq t}\|\mathcal{Z}^\beta\partial_y^2(  {\bf{v}},  p, p^{-1} )(s)\|_{   L_{x}^{\infty}(L^2_y)}\|\mathcal{Z}^\gamma({\bf{v}}, p^{-1}, \partial_y v_2, \nabla p, \partial_y p^{-1})\|_{L_{t, {\bf{x}}}^\infty} \\
&\cdot\left(\int_0^t\|\mathcal{Z}^\iota({\bf{v}}, p^{-1}, \partial_y v_2, \nabla p, \partial_y p^{-1})\|_{L_{x}^{2}(L_y^\infty)}^2  ds\right)^\frac12\left(\int_0^t\|\partial_y^2 p\|_{m-2}^2 ds\right)^\frac12\\
\lesssim&\;\varepsilon (2\mu+\lambda)\sum_{i=0}^1\|\partial_y^i({\bf{v}},  p, p^{-1} )\|_{[m/2], \infty}^2\left(\int_0^t\|\partial_y^2(  {\bf{v}},  p, p^{-1})\|_{m-2}^2 ds\right)^\frac12 \left(\int_0^t\|\partial_y^2 p\|_{m-2}^2 ds\right)^\frac12\notag\\ &+\varepsilon (2\mu+\lambda)\left[\|\partial_y^2( {\bf{v}}, p, p^{-1})(0)\|_{[m/2]} +\left(\int_0^t\|\partial_y^2( {\bf{v}}, p, p^{-1})\|_{[m/2]+1}^2 ds\right)^\frac12\right]\\
&\cdot\sum_{i=0}^1\|\partial_y^i({\bf{v}},  p, p^{-1} )\|_{[m/2], \infty} \sum_{j=0}^2\left(\int_0^t\|\partial_y^j({\bf{v}}, p-1, p^{-1}-1)\|_{m-j}^2     ds\right)^\frac12 \left(\int_0^t\|\partial_y^2 p\|_{m-2}^2 ds\right)^\frac12.
\end{align*}
And by \eqref{2.1}, one has
\begin{align*}
I_2\lesssim\varepsilon (2\mu+\lambda)  \sum_{i,j=0}^1\|\partial_y^i({\bf{v}}, p, p^{-1})\|_{1, \infty}^2\left(\int_0^t\|\partial_y^j({\bf{v}},  p-1, p^{-1}-1)\|_{m-1}^2 ds\right)^\frac12\left(\int_0^t\|\partial_y^2 p\|_{m-2}^2 ds\right)^\frac12.
\end{align*}
As for the last term on the right hand side of \eqref{4.31},
\begin{equation}\label{4.32}
\begin{split}
I_3=& -2\gamma^{-1}\varepsilon (2\mu+\lambda)\sum_{\beta+\gamma+\iota=\alpha\atop  |\iota|\neq|\alpha|} \int_0^t\int_{\mathbb{R}_+^2}\mathcal{Z}^\beta p^{-1}\mathcal{Z}^\gamma {\bf{v}}\cdot\mathcal{Z}^\iota\nabla \partial_y^{2}p  \cdot\mathcal{Z}^\alpha \partial_y^2 p\ d{\bf{x}}ds\\
& -2\gamma^{-1}\varepsilon (2\mu+\lambda) \int_0^t\int_{\mathbb{R}_+^2}  p^{-1}{\bf{v}}\cdot\mathcal{Z}^\alpha\nabla \partial_y^{2}p\cdot\mathcal{Z}^\alpha \partial_y^2 p\ d{\bf{x}}ds.
\end{split}
\end{equation}
Then, the first part on the right hand of \eqref{4.32} satisfies
\begin{align*}
& -2\gamma^{-1}\varepsilon (2\mu+\lambda)\sum_{\beta+\gamma+\iota=\alpha\atop  |\iota|\neq|\alpha|} \int_0^t\int_{\mathbb{R}_+^2}\mathcal{Z}^\beta p^{-1}\mathcal{Z}^\gamma {\bf{v}}\cdot\mathcal{Z}^\iota\nabla \partial_y^{2}p  \cdot\mathcal{Z}^\alpha \partial_y^2 p\ d{\bf{x}}ds\\
\lesssim&\;\varepsilon (2\mu+\lambda)\left[\sum_{\beta+\gamma+\iota=\alpha\atop   |\beta|\geq|\gamma|, |\iota|}\sup_{0\leq s\leq t}\|\phi\mathcal{Z}^\iota \nabla\partial_y^2 p(s)\|_{L_{x}^{\infty}(L_y^2)}\|\phi^{-1}\mathcal{Z}^\gamma {\bf{v}}\|_{L_{t, {\bf{x}}}^\infty} \left(\int_0^t\|\mathcal{Z}^\beta p^{-1}\|_{L_{x}^{2}(L_y^\infty)}^2 ds\right)^\frac12\right.\\
&\left.+\sum_{\beta+\gamma+\iota=\alpha\atop  |\gamma||\geq|\beta, |\iota|}\sup_{0\leq s\leq t}\|\phi\mathcal{Z}^\iota \nabla\partial_y^2 p(s)\|_{L_{x}^{\infty}(L_y^2)}\|\mathcal{Z}^\beta p^{-1}\|_{L_{t, {\bf{x}}}^\infty} \left(\int_0^t\|\phi^{-1}\mathcal{Z}^\gamma {\bf{v}}\|_{L_{x}^{2}(L_y^\infty)}^2 ds\right)^\frac12\right.\\
&\left.+\sum_{\beta+\gamma+\iota=\alpha\atop |\beta,|\gamma|\leq |\iota|<|\alpha| }\|\mathcal{Z}^\beta p^{-1}\|_{L_{t, {\bf{x}}}^\infty}\|\phi^{-1}\mathcal{Z}^\gamma {\bf{v}}\|_{L_{t, {\bf{x}}}^\infty}\left(\int_0^t\int_{\mathbb{R}_+^2}|\phi\mathcal{Z}^\iota\nabla\partial_y^2 p|^2 d{\bf{x}}ds\right)^\frac12\right]\cdot\left(\int_0^t\|\partial_y^2 p\|_{m-2}^2 ds\right)^\frac12\\
\lesssim&\; \varepsilon (2\mu+\lambda)\|(p^{-1}, \partial_y {\bf{v}})\|_{[m/2]-1, \infty}\sum_{i=1}^2\left[\|\partial_y^i p(0)\|_{[m/2]+2}+\left(\int_0^t\|\partial_y^i p\|_{[m/2]+3}^2 ds\right)^\frac12\right]\\
&\cdot\sum_{j=0}^2\left(\int_0^t\|\partial_y^j({\bf{v}}, p^{-1})\|_{m-2}^2 ds\right)^\frac12\left(\int_0^t\|\partial_y^2 p\|_{m-2}^2 ds\right)^\frac12\\
&+\varepsilon (2\mu+\lambda)\|p^{-1}\|_{[m/2]-1, \infty}\|\partial_y {\bf{v}}\|_{[m/2]-1, \infty} \left[\int_0^t\left(\|\partial_y p\|_{m-1}^2+ \|\partial_y^2 p\|_{m-2}^2\right)  ds\right]^\frac12\\
&\cdot\left(\int_0^t\|\partial_y^2 p\|_{m-2}^2 ds\right)^\frac12.
\end{align*}
And it is convenient to rewrite the last term on the right hand side of \eqref{4.32} as two parts.
\begin{equation}\label{4.33}
\begin{split}
& -2\gamma^{-1}\varepsilon (2\mu+\lambda) \int_0^t\int_{\mathbb{R}_+^2}  p^{-1}{\bf{v}}\cdot\mathcal{Z}^\alpha\nabla \partial_y^{2}p\cdot\mathcal{Z}^\alpha \partial_y^2 p\ d{\bf{x}}ds\\
= & -2\gamma^{-1}\varepsilon (2\mu+\lambda) \int_0^t\int_{\mathbb{R}_+^2}  p^{-1}{\bf{v}}\cdot[\mathcal{Z}^\alpha, \nabla] \partial_y^{2}p\cdot\mathcal{Z}^\alpha \partial_y^2 p\ d{\bf{x}}ds\\
& -2\gamma^{-1}\varepsilon (2\mu+\lambda) \int_0^t\int_{\mathbb{R}_+^2}  p^{-1}{\bf{v}}\cdot\nabla\mathcal{Z}^\alpha \partial_y^{2}p\cdot\mathcal{Z}^\alpha \partial_y^2 p\ d{\bf{x}}ds.
\end{split}
\end{equation}
Then, for the first term on the right hand side of \eqref{4.33},  by \eqref{2.2}, one obtains
\begin{align*}
& -2\gamma^{-1}\varepsilon (2\mu+\lambda) \int_0^t\int_{\mathbb{R}_+^2}  p^{-1}{\bf{v}}\cdot[\mathcal{Z}^\alpha, \nabla] \partial_y^{2}p\cdot\mathcal{Z}^\alpha \partial_y^2 p\ d{\bf{x}}ds\\
=& -2\gamma^{-1}\varepsilon (2\mu+\lambda)\sum_{k=0}^{m-3} \int_0^t\int_{\mathbb{R}_+^2} \phi^{k, m-2} (y)p^{-1}{\bf{v}}\cdot\partial_y\mathcal{Z}_2^{k} \partial_y^{2}p\cdot\mathcal{Z}^\alpha \partial_y^2 p\ d{\bf{x}}ds\\
\lesssim&\;\varepsilon (2\mu+\lambda)\sum_{k=0}^{m-3}
\|p^{-1}\|_{L_{t,{\bf{x}}}^\infty}\|\phi^{-1}{\bf{v}}\|_{L_{t,{\bf{x}}}^\infty}\left
(\int_0^t\int_{\mathbb{R}_+^2}|\phi\partial_y\mathcal{Z}_2^{k } \partial_y^{2}p|^2 d{\bf{x}}ds\right)^\frac12\left(\int_0^t\|\partial_y^2 p\|_{m-2}^2 ds\right)^\frac12\\
\lesssim&\;\varepsilon (2\mu+\lambda)\|p^{-1} \|_{L_{t,{\bf{x}}}^\infty}\|\partial_y {\bf{v}} \|_{L_{t,{\bf{x}}}^\infty}\int_0^t\|\partial_y^2 p\|_{m-2}^2 ds.
\end{align*}
By integration by parts, the second term on the right hand side of \eqref{4.33} is handled in the following way.
\begin{align*}
& -2\gamma^{-1}\varepsilon (2\mu+\lambda) \int_0^t\int_{\mathbb{R}_+^2}  p^{-1}{\bf{v}}\cdot\nabla\mathcal{Z}^\alpha \partial_y^{2}p\cdot\mathcal{Z}^\alpha \partial_y^2 p\ d{\bf{x}}ds\\
\lesssim &\;\varepsilon (2\mu+\lambda) \| ({\bf{v}}, p^{-1})\|_{L_{t, {\bf{x}}}^\infty} \| \nabla({\bf{v}}, p^{-1})\|_{L_{t, {\bf{x}}}^\infty} \int_0^t\|\partial_y^2 p\|_{m-2}^2 ds.
\end{align*}
Applying the Young's inequality to the last term on the right hand side of \eqref{4.28} yields that
\begin{align*}
&-2\varepsilon (2\mu+\lambda)\int_0^t\int_{\mathbb{R}_+^2}\mathcal{Z}^\alpha \partial_y^2\partial_xv_1\cdot\mathcal{Z}^\alpha \partial_y^2p\ d{\bf{x}}ds\\
\leq & 2\varepsilon ^2(2\mu+\lambda)^2\int_0^t\|\partial_y^2 v_1\|_{m-1}^2 ds+\frac12\int_0^t\|\partial_y^2 p\|_{m-2}^2ds.
\end{align*}

Next, it is left to estiamte the terms on the right hand side of \eqref{4.27}. For the first term on the right hand side of \eqref{4.27}, by  \eqref{2.1},  one derives
\begin{align*}
\int_0^t\int_{\mathbb{R}_+^2} |\mathcal{Z}^\alpha\partial_y(\rho \partial_t v_2)|^2 d{\bf{x}}ds
\lesssim&\int_0^t\int_{\mathbb{R}_+^2} |\mathcal{Z}^\alpha(\partial_y\rho \partial_t v_2)|^2 d{\bf{x}}ds+	\int_0^t\int_{\mathbb{R}_+^2} |\mathcal{Z}^\alpha(\rho \partial_t \partial_yv_2)|^2 d{\bf{x}}ds\\
\lesssim&\sum_{i, j=0}^1\|\partial_y^i(\rho, v_2)\|_{1, \infty}^2\int_0^t\|\partial_y^j(\rho-1, v_2)\|_{m-1}^2 ds.
\end{align*}
Similarly,
\begin{align*}
\int_0^t\int_{\mathbb{R}_+^2} |\mathcal{Z}^\alpha\partial_y (b_1\partial_x b_2)|^2 d{\bf{x}}ds
\lesssim\sum_{i, j=0}^1\|\partial_y^i(b_1, b_2)\|_{1, \infty}^2\int_0^t\|\partial_y^j(b_1, \tilde{b}_2)\|_{m-1}^2 ds.
\end{align*}
To estimate the third term on the right hand side of \eqref{4.27}, we decompose it into two terms.
\begin{equation}\label{4.34}
\begin{split}
&\int_0^t\int_{\mathbb{R}_+^2} |\mathcal{Z}^\alpha\partial_y (\rho {\bf{v}}\cdot \nabla v_2)|^2 d{\bf{x}}ds\\
= &\int_0^t\int_{\mathbb{R}_+^2} |\mathcal{Z}^\alpha[\partial_y(\rho {\bf{v}})\cdot \nabla v_2]|^2 d{\bf{x}}ds
+\int_0^t\int_{\mathbb{R}_+^2} |\mathcal{Z}^\alpha(\rho {\bf{v}}\cdot \nabla \partial_yv_2)|^2 d{\bf{x}}ds.
\end{split}
\end{equation}
For the first term on the right hand side of \eqref{4.34}, by \eqref{2.1}, one has
\begin{align*}
&\int_0^t\int_{\mathbb{R}_+^2} |\mathcal{Z}^\alpha[\partial_y(\rho {\bf{v}})\cdot  \nabla v_2]|^2 d{\bf{x}}ds\\
\lesssim&\int_0^t\int_{\mathbb{R}_+^2} |\mathcal{Z}^\alpha(\partial_y\rho {\bf{v}}\cdot  \nabla v_2)|^2 d{\bf{x}}ds+\int_0^t\int_{\mathbb{R}_+^2} |\mathcal{Z}^\alpha(\rho \partial_y{\bf{v}}\cdot  \nabla v_2)|^2 d{\bf{x}}ds\\
\lesssim & \sum_{i, j=0}^1 \|\partial_y^i (\rho, {\bf{v}})  \|_{1,\infty}^4\int_0^t\|\partial_y^j (\rho, {\bf{v}}) \|_{m-1}^2 ds.
\end{align*}
By a similar argument, the second term on the right hand side of \eqref{4.34} is estimated.
\begin{align*}
\int_0^t\int_{\mathbb{R}_+^2} |\mathcal{Z}^\alpha(\rho {\bf{v}}\cdot \nabla \partial_yv_2)|^2 d{\bf{x}}ds
\lesssim &\; \|(\rho, v_1, \partial_x\partial_y v_2)\|_{L_{t, {\bf{x}}}^\infty}^4\int_0^t\|(\rho, v_1, \partial_x\partial_y v_2)\|_{m-2}^2 ds\\
& +\|(\rho, \phi^{-1}v_2, \phi\partial_y^2 v_2)\|_{L_{t, {\bf{x}}}^\infty}^4\int_0^t\|(\rho, \phi^{-1}v_2, \phi\partial_y^2 v_2)\|_{m-2}^2 ds\\
\lesssim &\; \|(\rho, v_1,  \partial_y v_2)\|_{1,\infty}^4\int_0^t\|(\rho, v_1,  \partial_y  v_2)\|_{m-1}^2 ds.
\end{align*}

As for the fourth term on the right hand side of \eqref{4.27},  one has
\begin{equation}\label{4.35}
\begin{split}
\int_0^t\int_{\mathbb{R}_+^2} |\mathcal{Z}^\alpha\partial_y (b_1\partial_y b_1)|^ 2 d{\bf{x}}ds
=\int_0^t\int_{\mathbb{R}_+^2} |\mathcal{Z}^\alpha(\partial_yb_1\partial_y  b_1)|^ 2 d{\bf{x}}ds+\int_0^t\int_{\mathbb{R}_+^2} |\mathcal{Z}^\alpha (b_1\partial_y^2 b_1)|^ 2 d{\bf{x}}ds.
\end{split}
\end{equation}
Where the first term on the right hand of \eqref{4.35} satisfies
\begin{align*}
\int_0^t\int_{\mathbb{R}_+^2} |\mathcal{Z}^\alpha(\partial_yb_1\partial_y  b_1)|^ 2 d{\bf{x}}ds\lesssim \|\partial_y b_1\|_{L_{t, {\bf{x}}}^\infty}^2\int_0^t\|\partial_y b_1\|_{m-2}^2 ds,
\end{align*}
and the second term on the right hand side of \eqref{4.35} is estimated as follows by the Sobolev embedding inequality.
\begin{align*}
\int_0^t\int_{\mathbb{R}_+^2} |\mathcal{Z}^\alpha (b_1\partial_y^2 b_1)|^ 2 d{\bf{x}}ds
\lesssim&\sum_{\beta+\gamma=\alpha\atop|\beta|\geq|\gamma|}\sup_{0\leq s\leq t}\|\mathcal{Z}^\gamma\partial_y^2 b_1(s)\|_{L_{x}^{\infty}(L_y^2)}^2\int_0^t\|\mathcal{Z}^\beta b_1\|_{L_{x}^{2}(L_y^\infty)}^2 ds\\
&+\sum_{\beta+\gamma=\alpha\atop|\beta|<|\gamma|}\|\mathcal{Z}^\beta b_1\|_{L_{t,x}^\infty}^2\int_0^t\int_{\mathbb{R}_+^2}|\mathcal{Z}^\gamma\partial_y^2 b_1|^2 d{\bf{x}}ds\\
\lesssim&\left(\|\partial_y^2 b_1(0)\|_{[m/2]}+\int_0^t\|\partial_y^2 b_1\|_{[m/2]+1}^2 ds\right) \int_0^t\|(b_1,\partial_y b_1)\|_{m-2}^2 ds\\
& +\|b_1\|_{[m/2]-1, \infty}^2\int_0^t\|\partial_y^2 b_1\|_{m-2}^2 ds .
\end{align*}
Finally, it is direct to estimate the last two terms on the right hand side of \eqref{4.27}.
\begin{align*}
&\varepsilon^2\mu^2\int_0^t\int_{\mathbb{R}_+^2} |\mathcal{Z}^\alpha\partial_y \partial_x^2 v_2|^2 d{\bf{x}}ds +\varepsilon^2 \mu^2\int_0^t\int_{\mathbb{R}_+^2}|\mathcal{Z}^\alpha\partial_y^2 \partial_x v_1|^2 d{\bf{x}}ds\\
\lesssim&\;\varepsilon^2\mu^2\int_0^t\|\partial_y v_2\|_{m}^2 ds+\varepsilon^2\mu^2\int_0^t\|\partial_y^2 v_1\|_{m-1}^2 ds.
\end{align*}
Combining all of the estimates in this subsection, we obtain
\begin{align*}
&\int_0^t \left( \| \partial_y^2 p\|_{m-2}^2  +\varepsilon^2(2\mu+\lambda)^2
\| \partial_y^3  v_2\|_{m-2}^2 \right) ds +\varepsilon(2\mu+\lambda)\gamma^{-1}\sum_{|\alpha|\leq m-2}\int_{\mathbb{R}_+^2} p^{-1}(t)|\mathcal{Z}^\alpha \partial_y^2 p(t)|^2 d{\bf{x}}\\
\lesssim&\; \varepsilon(2\mu+\lambda)\gamma^{-1}\sum_{|\alpha|\leq m-2}\int_ {\mathbb{R}_+^2} p_0^{-1}|\mathcal{Z}^\alpha \partial_y^2 p_0|^2 d{\bf{x}}+\varepsilon^2\mu^2\int_0^t \|\partial_y v_2\|_{m}^2 ds \\
&+  \varepsilon^2(2\mu+\lambda)^2\int_0^t   \| \partial_y^2   v_1\|_{m-1}^2 ds    +\|b_1\|_{[m/2]-1, \infty}^2\int_0^t\|\partial_y^2 b_1\|_{m-2}^2 ds\\
&+\varepsilon(2\mu+\lambda)\left[\left(1+\sum_{i =0}^1\|\partial_y^i({\bf{v}}, p, p^{-1})\|_{[m/2], \infty}^2\right)^2
+  \sum_{i=1}^2\|\partial_y^i( {\bf{v}},  p, p^{-1})(0)\|_{[m/2]+2}^2\right.\\
&\left.+\sum_{i=1}^2\int_0^t\|\partial_y^i({\bf{v}},   p-1, p^{-1}-1)\|_{[m/2]+3}^2ds\right]\cdot\sum_{j=0}^2\int_0^t \|\partial_y^j({\bf{v}}, p-1, p^{-1}-1)\|_{m-j}^2   ds\\
&+\left(1+\sum_{i=0}^1\|\partial_y^i(\rho, {\bf{v}},  {\bf{B}})\|_{[m/2], \infty}^2\right)^2\cdot\sum_{j=0}^1\int_0^t\|\partial_y^j(\rho-1, {\bf{v}},  {\bf{B}}-\overset{\rightarrow}{e_y} )\|_{m-1}^2 ds\\
&+\left(\|\partial_y^2 b_1(0)\|_{[m/2]}^2+\int_0^t\|\partial_y^2 b_1\|_{[m/2]+1}^2 ds\right)\int_0^t\|(b_1,\partial_y b_1)\|_{m-2}^2 ds.
\end{align*}
It is noted that the smallness of $\|b_1\|_{[m/2]-1, \infty}$ will be required for later analysis.
\subsubsection{The Conormal Estimate for  $  \partial_y^2 b_1$}
This subsection is devoted to the conormal estimates of $\partial_y^2b_1$. Similarly to \eqref{4.14}, we have	
\begin{equation}\label{4.36}
\begin{split}
&	\int_0^t\int_{\mathbb{R}_+^2}\left( |\mathcal{Z}^\alpha\partial_y^2 b_1|^2  +\varepsilon^2\mu^2  |\mathcal{Z}^\alpha\partial_y^{3}  v_1|^2\right) d{\bf{x}}ds +2\varepsilon\mu\int_0^t\int_{\mathbb{R}_+^2} \mathcal{Z}^\alpha\partial_y^{3}  v_1 \cdot \mathcal{Z}^\alpha\partial_y^2 b_1 d{\bf{x}}ds\\
\lesssim&\int_0^t\int_{\mathbb{R}_+^2} |\mathcal{Z}^\alpha \partial_y(\rho\partial_t v_1)|^2 d{\bf{x}}ds+\int_0^t\int_{\mathbb{R}_+^2} |\mathcal{Z}^\alpha \partial_y(b_2\partial_x b_2)|^2 d{\bf{x}}ds\\
&+\int_0^t\int_{\mathbb{R}_+^2} |\mathcal{Z}^\alpha\partial_y (\rho {\bf{v}}\cdot\nabla v_1)|^2 d{\bf{x}}ds+\int_0^t\int_{\mathbb{R}_+^2} |\mathcal{Z}^\alpha \partial_x\partial_y p|^2 d{\bf{x}}ds \\
& +\int_0^t\int_{\mathbb{R}_+^2} |\mathcal{Z}^\alpha\partial_y (\tilde b_2\partial_y b_1)|^2 d{\bf{x}}ds+\varepsilon^2(2\mu+\lambda)^2\int_0^t\int_{\mathbb{R}_+^2}|\mathcal{Z}^\alpha \partial_x^2\partial_y v_1|^2d{\bf{x}}ds\\
&+\varepsilon^2(\mu+\lambda)^2\int_0^t\int_{\mathbb{R}_+^2} |\mathcal{Z}^\alpha \partial_x\partial_y^2v_2|^2 d{\bf{x}}ds.
\end{split}
\end{equation}
By using \eqref{4.1}, the second term on the left hand side of \eqref{4.36} can be rewritten as
\begin{equation}\label{4.37}
\begin{split}
&- 	2\varepsilon\mu\int_0^t\int_{\mathbb{R}_+^2} \mathcal{Z}^\alpha\partial_y^{3} v_1\cdot \mathcal{Z}^\alpha\partial_y^2 b_1d{\bf{x}}ds \\
=&-2\varepsilon\mu\int_0^t\int_{\mathbb{R}_+^2} \partial_t \mathcal{Z}^\alpha\partial_y^2	b_1\cdot \mathcal{Z}^\alpha\partial_y^2 b_1\ d{\bf{x}}ds+2\varepsilon\mu\int_0^t\int_{\mathbb{R}_+^2} \mathcal{Z}^\alpha\partial_y^2(\tilde{b}_2\partial_y v_1)\cdot \mathcal{Z}^\alpha\partial_y^2 b_1\ d{\bf{x}}ds\\
&-2\varepsilon\mu\int_0^t\int_{\mathbb{R}_+^2} \mathcal{Z}^\alpha\partial_y^2(v_2\partial_y b_1)\cdot \mathcal{Z}^\alpha\partial_y^2 b_1\ d{\bf{x}}ds-2\varepsilon\mu\int_0^t\int_{\mathbb{R}_+^2} \mathcal{Z}^\alpha\partial_y^2(b_1\partial_y v_2)\cdot \mathcal{Z}^\alpha\partial_y^2 b_1\ d{\bf{x}}ds\\ &-2\varepsilon\mu\int_0^t\int_{\mathbb{R}_+^2} \mathcal{Z}^\alpha\partial_y^2(v_1\partial_x b_1)\cdot \mathcal{Z}^\alpha\partial_y^2 b_1\ d{\bf{x}}ds.
\end{split}
\end{equation}
Below, we handle term by term on the right hand side of \eqref{4.37}. First,
\begin{align*}
-2\varepsilon\mu\int_0^t\int_{\mathbb{R}_+^2} \partial_t \mathcal{Z}^\alpha\partial_y^2	b_1\cdot\mathcal{Z}^\alpha\partial_y^2 b_1\ d{\bf{x}}ds
=-\varepsilon\mu\int_{\mathbb{R}_+^2} |\mathcal{Z}^\alpha\partial_y^2 b_1(t)|^2 d{\bf{x}}+\varepsilon\mu\int_{\mathbb{R}_+^2} |\mathcal{Z}^\alpha\partial_y^2 b_1(0)|^2 d{\bf{x}}.
\end{align*}
The second term on the right hand side of \eqref{4.37} includes three terms.
\begin{equation}\label{4.38}
\begin{split}
&2\varepsilon\mu\int_0^t\int_{\mathbb{R}_+^2}  \mathcal{Z}^\alpha\partial_y^2(\tilde{b}_2\partial_y v_1)\cdot\mathcal{Z}^\alpha\partial_y^2 b_1\ d{\bf{x}}ds\\
=&\;2\varepsilon\mu\int_0^t\int_{\mathbb{R}_+^2}  \mathcal{Z}^\alpha(\partial_y^2\tilde{b}_2\partial_y v_1)\cdot\mathcal{Z}^\alpha\partial_y^2 b_1\ d{\bf{x}}ds+4\varepsilon\mu\int_0^t\int_{\mathbb{R}_+^2}  \mathcal{Z}^\alpha(\partial_y\tilde{b}_2\partial_y^2 v_1)\cdot\mathcal{Z}^\alpha\partial_y^2 b_1\ d{\bf{x}}ds\\
&+2\varepsilon\mu\int_0^t\int_{\mathbb{R}_+^2}  \mathcal{Z}^\alpha(\tilde{b}_2\partial_y^3 v_1)\cdot\mathcal{Z}^\alpha\partial_y^2 b_1\ d{\bf{x}}ds.
\end{split}
\end{equation}
For the first term on the right hand side of \eqref{4.38}, since $\nabla\cdot {\bf{B}}=0$, by \eqref{2.1}, we have
\begin{align*}
&2\varepsilon\mu\int_0^t\int_{\mathbb{R}_+^2}  \mathcal{Z}^\alpha(\partial_y^2\tilde{b}_2\partial_y v_1)\cdot\mathcal{Z}^\alpha\partial_y^2 b_1\ d{\bf{x}}ds\\
=&-2\varepsilon\mu\int_0^t\int_{\mathbb{R}_+^2}  \mathcal{Z}^\alpha(\partial_y\partial_xb_1\partial_y v_1)\cdot\mathcal{Z}^\alpha\partial_y^2 b_1\ d{\bf{x}}ds\\
\lesssim &\; \varepsilon\mu\|\partial_y b_1\|_{1, \infty}\left(\int_0^t\|\partial_y v_1\|_{m-2}^2 ds\right)^\frac12\left(\int_0^t\|\partial_y^2 b_1\|_{m-2}^2 ds\right)^\frac12\\
&+\varepsilon\mu\|\partial_y v_1\|_{L_{t, {\bf{x}}}^\infty}\left(\int_0^t\|\partial_y b_1\|_{m-1}^2 ds\right)^\frac12\left(\int_0^t\|\partial_y^2 b_1\|_{m-2}^2 ds\right)^\frac12.
\end{align*}
Since $\nabla\cdot {\bf{B}}=0$, by the Sobolev embedding inequality, the second term on the right hand side of \eqref{4.38} is estimated as follows.
\begin{align*}
&4\varepsilon\mu\int_0^t\int_{\mathbb{R}_+^2}  \mathcal{Z}^\alpha(\partial_y\tilde{b}_2\partial_y^2 v_1)\cdot\mathcal{Z}^\alpha\partial_y^2 b_1\ d{\bf{x}}ds\\
\lesssim&\;\varepsilon\mu\sum_{\beta+\gamma=\alpha\atop|\beta|>|\gamma|}\sup_{0\le s\leq t}\|\mathcal{Z}^\gamma\partial_y^2 v_1(s)\|_{L_{x}^{\infty}(L_y^2)}\left(\int_0^t\|\mathcal{Z}^\beta \partial_y\tilde b_2\|_{L_{x}^{2}(L_y^\infty)}^2 ds\right)^\frac12\left(\int_0^t\|\partial_y^2 b_1\|_{m-2}^2 ds\right)^\frac12\\
&+\varepsilon\mu\sum_{\beta+\gamma=\alpha\atop|\beta|\leq|\gamma|}\|\mathcal{Z}^\beta \partial_y\tilde b_2\|_{L_{t, {\bf{x}}}^\infty}\left(\int_0^t\int_{\mathbb{R}_+^2}|\mathcal{Z}^\gamma\partial_y^2 v_1|^2 d{\bf{x}}ds\right)^\frac12\left(\int_0^t\|\partial_y^2 b_1\|_{m-2}^2 ds\right)^\frac12\\
\lesssim&\;\varepsilon\mu\left[\|\partial_y^2v_1(0)\|_{[m/2]}
+\left(\int_0^t\|\partial_y^2v_1\|_{[m/2]+1}^2 ds\right)^\frac12\right] \left(\int_0^t\|(b_1,\partial_y  b_1)\|_{m-1}^2 ds\right)^\frac12\\
&\cdot \left(\int_0^t\|\partial_y^2 b_1\|_{m-2}^2 ds\right)^\frac12 +\varepsilon\mu\|  b_1\|_{[m/2], \infty}\left(\int_0^t\|\partial_y^2v_1\|_{m-2}^2 ds\right)^\frac12\left(\int_0^t\|\partial_y^2 b_1\|_{m-2}^2 ds\right)^\frac12.
\end{align*}
As for the last term on the right hand side of \eqref{4.38},
\begin{equation}\label{4.39}
\begin{split}
&2\varepsilon\mu\int_0^t\int_{\mathbb{R}_+^2}  \mathcal{Z}^\alpha(\tilde{b}_2\partial_y^3 v_1)\cdot\mathcal{Z}^\alpha\partial_y^2 b_1\ d{\bf{x}}ds\\
= &\;2\varepsilon\mu\sum_{\beta+\gamma=\alpha\atop |\beta|\geq 1}\int_0^t\int_{\mathbb{R}_+^2}  \mathcal{Z}^\beta \tilde{b}_2\mathcal{Z}^\gamma
\partial_y^3 v_1\cdot\mathcal{Z}^\alpha\partial_y^2 b_1\ d{\bf{x}}ds+2\varepsilon\mu\int_0^t\int_{\mathbb{R}_+^2}  \tilde{b}_2\mathcal{Z}^\alpha \partial_y^3 v_1\cdot\mathcal{Z}^\alpha\partial_y^2 b_1\ d{\bf{x}}ds.
\end{split}
\end{equation}
Then, the first term on the right hand side of \eqref{4.39} can be estimated as follows by using \eqref{2.3}, \eqref{2.4} and the Sobolev embedding inequality.

\begin{align*}
\begin{split}
&2\varepsilon\mu\sum_{\beta+\gamma=\alpha\atop |\beta|\geq 1}\int_0^t\int_{\mathbb{R}_+^2}  \mathcal{Z}^\beta \tilde{b}_2\mathcal{Z}^\gamma
\partial_y^3 v_1\cdot\mathcal{Z}^\alpha\partial_y^2 b_1\ d{\bf{x}}ds\\
\lesssim&\;\varepsilon\mu\sum_{\beta+\gamma=\alpha\atop |\beta|\geq|\gamma|}\sup_{0\leq s\leq t}\|\phi\mathcal{Z}^\gamma \partial_y^3 v_1(s)\|_{L_{x}^{\infty}(L_y^2)}\left(\int_0^t\|\phi^{-1}\mathcal{Z}^\beta\tilde b_2\|_{L_{x}^{2}(L_y^\infty)}^2 ds\right)^\frac12\left(\int_0^t\|\partial_y^2 b_1\|_{m-2}^2 ds\right)^\frac12\\
&+\varepsilon\mu\sum_{\beta+\gamma=\alpha\atop 1\leq |\beta|<|\gamma|<|\alpha|}\|\phi^{-1}\mathcal{Z}^\beta\tilde b_2\|_{L_{t, {\bf{x}}}^\infty}\left(\int_0^t\int_{\mathbb{R}_+^2}|\phi \mathcal{Z}^\gamma\partial_y^3 v_1|^2 d{\bf{x}}ds\right)^\frac12\left(\int_0^t\|\partial_y^2 b_1\|_{m-2}^2 ds\right)^\frac12\\
\lesssim &\;\varepsilon\mu\left[\|\partial_y^2 v_1(0)\|_{[m/2]+1}^2+\left(\int_0^t\|\partial_y^2 v_1\|_{[m/2]+2}^2 ds\right)^\frac12\right] \left(\int_0^t\|(b_1,\partial_y b_1)\|_{m-1}^2 ds\right)^\frac12\\
&\cdot \left(\int_0^t\|\partial_y^2 b_1\|_{m-2}^2 ds\right)^\frac12 +\varepsilon\mu\| b_1\|_{[m/2], \infty}\left(\int_0^t\|\partial_y^2 v_1\|_{m-2}^2 ds\right)^\frac12\left(\int_0^t\|\partial_y^2 b_1\|_{m-2}^2 ds\right)^\frac12.
\end{split}
\end{align*}
For the second term on the right hand side of \eqref{4.39}, we have
\begin{align*}
&2\varepsilon\mu\int_0^t\int_{\mathbb{R}_+^2}  \tilde{b}_2\mathcal{Z}^\alpha\partial_y^3 v_1\cdot\mathcal{Z}^\alpha\partial_y^2 b_1\ d{\bf{x}}ds\\
\leq&\;2\varepsilon\mu\|\phi^{-1}\tilde b_2\|_{L_{t, {\bf{x}}}^\infty}\left(\int_0^t\int_ {\mathbb{R}_+^2}|\mathcal{Z}^\alpha\partial_y^3 v_1|^2d{\bf{x}}ds\right)^\frac12\left(\int_0^t\int_{\mathbb{R}_+^2}|\phi\mathcal{Z}^\alpha\partial_y^2 b_1|^2 d{\bf{x}}ds\right)^\frac12\\
\leq&\;\frac{\varepsilon^2\mu^2}{4}\int_0^t\| \partial_y^3 v_1\|_{m-2}^2  ds+ 4\|b_1\|_{1, \infty}^2\int_0^t\|\partial_y b_1\|_{m-1}^2 ds.
\end{align*}
The third term on the right hand side of \eqref{4.37} is handled by similar arguments.
\begin{equation}\label{4.40}
\begin{split}
&-2\varepsilon\mu\int_0^t\int_{\mathbb{R}_+^2} \mathcal{Z}^\alpha\partial_y^2(v_2\partial_y b_1)\cdot \mathcal{Z}^\alpha\partial_y^2 b_1\ d{\bf{x}}ds\\
=&-2\varepsilon\mu\int_0^t\int_{\mathbb{R}_+^2} \mathcal{Z}^\alpha(\partial_y^2v_2\partial_y b_1+2\partial_yv_2\partial_y^2 b_1)\cdot \mathcal{Z}^\alpha\partial_y^2 b_1\ d{\bf{x}}ds\\
&-2\varepsilon\mu\int_0^t\int_{\mathbb{R}_+^2} \mathcal{Z}^\alpha( v_2\partial_y^3  b_1)\cdot \mathcal{Z}^\alpha\partial_y^2 b_1\ d{\bf{x}}ds.
\end{split}
\end{equation}
Here, by the Sobolev embedding inequality, the first term on the right hand side of \eqref{4.40} satisfies the following estimate.
\begin{align*}
&-2\varepsilon\mu\int_0^t\int_{\mathbb{R}_+^2} \mathcal{Z}^\alpha(\partial_y^2v_2\partial_y b_1+2\partial_yv_2\partial_y^2 b_1)\cdot \mathcal{Z}^\alpha\partial_y^2 b_1\ d{\bf{x}}ds\\
\lesssim&\;\varepsilon\mu\left[\sum_{\beta+\gamma=\alpha\atop |\beta|\leq|\gamma|}\sup_{0\leq s\leq t}\|\mathcal{Z}^\beta\partial_y^2(v_2,  b_1)(s)\|_{L_{x}^{\infty}(L_y^2)}\left(\int_0^t\|\mathcal{Z}^\gamma\partial_y(v_2,  b_1)\|_{L_{x}^{2}(L_y^\infty)}^2 ds\right)^\frac12\right.\\
&\left.+ \sum_{\beta+\gamma=\alpha\atop |\beta| >|\gamma|}\|\mathcal{Z}^\gamma\partial_y(v_2,  b_1)\|_{L_{t, {\bf{x}}}^\infty}\left(\int_0^t\int_ {\mathbb{R}_+^2}|\mathcal{Z}^\beta \partial_y^2(v_2,  b_1)|^2 d{\bf{x}}ds\right)^\frac12\right]\left(\int_0^t\|\partial_y^2 b_1\|_{m-2}^2 ds\right)^\frac12\\
\lesssim&\;\varepsilon\mu\left[\|\partial_y^2(v_2,  b_1)(0)\|_{[m/2]}+\left(\int_0^t\|\partial_y^2(v_2,  b_1)\|_{[m/2]+1}^2 ds\right)^\frac12\right]\\
&\cdot\sum_{i=1}^2\left(\int_0^t\|\partial_y^i(v_2, b_1)\|_{m-2}^2 ds\right)^\frac12 \left(\int_0^t\|\partial_y^2 b_1\|_{m-2}^2 ds\right)^\frac12 \\
 &+\varepsilon\mu\| \partial_y(v_2,  b_1)\|_{[m/2]-1, \infty}\left(\int_0^t\| \partial_y^2(v_2,  b_1)\|_{m-2}^2 ds\right)^\frac12\left(\int_0^t\|\partial_y^2 b_1\|_{m-2}^2 ds\right)^\frac12.
\end{align*}
And it is helpful to rewrite the second term on the right hand side of \eqref{4.40} into two parts.
\begin{equation}\label{4.41}
\begin{split}
-2\varepsilon\mu\int_0^t\int_{\mathbb{R}_+^2} \mathcal{Z}^\alpha( v_2\partial_y^3  b_1)\cdot \mathcal{Z}^\alpha\partial_y^2 b_1\ d{\bf{x}}ds
=&-2\varepsilon\mu\sum_{\beta+\gamma=\alpha\atop |\beta|\geq1}\int_0^t\int_{\mathbb{R}_+^2} \mathcal{Z}^\beta v_2 \mathcal{Z}^\gamma\partial_y^3  b_1\cdot \mathcal{Z}^\alpha\partial_y^2 b_1\ d{\bf{x}}ds\\
&-2\varepsilon\mu\int_0^t\int_{\mathbb{R}_+^2}  v_2\mathcal{Z}^\alpha\partial_y^3  b_1\cdot \mathcal{Z}^\alpha\partial_y^2 b_1\ d{\bf{x}}ds.
\end{split}
\end{equation}
Then, by the Sobolev embedding inequality, \eqref{2.3} and \eqref{2.4}, the first term on the right hand side of \eqref{4.41} is estimated.
\begin{align*}
&-2\varepsilon\mu\sum_{\beta+\gamma=\alpha\atop |\beta|\geq1}\int_0^t\int_{\mathbb{R}_+^2} \mathcal{Z}^\beta v_2 \mathcal{Z}^\gamma\partial_y^3  b_1\cdot \mathcal{Z}^\alpha\partial_y^2 b_1\ d{\bf{x}}ds\\
\lesssim&\;\varepsilon\mu\left[\sum_{\beta+\gamma=\alpha\atop |\beta|\geq|\gamma|}\sup_{0\leq s\leq t}\|\phi \mathcal{Z}^\gamma\partial_y^3b_1(s)\|_{L_{x}^{\infty}(L_y^2)}\left(\int_0^t\|\phi^{-1}\mathcal{Z}^\beta v_2\|_{L_{x}^{2}(L_y^\infty)}^2 ds\right)^\frac12\right.\\
&\left.+\sum_{\beta+\gamma=\alpha\atop |\beta|<|\gamma|<|\alpha|}\|\phi^{-1}\mathcal{Z}^\beta v_2\|_{L_{t, {\bf{x}}}^\infty}\left(\int_0^t\int_{\mathbb{R}_+^2}|\phi\mathcal{Z}^\gamma \partial_y^3 b_1|^2 d{\bf{x}}ds\right)^\frac12 \right]\left(\int_0^t\|\partial_y^2 b_1\|_{m-2}^2 ds\right)^\frac12\\
\lesssim&\;\varepsilon\mu\left[\|\partial_y^2 b_1(0)\|_{[m/2]+1}
+\left(\int_0^t\|\partial_y^2 b_1 \|_{[m/2]+2}^2 ds\right)^\frac12\right] \left(\int_0^t\|(\partial_y v_2, \partial_y^2 v_2)\|_{m-2}^2 ds\right)^\frac12\\
&\cdot\left(\int_0^t\|\partial_y^2 b_1\|_{m-2}^2 ds\right)^\frac12+\varepsilon\mu\|\partial_y v_2\|_{[m/2]-1, \infty}\int_0^t\|\partial_y^2 b_1\|_{m-2}^2 ds.
\end{align*}
For the second term on the right hand side of \eqref{4.41},  by \eqref{2.2} and integration by parts, we have

\begin{align*}
\begin{split}
&-2\varepsilon\mu\int_0^t\int_{\mathbb{R}_+^2}  v_2\mathcal{Z}^\alpha\partial_y^3  b_1\cdot \mathcal{Z}^\alpha\partial_y^2 b_1\ d{\bf{x}}ds\\
\lesssim&-2\varepsilon\mu\int_0^t\int_{\mathbb{R}_+^2} \left( v_2[\mathcal{Z}^\alpha,\partial_y]\partial_y^2  b_1+v_2\partial_y\mathcal{Z}^\alpha\partial_y^2b_1\right)\cdot \mathcal{Z}^\alpha\partial_y^2 b_1\ d{\bf{x}}ds\\
\lesssim&\;\varepsilon\mu\sum_{k=0}^{m-3}\|\phi^{-1}v_2\|_{L_{t, {\bf{x}}}^\infty}\left(\int_0^t\int_{\mathbb{R}_+^2}  |\phi\partial_y\mathcal{Z}_2^k\partial_y^2 b_1|^2 d{\bf{x}}ds\right)^\frac12\left(\int_0^t\|\partial_y^2 b_1\|_{m-2}^2 ds\right)^\frac12\\
&+\varepsilon\mu\|\partial_yv_2\|_{L_{t, {\bf{x}}}^\infty}\int_0^t\|\partial_y^2 b_1\|_{m-2}^2 ds\\
\lesssim&\;\varepsilon\mu\|\partial_yv_2\|_{L_{t, {\bf{x}}}^\infty}\int_0^t\|\partial_y^2 b_1\|_{m-2}^2 ds.
\end{split}
\end{align*}
Next, we deal with the fourth term on the right hand side of \eqref{4.37}.
\begin{equation}\label{4.42}
\begin{split}
&-2\varepsilon\mu\int_0^t\int_{\mathbb{R}_+^2} \mathcal{Z}^\alpha\partial_y^2(b_1\partial_y v_2)\cdot\mathcal{Z}^\alpha\partial_y^2 b_1\ d{\bf{x}}ds\\
= &-2\varepsilon\mu\int_0^t\int_{\mathbb{R}_+^2} \mathcal{Z}^\alpha(\partial_y^2b_1\partial_y v_2+2\partial_yb_1\partial_y^2 v_2)\cdot\mathcal{Z}^\alpha\partial_y^2 b_1\ d{\bf{x}}ds\\
&-2\varepsilon\mu\int_0^t\int_{\mathbb{R}_+^2} \mathcal{Z}^\alpha (b_1\partial_y^3 v_2)\cdot\mathcal{Z}^\alpha\partial_y^2 b_1\ d{\bf{x}}ds.
\end{split}
\end{equation}
For the first term on the right hand side of \eqref{4.42}, by the Sobolev embedding inequality,  we have
\begin{align*}
&-2\varepsilon\mu\int_0^t\int_{\mathbb{R}_+^2} \mathcal{Z}^\alpha(\partial_y^2b_1\partial_y v_2+2\partial_yb_1\partial_y^2 v_2)\cdot\mathcal{Z}^\alpha\partial_y^2 b_1\ d{\bf{x}}ds\\
\lesssim&\;\varepsilon\mu\left[\sum_{\beta+\gamma=\alpha\atop |\beta|<|\gamma|}\sup_{0\leq s\leq t}\|\mathcal{Z}^\beta \partial_y^2 (v_2,   b_1)(s)\|_{L_{x}^{\infty}(L_y^2)}\left(\int_0^t\|\mathcal{Z}^\gamma \partial_y( v_2,   b_1)\|_{L_{x}^{2}(L_y^\infty)}^2 ds\right)^\frac12\right. \\
&\left.+ \sum_{\beta+\gamma=\alpha\atop |\beta|\geq|\gamma|}\|\mathcal{Z}^\gamma \partial_y( v_2,  b_1)\|_{L_{t, {\bf{x}}}^\infty}\left(\int_0^t\int_{\mathbb{R}_+^2}|\mathcal{Z}^\beta\partial_y^2( v_2,   b_1)|^2 d{\bf{x}}ds\right)^\frac12\right]\left(\int_0^t\|\partial_y^2 b_1\|_{m-2}^2 ds\right)^\frac12\\
\lesssim&\;\varepsilon\mu\left[\|\partial_y^2( v_2,   b_1)(0)\|_{[m/2]}+\left(\int_0^t\|\partial_y^2( v_2,   b_1)\|_{[m/2]+1}^2 ds\right)^\frac12\right] \\
&\cdot\sum_{i=1}^2\left(\int_0^t\|\partial_y^i(  v_2, b_1 )\|_{m-2}^2 ds\right)^\frac12\left(\int_0^t\|\partial_y^2 b_1\|_{m-2}^2 ds\right)^\frac12\\
&+\varepsilon\mu\|\partial_y( v_2,   b_1)\|_{[m/2]-1, \infty}\left(\int_0^t\|\partial_y^2( v_2,  b_1)\|_{m-2}^2 ds\right)^\frac12\left(\int_0^t\|\partial_y^2 b_1\|_{m-2}^2 ds\right)^\frac12.
\end{align*}
As for the second term on the right hand side of \eqref{4.42}, one has
\begin{equation}\label{4.43}
\begin{split}
&-2\varepsilon\mu\int_0^t\int_{\mathbb{R}_+^2} \mathcal{Z}^\alpha (b_1\partial_y^3 v_2)\cdot\mathcal{Z}^\alpha\partial_y^2 b_1\ d{\bf{x}}ds\\
= &-2\varepsilon\mu\sum_{\beta+\gamma=\alpha\atop |\beta|\geq 1}\int_0^t\int_{\mathbb{R}_+^2} \mathcal{Z}^\beta b_1\mathcal{Z}^\gamma\partial_y^3 v_2\cdot\mathcal{Z}^\alpha\partial_y^2 b_1\ d{\bf{x}}ds\\
&-2\varepsilon\mu\int_0^t\int_{\mathbb{R}_+^2}  b_1\mathcal{Z}^\alpha\partial_y^3 v_2\cdot\mathcal{Z}^\alpha\partial_y^2 b_1\ d{\bf{x}}ds.
\end{split}
\end{equation}
Again, by \eqref{4.3}, the first term on the right hand side of \eqref{4.43} is rewritten as follows.
\begin{equation}\label{4.44}
\begin{split}
&-2\varepsilon\mu\sum_{\beta+\gamma=\alpha\atop |\beta|\geq 1}\int_0^t\int_{\mathbb{R}_+^2} \mathcal{Z}^\beta b_1\mathcal{Z}^\gamma\partial_y^3 v_2\cdot\mathcal{Z}^\alpha\partial_y^2 b_1\ d{\bf{x}}ds\\
=&\;2\varepsilon\mu\sum_{\beta+\gamma=\alpha\atop |\beta|\geq 1}\int_0^t\int_{\mathbb{R}_+^2} \mathcal{Z}^\beta b_1\mathcal{Z}^\gamma\partial_y^2\partial_x v_1\cdot\mathcal{Z}^\alpha\partial_y^2 b_1\ d{\bf{x}}ds\\
&+	2\gamma^{-1}\varepsilon\mu\sum_{\beta+\gamma=\alpha\atop |\beta|\geq 1}\int_0^t\int_{\mathbb{R}_+^2} \mathcal{Z}^\beta b_1\mathcal{Z}^\gamma\partial_y^2(p^{-1}\partial_t p) \cdot\mathcal{Z}^\alpha\partial_y^2 b_1\ d{\bf{x}}ds\\
&+2\gamma^{-1}\varepsilon\mu\sum_{\beta+\gamma=\alpha\atop |\beta|\geq 1}\int_0^t\int_{\mathbb{R}_+^2} \mathcal{Z}^\beta b_1\mathcal{Z}^\gamma\partial_y^2(p^{-1}{\bf{v}}\cdot\nabla p) \cdot\mathcal{Z}^\alpha\partial_y^2 b_1\ d{\bf{x}}ds.
\end{split}
\end{equation}
For the first term on the right hand side of \eqref{4.44}, by the Sobolev embedding inequality, one has
\begin{align*}
&2\varepsilon\mu\sum_{\beta+\gamma=\alpha\atop |\beta|\geq 1}\int_0^t\int_{\mathbb{R}_+^2} \mathcal{Z}^\beta b_1\mathcal{Z}^\gamma\partial_y^2\partial_x v_1\cdot\mathcal{Z}^\alpha\partial_y^2 b_1\ d{\bf{x}}ds\\
\lesssim&\;\varepsilon\mu\sum_{\beta+\gamma=\alpha\atop|\beta|\geq|\gamma|}\sup_{0\leq s\leq t}\|\mathcal{Z}^\gamma\partial_y^2\partial_x v_1(s)\|_{L_{x}^{\infty}(L_y^2)}\left(\int_0^t\|\mathcal{Z}^\beta b_1\|_{L_{x}^{2}(L_y^\infty)}^2 ds\right)^\frac12\left(\int_0^t\|\partial_y^2 b_1\|_{m-2}^2 ds\right)^\frac12\\
&+\varepsilon\mu\sum_{\beta+\gamma=\alpha\atop 1\leq|\beta|<\gamma|<|\alpha|} \|\mathcal{Z}^\beta b_1\|_{L_{t, {\bf{x}}}^\infty}\left(\int_0^t\int_{\mathbb{R}_+^2}|\mathcal{Z}^\gamma\partial_y^2\partial_x v_1|^2 d{\bf{x}}ds\right)^\frac12\left(\int_0^t\|\partial_y^2 b_1\|_{m-2}^2 ds\right)^\frac12\\
\lesssim&\;\varepsilon\mu\left[\|\partial_y^2 v_1(0)\|_{[m/2]+1}+\left(\int_0^t\|\partial_y^2 v_1\|_{[m/2]+2}^2 ds\right)^\frac12\right] \left(\int_0^t\|(b_1,\partial_y b_1)\|_{m-2}^2 ds\right)^\frac12\\
&\cdot \left(\int_0^t\| \partial_y^2 b_1\|_{m-2}^2 ds\right)^\frac12 +\varepsilon\mu\|b_1\|_{[m/2]-1, \infty}\left(\int_0^t\|\partial_y^2 v_1\|_{m-2}^2 ds\right)^\frac12\left(\int_0^t\|\partial_y^2 b_1\|_{m-2}^2 ds\right)^\frac12.
\end{align*}
And the second term on the right hand side of \eqref{4.44} includes three parts.
\begin{equation}\label{4.45}
\begin{split}
&2\gamma^{-1}\varepsilon\mu\sum_{\beta+\gamma=\alpha\atop |\beta|\geq 1}\int_0^t\int_{\mathbb{R}_+^2} \mathcal{Z}^\beta b_1\mathcal{Z}^\gamma\partial_y^2(p^{-1}\partial_t p) \cdot\mathcal{Z}^\alpha\partial_y^2 b_1\ d{\bf{x}}ds\\
=&2\gamma^{-1}\varepsilon\mu\sum_{\beta+\gamma=\alpha\atop |\beta|\geq 1}\sum_{\gamma_1+\gamma_2=\gamma}\int_0^t\int_{\mathbb{R}_+^2} \mathcal{Z}^\beta b_1\mathcal{Z}^{\gamma_1}\partial_y^2p^{-1}\mathcal{Z}^{\gamma_2}\partial_t p \cdot\mathcal{Z}^\alpha\partial_y^2 b_1\ d{\bf{x}}ds\\
&+4\gamma^{-1}\varepsilon\mu\sum_{\beta+\gamma=\alpha\atop |\beta|\geq 1}\sum_{\gamma_1+\gamma_2=\gamma}\int_0^t\int_{\mathbb{R}_+^2} \mathcal{Z}^\beta b_1\mathcal{Z}^{\gamma_1}\partial_yp^{-1}\mathcal{Z}^{\gamma_2}\partial_t \partial_yp \cdot\mathcal{Z}^\alpha\partial_y^2 b_1\ d{\bf{x}}ds\\
&+2\gamma^{-1}\varepsilon\mu\sum_{\beta+\gamma=\alpha\atop |\beta|\geq 1}\sum_{\gamma_1+\gamma_2=\gamma}\int_0^t\int_{\mathbb{R}_+^2} \mathcal{Z}^\beta b_1\mathcal{Z}^{\gamma_1} p^{-1}\mathcal{Z}^{\gamma_2}\partial_t \partial_y^2p \cdot\mathcal{Z}^\alpha\partial_y^2 b_1\ d{\bf{x}}ds.
\end{split}
\end{equation}
For the first term on the right hand side of \eqref{4.45}, it follows from the Sobolev embedding inequality that
\begin{align*}
&2\gamma^{-1}\varepsilon\mu\sum_{\beta+\gamma=\alpha\atop |\beta|\geq 1}\sum_{\gamma_1+\gamma_2=\gamma}\int_0^t\int_{\mathbb{R}_+^2} \mathcal{Z}^\beta b_1\mathcal{Z}^{\gamma_1}\partial_y^2p^{-1}\mathcal{Z}^{\gamma_2}\partial_t p \cdot\mathcal{Z}^\alpha\partial_y^2 b_1\ d{\bf{x}}ds\\
\lesssim&\;\varepsilon\mu\left[\sum_{\beta+\gamma_1+\gamma_2=\alpha\atop |\beta|\geq |\gamma_1|, |\gamma_2|}\sup_{0\leq s\leq t}\|\mathcal{Z}^{\gamma_1}\partial_y^2p^{-1}(s)\|_{L_{x}^{\infty}(L_y^2)}\|\mathcal{Z}^{\gamma_2}\partial_t p\|_{L_{t, {\bf{x}}}^\infty}\left(\int_0^t\|\mathcal{Z}^\beta b_1\|_{L_{x}^{2}(L_y^\infty)}^2 ds\right)^\frac12\right.\\
&\left.+\sum_{\beta+\gamma_1+\gamma_2=\alpha\atop |\gamma_1|\geq  |\beta|,|\gamma_2|}\|\mathcal{Z}^\beta b_1\|_{L_{t, {\bf{x}}}^\infty}\|\mathcal{Z}^{\gamma_2}\partial_t p\|_{L_{t, {\bf{x}}}^\infty}\left(\int_0^t\int_{\mathbb{R}_+^2}|\mathcal{Z}^{\gamma_1}\partial_y^2p^{-1}|^2d{\bf{x}}ds\right)^\frac12\right.\\
&\left.+\sum_{\beta+\gamma_1+\gamma_2=\alpha\atop |\gamma_2|\geq|\beta|, |\gamma_1|}\sup_{0\leq s\leq t}\|\mathcal{Z}^{\gamma_1}\partial_y^2p^{-1}(s)\|_{L_{x}^{\infty}(L_y^2)}\|\mathcal{Z}^\beta b_1\|_{L_{t, {\bf{x}}}^\infty}\left(\int_0^t\|\mathcal{Z}^{\gamma_2}\partial_t p\|_{L_{x}^{2}(L_y^\infty)}^2 ds\right)^\frac12
\right]\\
&\cdot\left(\int_0^t\|\partial_y^2 b_1\|_{m-2}^2 ds\right)^\frac12\\
\lesssim&\;\varepsilon\mu\|(b_1, p)\|_{[m/2], \infty}\left[\|\partial_y^2 p^{-1}(0)\|_{[m/2]}+\left(\int_0^t\|\partial_y^2 p^{-1} \|_{[m/2]+1}^2 ds\right)^\frac12\right]\\
&\cdot\sum_{i=0}^1\left(\int_0^t\|\partial_y^i(b_1, p)\|_{m-2}^2 ds\right)^\frac12\left(\int_0^t\|\partial_y^2 b_1\|_{m-2}^2 ds\right)^\frac12\\
&+\varepsilon\mu\|b_1\|_{[m/2]-1, \infty}\|p\|_{[m/2], \infty}\left(\int_0^t\|\partial_y^2 p^{-1}\|_{m-2}^2 ds\right)^\frac12\left(\int_0^t\|\partial_y^2 b_1\|_{m-2}^2 ds\right)^\frac12.
\end{align*}
For the second term on the right hand side of \eqref{4.45}, we have
\begin{align*}
&4\gamma^{-1}\varepsilon\mu\sum_{\beta+\gamma=\alpha\atop |\beta|\geq 1}\sum_{\gamma_1+\gamma_2=\gamma}\int_0^t\int_{\mathbb{R}_+^2} \mathcal{Z}^\beta b_1\mathcal{Z}^{\gamma_1}\partial_yp^{-1}\mathcal{Z}^{\gamma_2}\partial_t \partial_yp \cdot\mathcal{Z}^\alpha\partial_y^2 b_1\ d{\bf{x}}ds\\
\lesssim&\; \varepsilon\mu\left[\sum_{\beta+\gamma_1+\gamma_2=\alpha\atop |\beta|\geq|\gamma_1|, |\gamma_2|}\|\mathcal{Z}^{\gamma_1}\partial_yp^{-1}\|_{L_{t, {\bf{x}}}^\infty}\|\mathcal{Z}^{\gamma_2}\partial_t \partial_yp\|_{L_{t, {\bf{x}}}^\infty}\left(\int_0^t\int_{\mathbb{R}_+^2}|\mathcal{Z}^\beta b_1|^2 d{\bf{x}}ds\right)^\frac12 \right.\\
&\left.+  \sum_{\beta+\gamma_1+\gamma_2=\alpha\atop |\gamma_1|\geq |\beta|, |\gamma_2|}\|\mathcal{Z}^\beta b_1\|_{L_{t, {\bf{x}}}^\infty}\|\mathcal{Z}^{\gamma_2}\partial_t \partial_yp\|_{L_{t, {\bf{x}}}^\infty}\left(\int_0^t\int_{\mathbb{R}_+^2}|\mathcal{Z}^{\gamma_1}\partial_yp^{-1}|^2 d{\bf{x}}ds\right)^\frac12\right. \\
&\left.+  \sum_{\beta+\gamma_1+\gamma_2=\alpha\atop |\gamma_2|\geq |\beta|, |\gamma_1|}\|\mathcal{Z}^\beta b_1\|_{L_{t, {\bf{x}}}^\infty}\|\mathcal{Z}^{\gamma_1}\partial_yp^{-1}\|_{L_{t, {\bf{x}}}^\infty}\left(\int_0^t\int_{\mathbb{R}_+^2}|\mathcal{Z}^{\gamma_2}\partial_t \partial_yp|^2 d{\bf{x}}ds\right)^\frac12\right]\cdot\left(\int_0^t\|\partial_y^2 b_1\|_{m-2}^2 ds\right)^\frac12\\
\lesssim&\;\varepsilon\mu  \|(b_1, \partial_y p, \partial_y p^{-1} )\|_{[m/2],\infty}^2\left(\int_0^t\|(b_1,\partial_y p, \partial_y p^{-1} )\|_{m-1}^2 ds\right)^\frac12\left(\int_0^t\|\partial_y^2 b_1\|_{m-2}^2 ds\right)^\frac12.
\end{align*}
For the last term on the right hand side of \eqref{4.45} satisfies the similar estimate as the first term.

By a similar argument as for the second term on the right hand side of \eqref{4.44}, the last term on the right hand side of \eqref{4.44} is estimated as follows.
\begin{align*}
&2\gamma^{-1}
\varepsilon\mu\sum_{\beta+\gamma=\alpha\atop |\beta|\geq 1}\int_0^t\int_{\mathbb{R}_+^2} \mathcal{Z}^\beta b_1\mathcal{Z}^\gamma\partial_y^2(p^{-1}{\bf{v}}\cdot\nabla p) \cdot\mathcal{Z}^\alpha\partial_y^2 b_1\ d{\bf{x}}ds\\
\lesssim&\;\varepsilon\mu\sum_{i=0}^1\|\partial_y^i({\bf{v}},b_1,  p, p^{-1})\|_{[m/2], \infty}^3\sum_{i=0}^2\left(\int_0^t\|\partial_y^i({\bf{v}}, p-1, p^{-1}-1)\|_{m-2}^2 ds\right)^\frac12\\
&\cdot\left(\int_0^t\|\partial_y^2 b_1\|_{m-2}^2 ds\right)
^\frac12 +\varepsilon\mu\sum_{i=0}^1\|\partial_y^i({\bf{v}},b_1,  p, p^{-1})\|_{[m/2], \infty}^2\left[\|\partial_y^2({\bf{v}}, p, p^{-1})(0)\|_{[m/2]+1}\right.\\
&\left.+\left(\int_0^t\|\partial_y^2({\bf{v}},p-1, p^{-1}-1)
\|_{[m/2]+2}^2ds\right)^\frac12\right]\cdot\sum_{j=0}^1\left(\int_0^t\|\partial_y^j({\bf{v}}, b_1, p-1, p^{-1}-1)\|_{m-2}^2 ds\right)^\frac12\\
&\cdot\left(\int_0^t\|\partial_y^2 b_1\|_{m-2}^2 ds\right)^\frac12.
\end{align*}
By using Young's inequality, the second term on the right hand side of \eqref{4.43} is estimated.
\begin{align*}
-2\varepsilon\mu\int_0^t\int_{\mathbb{R}_+^2}  b_1\mathcal{Z}^\alpha\partial_y^3 v_2\cdot\mathcal{Z}^\alpha\partial_y^2 b_1\ d{\bf{x}}ds\leq 2\varepsilon^2\mu^2\|b_1\|_{L_{t, {\bf{x}}}^\infty}^2 \int_0^t  \| \partial_y^3 v_2\|_{m-2}^2 ds+\frac12 \int_0^t\|\partial_y^2 b_1\|_{m-2}^2 ds.
\end{align*}
For the last term on the right hand side of \eqref{4.37}, one obtains

\begin{equation}\label{4.46}
\begin{split}
-2\varepsilon\mu\int_0^t\int_{\mathbb{R}_+^2} \mathcal{Z}^\alpha\partial_y^2(v_1\partial_x b_1)\cdot\mathcal{Z}^\alpha\partial_y^2 b_1\ d{\bf{x}}ds
=&-2\varepsilon\mu\int_0^t\int_{\mathbb{R}_+^2} \mathcal{Z}^\alpha(\partial_y^2v_1\partial_x b_1)\cdot\mathcal{Z}^\alpha\partial_y^2 b_1\ d{\bf{x}}ds\\
&-4\varepsilon\mu\int_0^t\int_{\mathbb{R}_+^2} \mathcal{Z}^\alpha(\partial_yv_1\partial_y\partial_x b_1)\cdot\mathcal{Z}^\alpha\partial_y^2 b_1\ d{\bf{x}}ds\\
&-2\varepsilon\mu\int_0^t\int_{\mathbb{R}_+^2} \mathcal{Z}^\alpha(v_1\partial_y^2\partial_x b_1)\cdot\mathcal{Z}^\alpha\partial_y^2 b_1\ d{\bf{x}}ds.
\end{split}
\end{equation}
The first term on the right hand side of \eqref{4.46} is estimated by the Sobolev embedding inequality.
\begin{align*}
&-2\varepsilon\mu\int_0^t\int_{\mathbb{R}_+^2} \mathcal{Z}^\alpha(\partial_y^2v_1\partial_x b_1)\cdot\mathcal{Z}^\alpha\partial_y^2 b_1\ d{\bf{x}}ds\\
\lesssim&\;\varepsilon\mu\sum_{\beta+\gamma=\alpha\atop |\beta|\geq|\gamma|}\|\mathcal{Z}^\gamma \partial_x b_1\|_{L_{t, {\bf{x}}}^\infty}\left(\int_0^t\int_{\mathbb{R}_+^2}|\mathcal{Z}^\beta\partial_y^2 v_1|^2 d{\bf{x}}ds\right)^\frac12\left(\int_0^t\|\partial_y^2 b_1\|_{m-2}^2 ds\right)^\frac12\\
&+\varepsilon\mu\sum_{\beta+\gamma=\alpha\atop |\beta|<|\gamma|}\sup_{0\leq s\leq t}\|\mathcal{Z}^\beta \partial_y^2 v_1(s)\|_{L_{x}^{\infty}(L_y^2)}\left(\int_0^t\|\mathcal{Z}^\gamma \partial_x b_1\|_{L_{x}^{2}(L_y^\infty)}^2 ds\right)^\frac12\cdot\left(\int_0^t\|\partial_y^2 b_1\|_{m-2}^2 ds\right)^\frac12\\
\lesssim&\;\varepsilon\mu\|b_1\|_{[m/2], \infty}\left(\int_0^t\|\partial_y^2 v_1\|_{m-2}^2 ds\right)^\frac12\left(\int_0^t\|\partial_y^2 b_1\|_{m-2}^2 ds\right)^\frac12\\
&+\varepsilon\mu\left[\|\partial_y^2 v_1(0)\|_{[m/2]}+\left(\int_0^t\|\partial_y^2 v_1\|_{[m/2]+1}^2 ds\right)^\frac12\right] \\
&\cdot\left(\int_0^t\| (b_1,\partial_y b_1)\|_{m-1}^2 ds\right)^\frac12 \left(\int_0^t\|\partial_y^2 b_1\|_{m-2}^2 ds\right)^\frac12.
\end{align*}
By \eqref{2.1}, the second term on the right hand side of \eqref{4.46} satisfies the estimate as follows.
\begin{align*}
&-4\varepsilon\mu\int_0^t\int_{\mathbb{R}_+^2} \mathcal{Z}^\alpha(\partial_yv_1\partial_y\partial_x b_1)\cdot\mathcal{Z}^\alpha\partial_y^2 b_1\ d{\bf{x}}ds\\
\lesssim&\;\varepsilon\mu\|\partial_y v_1\|_{L_{t, {\bf{x}}}^\infty}\left(\int_0^t\|\partial_y b_1\|_{m-1}^2 ds\right)^\frac12
\left(\int_0^t\|\partial_y^2 b_1\|_{m-2}^2 ds\right)^\frac12\\
&+\varepsilon\mu\|\partial_y b_1\|_{1,\infty}\left(\int_0^t\|\partial_y v_1\|_{m-2}^2 ds\right)^\frac12
\left(\int_0^t\|\partial_y^2 b_1\|_{m-2}^2 ds\right)^\frac12.
\end{align*}
For the last term on the right hand side of \eqref{4.46}, by integration by parts, \eqref{2.3} and \eqref{2.4}, we have

\begin{align*}
\begin{split}
&-2\varepsilon\mu\int_0^t\int_{\mathbb{R}_+^2} \mathcal{Z}^\alpha(v_1\partial_y^2\partial_x b_1)\cdot\mathcal{Z}^\alpha\partial_y^2 b_1\ d{\bf{x}}ds\\
\lesssim&\;\varepsilon\mu\sum_{\beta+\gamma=\alpha\atop |\beta|\geq|\gamma|} \|\phi\mathcal{Z}^\gamma\partial_y^2 \partial_x b_1\|_{L_{t, {\bf{x}}}^{\infty }}\left(\int_0^t\int_{\mathbb{R}_+^2}|\phi^{-1}\mathcal{Z}^\beta v_1|^2 d{\bf{x}}ds\right)^\frac12\left(\int_0^t\|\partial_y^2 b_1\|_{m-2}^2 ds\right)^\frac12\\
&+\varepsilon\mu\sum_{\beta+\gamma=\alpha\atop |\beta|<|\gamma|<|\alpha|}\|\phi^{-1}\mathcal{Z}^\beta v_1\|_{L_{t, {\bf{x}}}^\infty}\left(\int_0^t\int_{\mathbb{R}_+^2}|\phi\mathcal{Z}^\gamma \partial_y^2 \partial_x b_1|^2 d{\bf{x}}ds\right)^\frac12\left(\int_0^t\|\partial_y^2 b_1\|_{m-2}^2 ds\right)^\frac12]\\
&+\varepsilon\mu\|\phi^{-1}\partial_x v_1\|_{L_{t, {\bf{x}}}^\infty}\left(\int_0^t\int_{\mathbb{R}_+^2}|\phi\mathcal{Z}^\alpha\partial_y^2 b_1|^2 d{\bf{x}}ds\right)^\frac12\left(\int_0^t\| \partial_y^2 b_1\|_{m-2}^2  ds\right)^\frac12\\
\lesssim&\;\varepsilon\mu\|\partial_y b_1\|_{[m/2]+1, \infty}\left(\int_0^t\|\partial_y v_1\|_{m-2}^2 ds\right)^\frac12 \left(\int_0^t\|\partial_y^2 b_1\|_{m-2}^2 ds\right)^\frac12\\
&+\varepsilon\mu\|\partial_y v_1\|_{[m/2], \infty}\left(\int_0^t\|\partial_y b_1\|_{m-1}^2 ds\right)^\frac12 \left(\int_0^t\|\partial_y^2 b_1\|_{m-2}^2 ds\right)^\frac12.
\end{split}
\end{align*}

Next, it is left to consider the terms on the right hand side of \eqref{4.36}. For the first two terms on the right hand side of \eqref{4.36}, by \eqref{2.1}, the following estimate holds true.
\begin{align*}
&\int_0^t\int_{\mathbb{R}_+^2} |\mathcal{Z}^\alpha \partial_y(\rho\partial_t v_1)|^2 d{\bf{x}}ds+\int_0^t\int_{\mathbb{R}_+^2} |\mathcal{Z}^\alpha \partial_y(b_2\partial_x b_2)|^2 d{\bf{x}}ds \\
\lesssim&\int_0^t\int_{\mathbb{R}_+^2} |\mathcal{Z}^\alpha (\partial_y\rho\partial_t v_1)|^2 d{\bf{x}}ds+\int_0^t\int_{\mathbb{R}_+^2} |\mathcal{Z}^\alpha (\rho\partial_t\partial_y v_1)|^2 d{\bf{x}}ds\\
&+\int_0^t\int_{\mathbb{R}_+^2} |\mathcal{Z}^\alpha( \partial_yb_2\partial_x b_2)|^2 d{\bf{x}}ds+\int_0^t\int_{\mathbb{R}_+^2} |\mathcal{Z}^\alpha(b_2\partial_x \partial_y b_2)|^2 d{\bf{x}}ds\\
\lesssim&\sum_{i, j=0}^1\|\partial_y^i(\rho, v_1, b_2)\|_{1, \infty}^2\int_0^t\|\partial_y^j(\rho-1, v_1, \tilde{b}_2)\|_{m-1}^2 ds.
\end{align*}
For the third term on the right hand side of \eqref{4.36}, we have
\begin{equation}\label{4.47}
\begin{split}
&\int_0^t\int_{\mathbb{R}_+^2} |\mathcal{Z}^\alpha\partial_y (\rho {\bf{v}}\cdot\nabla v_1)|^2 d{\bf{x}}ds\\
\lesssim&\int_0^t\int_{\mathbb{R}_+^2} |\mathcal{Z}^\alpha\partial_y (\rho v_1\partial_x v_1)|^2 d{\bf{x}}ds +\int_0^t\int_{\mathbb{R}_+^2} |\mathcal{Z}^\alpha\partial_y (\rho v_2\partial_y v_1)|^2 d{\bf{x}}ds.
\end{split}
\end{equation}
Then, the first term on the right hand side of \eqref{4.47} is estimated by \eqref{2.1}.

\begin{equation*}
\begin{split}
&\int_0^t\int_{\mathbb{R}_+^2} |\mathcal{Z}^\alpha\partial_y (\rho v_1\partial_x v_1)|^2 d{\bf{x}}ds\\
\lesssim& \int_0^t\int_{\mathbb{R}_+^2} |\mathcal{Z}^\alpha(\partial_y \rho v_1\partial_x v_1)|^2 d{\bf{x}}ds+\int_0^t\int_{\mathbb{R}_+^2} |\mathcal{Z}^\alpha(\rho \partial_y v_1\partial_x v_1)|^2 d{\bf{x}}ds\\
&+\int_0^t\int_{\mathbb{R}_+^2} |\mathcal{Z}^\alpha(\rho v_1\partial_x \partial_y v_1)|^2 d{\bf{x}}ds\\
\lesssim&\sum_{i, j=0}^1\|\partial_y^i(\rho, v_1)\|_{1, \infty}^4\int_0^t\|\partial_y^j(\rho-1, v_1 )\|_{m-1}^2 ds.
\end{split}
\end{equation*}
Similarly, for the second term on the right hand side of \eqref{4.47}, we obtain
\begin{align*}
&\int_0^t\int_{\mathbb{R}_+^2} |\mathcal{Z}^\alpha\partial_y (\rho v_2\partial_y v_1)|^2 d{\bf{x}}ds\\
\lesssim&\int_0^t\int_{\mathbb{R}_+^2} |\mathcal{Z}^\alpha (\partial_y\rho v_2\partial_y v_1)|^2 d{\bf{x}}ds+\int_0^t\int_{\mathbb{R}_+^2} |\mathcal{Z}^\alpha (\rho \partial_yv_2\partial_y v_1)|^2 d{\bf{x}}ds\\
&+\int_0^t\int_{\mathbb{R}_+^2} |\mathcal{Z}^\alpha (\rho \phi^{-1}v_2\phi\partial_y^2 v_1)|^2 d{\bf{x}}ds\\
\lesssim&\sum_{i, j=0}^1\|\partial_y^i(\rho, {\bf{v}})\|_{1, \infty}^4\int_0^t\|\partial_y^j(\rho-1, {\bf{v}} )\|_{m-1}^2 ds.
\end{align*}
For the last third term on the right hand side of \eqref{4.36}, since $\nabla\cdot {\bf{B}}=0$,  one derives
\begin{align*}
&\int_0^t\int_{\mathbb{R}_+^2} |\mathcal{Z}^\alpha\partial_y (\tilde b_2\partial_y b_1)|^2 d{\bf{x}}ds\\
\lesssim&\int_0^t\int_{\mathbb{R}_+^2} |\mathcal{Z}^\alpha(\partial_y \tilde b_2\partial_y b_1)|^2 d{\bf{x}}ds+\int_0^t\int_{\mathbb{R}_+^2} |\mathcal{Z}^\alpha  (\tilde b_2\partial_y^2 b_1)|^2 d{\bf{x}}ds\\
\lesssim&\;\|\partial_y \tilde b_2\|_{L_{t, {\bf{x}}}^\infty}^2\int_0^t\|\partial_y b_1\|_{m-2}^2 ds+\|\partial_y b_1\|_{L_{t, {\bf{x}}}^\infty}^2\int_0^t\|\partial_y \tilde b_2\|_{m-2}^2 ds\\
&+\|\phi^{-1}  \tilde b_2\|_{L_{t, {\bf{x}}}^\infty}^2\int_0^t\|\phi\partial_y^2 b_1\|_{m-2}^2 ds+\|\phi\partial_y ^2b_1\|_{L_{t, {\bf{x}}}^\infty}^2\int_0^t\| \phi^{-1} \tilde b_2\|_{m-2}^2 ds\\
\lesssim&\;\| b_1\|_{1,\infty}^2\int_0^t\|\partial_y b_1\|_{m-1}^2 ds+\|\partial_y b_1\|_{ 1, \infty}^2\int_0^t\| b_1\|_{m-1}^2 ds.
\end{align*}
Finally, it is direct to estimate the last two terms on the right hand side of \eqref{4.36}.
\begin{align*}
&\varepsilon^2(2\mu+\lambda)^2\int_0^t\int_{\mathbb{R}_+^2}|\mathcal{Z}^\alpha \partial_x^2\partial_y v_1|^2 d{\bf{x}}ds+\varepsilon^2(\mu+\lambda)^2\int_0^t\int_{\mathbb{R}_+^2} |\mathcal{Z}^\alpha \partial_x\partial_y^2v_2|^2 d{\bf{x}}ds\\
\lesssim&\;\varepsilon^2(2\mu+\lambda)^2\int_0^t\|\partial_y v_1\|_{m}^2 ds+\varepsilon^2(\mu+\lambda)^2\int_0^t\|\partial_y^2 v_2\|_{m-1}^2 ds.
\end{align*}
Combining all of estimates in this subsection, we have
\begin{align*}
&\int_0^t \left(\| \partial_y^2 b_1\|_{m-2}^2  + \varepsilon^2\mu^2 \|  \partial_y^3 v_1\|_{m-2}^2 \right)ds+\varepsilon\mu\| \partial_y^2 b_1(t)\|_{m-2}^2 \\
\lesssim&\;\varepsilon\mu\| \partial_y^2 b_1(0)\|_{m-2}^2+\varepsilon^2\mu^2\|b_1\|_{L_{t,{\bf{x}}}^\infty}^2
	\int_0^t\|  \partial_y^3 v_2\|_{m-2}^2 ds+\varepsilon^2(2\mu+\lambda)^2\int_0^t\left(\|\partial_y v_1\|_{m}^2 +\|\partial_y^2v_2\|_{m-1}^2 \right) ds  \\ &+\varepsilon\mu\left[\left(1+\sum_{i=0}^1\|\partial_y^i({\bf{v}}, b_1, p, p^{-1})\|_{[m/2]+1, \infty}^2\right)^2 +  \sum_{i=1}^2\|\partial_y^i({\bf{v}}, b_1, p, p^{-1})(0)\|_{[m/2]+1}^2\right.\\
&\left.+\sum_{i=1}^2\int_0^t\|\partial_y^i({\bf{v}}, b_1, p-1, p^{-1}-1)\|_{[m/2]+2}^2 ds \right]\cdot\sum_{j=0}^2 \int_0^t \|\partial_y^j({\bf{v}}, b_1, p-1, p^{-1}-1)\|_{m-j}^2 ds \\
 &+\left(1+\sum_{i=0}^1\|\partial_y^i(\rho, {\bf{v}},   {\bf{B}})\|_{1, \infty}^2\right)^2\sum_{j=0}^1\int_0^t\|\partial_y^j(\rho-1, {\bf{v}},   {\bf{B}}-\overset{\rightarrow}{e_y} )\|_{m-1}^2 ds.
\end{align*}

As a consequence, the proof of Lemma \ref{lem3} is complete by collecting the conormal  estimates of $\partial_y^2{\bf{v}}, \partial_y^2 b_1$ and $\partial_y^2 p$ together. It is noted that the smallness of $\|b_1\|_{L_{t,{\bf{x}}}^\infty}$ is also needed here to close the energy estimates. And this a priori assumption can be verified in Section 5.

\section{The Proof of Theorem \ref{Th1}}
In this section, we show the proof of Theorem \ref{Th1}. Based on the estimates in Lemmas \ref{lem1}, \ref{lem2} and \ref{lem3}, we have
\begin{equation}\label{5.1}
\begin{split}
& N_m(t)+ \sum_{|\alpha|+i\leq  m\atop i=1,2}\int_\Omega
\left(\gamma^{-1}\varepsilon(2\mu+\lambda) p^{-1}(t)|\mathcal{Z}^\alpha \partial_y^i p(t)|^2 d{\bf{x}}+\varepsilon\mu |\mathcal{Z}^\alpha\partial_y^i b_1(t)|^2 \right)d{\bf{x}}\\
\lesssim&\;N_m(0)+ \sum_{|\alpha|+i\leq m\atop i=1,2}\int_\Omega
\left(\gamma^{-1}\varepsilon(2\mu+\lambda) p_0^{-1}|\mathcal{Z}^\alpha \partial_y^i p_0|^2 d{\bf{x}}+\varepsilon\mu |\mathcal{Z}^\alpha\partial_y^i b_1(0)|^2 \right)d{\bf{x}}\\
&+\left[\left(1+\sum_{i=0}^1\|\partial_y^i(\rho, {\bf{v}},  {\bf{B}}, p, p^{-1} )\|_{[m/2]+1, \infty}^2\right)^2+\|\partial_y^2 b_1(0)\|_{[m/2]}^2+\int_0^t\|\partial_y^2 b_1\|_{[m/2]+1}^2 ds\right]\\
&\cdot\sum_{j=0}^1\int_0^t \|\partial_y^j(\rho-1, {\bf{v}},  {\bf{B}}-\overset{\rightarrow}{e_y}, p-1, p^{-1}-1)\|_{m-j}^2 ds\\
&+\varepsilon(2\mu+\lambda) \Bigg[\left(1+\sum_{i=0}^1\|\partial_y^i({\bf{v}}, b_1, p, p^{-1})\|_{[m/2]+1, \infty}\right)^2+\sum_{i=1}^2\|\partial_y^i({\bf{v}}, b_1, p, p^{-1})(0)\|_{[m/2]+2}^2\\
&+\int_0^t\|\partial_y^2({\bf{v}}, b_1, p-1, p^{-1}-1)(0)\|_{[m/2]+3}^2ds\Bigg]\cdot\sum_{j=0}^2\int_0^t  \|\partial_y^j({\bf{v}}, b_1, p-1, p^{-1}-1)\|_{m-j}^2 ds.
\end{split}
\end{equation}
By Lemma \ref{lem 2.2}, the third term on the right hand side of \eqref{5.1} can be estimated by the following way.
\begin{align*}
&\sum_{i=0}^1\|\partial_y^i(\rho, {\bf{v}},  {\bf{B}}, p, p^{-1} )\|_{[m/2]+1, \infty}\notag\\
\lesssim&\sum_{i=0}^1\|\partial_y^i(\rho, {\bf{v}},  {\bf{B}}, p, p^{-1} )(0)\|_{[m/2]+3}+	\sum_{i =1}^2\|\partial_y^i(\rho, {\bf{v}},  {\bf{B}}, p, p^{-1}) (0)\|_{[m/2]+2}\\
&+\int_0^t\left(	\sum_{i =0}^1\|\partial_y^i(\rho, {\bf{v}},  {\bf{B}}, p-1, p^{-1}-1)\|_{[m/2]+4}+	\sum_{i=1}^2\|\partial_y^i(\rho, {\bf{v}},  {\bf{B}}, p, p^{-1} )\|_{[m/2]+3}\right) ds\\
\lesssim&\;\mathcal{P}(N_m(0))+t\mathcal{P}(N_m(t)).
\end{align*}

For any $m\geq 9$,  by inserting the above inequality into \eqref{5.1},  we have
\begin{equation}\label{5.2}
\begin{split}
& N_m(t)+ \sum_{|\alpha|+i\leq m\atop i=1,2}\int_{\mathbb{R}_+^2}
\left(\gamma^{-1}\varepsilon(2\mu+\lambda) p^{-1}(t)|\mathcal{Z}^\alpha \partial_y^i p(t)|^2 d{\bf{x}}+\varepsilon\mu |\mathcal{Z}^\alpha\partial_y^i b_1(t)|^2 \right)d{\bf{x}}\\
 \lesssim&\;\;\mathcal{P}(N_m(0))+[t+\varepsilon(2\mu+\lambda)]\mathcal{P}(N_m(t)).
\end{split}
\end{equation}

Let  the time $t$ and $\varepsilon$ be suitablely small,  then we achieve that
\begin{align*}
N_m(t)+ \sum_{|\alpha|+i\leq m\atop i=1,2}\int_{\mathbb{R}_+^2}
\left(\gamma^{-1}\varepsilon(2\mu+\lambda) p^{-1}(t)|\mathcal{Z}^\alpha \partial_y^i p(t)|^2 d{\bf{x}}+\varepsilon\mu |\mathcal{Z}^\alpha\partial_y^i b_1|^2 \right)d{\bf{x}}
\lesssim\mathcal{P}(N_m(0)).
\end{align*}
Based on the above uniform conormal estimates, we can justify the inviscid limit by the similar arguments in \cite{MR12}. Then,
the proof of Theorem \ref{Th1} is complete.

\section*{Acknowledgement}
The research  of X. Cui was  partially supported by China Postdoctoral Science Foundation grant 2021M692088.
F.  Xie's research was supported by National Natural Science Foundation of China No.11831003, Shanghai Science and Technology Innovation Action Plan No. 20JC1413000 and Institute of Modern Analysis-A Frontier Research Center of Shanghai.

\end{document}